THE UNIVERSITY OF CHICAGO

# ON BLOCH-KATO CONJECTURE OF TAMAGAWA NUMBERS FOR HECKE CHARACTERS OF IMAGINARY QUADRATIC NUMBER FIELD

A DISSERTATION SUBMITTED TO

THE FACULTY OF THE DIVISION OF THE PHYSICAL SCIENCES

IN CANDIDACY FOR THE DEGREE OF

DOCTOR OF PHILOSOPHY

DEPARTMENT OF MATHEMATICS

BY

BING HAN

# ACKNOWLEDGMENTS

I want to first thank my advisor Professor Spencer Bloch for his kind guidance and constant encouragement throughout my graduate study. I benefit enourmously from his deep insight in many areas of mathematics. He has set a great example for me. I especially thank him for introducing me to the beautiful field of Arithmetic Algebraic Geometry, and all his effort to patiently read through earlier drafts of this thesis and many excellent suggestions. I also want to thank him for his confidence in me and his understanding.

I also want to dedicate this thesis to my family, who have all shown me great love and support during all these years. I especially thank my parents for their great sacrifice. I can not thank my dear wife Zhe Yang enough for all her love.

I am fortunate to study at University of Chciago, which supplies me with the best environment one can hope for to do research in mathematics. Many people contribute to my growth. I want to thank heartly Professor Niels Nygaard and Professor Madhav Nori for all their kind help and mathematical discussions. It was Niels who got me interested in number theory in the first place and supervised my first topic. Thanks are also due to many great staffs at our Math department.

Above all, I want to thank God for his blessing on me and my family. Thank Him for leading me to the results I presented in this thesis.

# TABLE OF CONTENTS

## ABSTRACT

In this thesis I prove the validity of Tamagawa number conjecture of Bloch-Kato for certain Hecke characters. I study the exponential map and local Tamagawa measure for all odd primes. I also study $p$-part of Shafarevich-Tate group for motives associated to Hecke characters, for all prime $p \neq 2$ or 3.

Here is the main result of this paper: (for more detail, see Chapter 2.2)

Let $K$ be an imaginary quadratic number field of class number 1 and discriminant $-d_K$, where $d_K > 0$. Fix $\psi$ of type $(1,0)$ and conductor $\mathfrak{f}$ which satisfies $\psi(\bar{\boldsymbol{\alpha}}) = \overline{\psi(\boldsymbol{\alpha})}$ for all ideals $\boldsymbol{\alpha}$ of $\mathcal{O}_K$. By the theory of complex multiplication there exists an elliptic curve $E$ over $\mathbf{Q}$ such that $End(E) = \mathcal{O}_K$ and the Grössencharacter associated to $E$ in the sense of Deuring is precisely $\psi$. For each prime $\mathfrak{p}$ of $\mathcal{O}_K$, let $\psi_{\mathfrak{p}} : Gal(K(\mathfrak{f}\mathfrak{p}^\infty)/K) \to K_{\mathfrak{p}}^*$ be the Weil realization of $\psi$ at $\mathfrak{p}$.

Fix integer $k, j$ such that $k - j > 1$ and $j \geq 0$. There exists motive $\mathcal{M}_{k,j}$ such that $L(\mathcal{M}_{k,j}, 0) = L(\overline{\psi^{j+k}}, k)$. (see [5] prop 2.1). Let $\Delta = Gal(K(E_{\mathfrak{p}})/K)$ which is a finite group of order prime to $p$ since by theory of complex multiplication we have $\Delta \hookrightarrow (\mathcal{O}_{\mathfrak{p}}/\mathfrak{p})^*$. Hence we can view $\Delta$ also as a subgroup of $Gal(K(\mathfrak{f}\mathfrak{p}^\infty)/K)$ and consider the restriction of $\psi_{\mathfrak{p}}$ on $\Delta$.

**Main Theorem.** *Fix integer $k - j > 1$ and $j \geq 0$. Then the $p-$ part of Tamagawa number conjecture for $\mathcal{M}_{k,j}$ is true for all prime $p \neq 2, 3$ which splits in $K/\mathbb{Q}$ and all $p$ which is nonsplit in $K/\mathbb{Q}$ such that $\psi_p^k \bar{\psi}_p^{-j}$ is nontrivial when restricted to $\Delta$.*

# CHAPTER 1

# INTRODUCTION

Zeta and $L$-functions are important and intriguing subjects that have been the focal point of study in number theory ever since the beginning. Classically there are Dedekind zeta functions and Dirichlet $L$-functions. Dirichlet applied his $L$-functions to the proof of existence of infinite primes in arithmetic progessions. Kummer showed the connection between Fermat's Last Theorem and class numbers of cyclotomic fields, which are closely related to special values of $L$-functions by class number formula. In the begining of this century Hecke $L$-functions and Artin $L$-functions were introduced. Hecke $L$-functions are special cases of automorphic $L$-functions. These are very general class of $L$-functions which have nice analytic properties such as analytic continuation and functional equations. On the other hand, through the advancement of algebraic geometry, especially cohomology theories and theories of motives, one can define Hasse-Weil $L$-functions for general motive which include Artin $L$-functions as special case. These $L$-functions encode an amazing amount of arithmetic information about the motive. Langland's conjecture essentially identifies these two kinds of $L$-functions.

Modern interest in the study of special values of $L$-functions at integer argument is revived by Birch-Swinnerton-Dyer conjecture, which is an analogue of the class number formula for abelian varieties:

**Conjecture 1 (Brich-Swinnerton-Dyer).** *Let $E$ be an elliptic curve over $\mathbb{Q}$. Then $L(E/\mathbb{Q}, s)$ has a zero at $s = 1$ of order equal to the rank $r$ of Mordell-Weil group $E/\mathbb{Q}$ and*

$$\lim_{s\to 1}(s-1)^{-r}L(E/\mathbb{Q}, s) = \#(\bigsqcup(E))R(E/\mathbb{Q})\Omega \prod_{p<\infty} c_p \#((E/\mathbb{Q})_{tors})^{-2}$$

*where* $\Omega = \int_{E(\mathbb{R})} |\omega|$ *is the real period of Neron differential* $\omega$ *of* $E$*,* $R(E/\mathbb{Q})$ *is a regulator term comming from the height pairing, and* $\bigsqcup(E)$ *is Shaferivich-Tate group which is conjectured to be finite torsion.*

Bloch had given a volume theoretic version of the above conjecture as $G = E$ (an elliptic curve) in the following generalization to the case $G$ being an extension of $E$ by a torus.

**Conjecture 2 (Tamagwa-Bloch).** *Let* $G$ *be an semi-simple algebraic group over* $\mathbb{Q}$*.*

$$\tau(G) = \frac{\#(Pic(G)_{tors})}{\#(\bigsqcup(G))}$$

where $\tau(G) = \int_{G(\mathbf{A}_{\mathbb{Q}})/G(\mathbb{Q})}(\omega, (\lambda_p))$ is Tamagawa number of G with respect to Tamagawa measure $(\omega, (\lambda_p))$. Here $\omega$ is a left invariant gauge form of $dim(G)$ and $(\lambda_p)$ is suitable convergence factor.

Deligne formulated conjecture on special values of $L$-function at critical point and its relation to period integrals, generalizing classical work of Euler among others. Inspired by work of Bloch, Beilinson viewed $K$ groups as rational motivic cohomology group and defined regulator map which represents chern character to real Deligne-Beilinson cohomology groups. He made several conjectures on special values of $L$-functions which generalized Deligne's conjecture. For a motive $\mathfrak{M}$ with coefficent in a number field $E$, these conjectures give the order of vanishing of $L(\mathfrak{M}^*(1), s)$ at $s = 0$. They predict interesting elements in $K$ groups whose image under the regulator map determine the first nonzero coefficent of Taylor series expansion of $L(\mathfrak{M}, s)$ around $s = 0$ up to $E^*$. In order to remove this ambiguity, Bloch and Kato made the Tamagawa number conjecture for pure motives with coefficents in $\mathbb{Q}$ as follows ([1] 5.15):

**Conjecture 3 (Bloch-Kato).** *Assume the triple* $(V, D, \phi)$ *comes from a motif. Let* $\mathfrak{M}$ *be a* $\mathbb{Z}$*-lattice in* $V$ *such that* $\mathfrak{M} \otimes \hat{\mathbb{Z}}$ *is Galois stable in* $V \otimes \mathbf{A}_f$*. Then* $\bigsqcup(\mathfrak{M})$ *is finite and*

$$Tam(\mathfrak{M}) = \frac{\#(H^0(\mathbb{Q}, \mathfrak{M}^* \otimes \mathbb{Q}/\mathbb{Z}(1)))}{\#(\bigsqcup(\mathfrak{M}))}$$

Here $(V, D)$ is a motivic pair (see [1] 5.15 for its definition) which is a pair of $\mathbb{Q}$ vector space with additional structures and comparision isomorphisms between them when tensored with $\mathbb{C}$ or $\mathbf{Q}_p$, satisfying some additional axioms. They usually arise from the rational Betti and De Rham cohomology of pure motives over $\mathbb{Q}$. $Tam(\mathfrak{M})$ and $\bigsqcup(\mathfrak{M})$ are defined analogously as for $G$ above. See section 2.1 for more information on their definition and $\phi$ (which is a certain Haar measure). The main difficulty is to give a reasonable definition of local points $\mathbb{A}(\mathbb{Q}_p)$ for $p \leqslant \infty$, global points $\mathbb{A}(\mathbb{Q})$ and Haar measures on $\mathbb{A}(\mathbb{Q}_p)$. For this purpose Bloch and Kato rely heavily on $p$-adic Hodge theory, especially Fontaine's $B_{\text{cris}}$ and $B_{\text{dR}}$. See section 2.1 for a brief recall of this conjecture of Bloch-Kato.

Over the past few years there have been several important new development. First Fontaine and Perrin-Riou ([4] ) reformulated Bloch-Kato conjecture in terms of measure on Fundamental line induced by canonical measure on Euler-Poincare line which enable them to extend the Tamagawa number conjecture for motives over $F$ with coefficent in $E$ (here $E$ and $F$ are number fields). It also reduces the verification of the conjecture for general motives to pure motives.

Secondly Kato ([8]) made a generalized Iwasawa Main conjecture which generalize the classical Iwasawa main conjecture and Bloch-Kato conjecture at the same time. Most importantly he conjectured the existence of $p$-adic zeta elements (non-archimedean analogues of zeta values in $\mathbb{C}$) which can be viewed as "$p$-adic Langlands program".

Thirdly, Perrin-Riou ([9]) has proposed a new framework for defining general $p$-adic $L$-functions for crystalline representation using Bloch-Kato's exponential map and Iwasawa theory. Moreover she established a direct link between Tamagawa numbers of motives with her conjectual construction of arithmetic $p$-adic $L$-function. We expect Kato's $p$-adic zeta elements which can be made into an abstract Euler system also give rise to Perrin-Riou's $p$-adic $L$-function. From their work, we see in order to prove Tamagawa number conjecture, we should look for magical elements in $K$ groups (resp. Galois cohomolgy groups) which are strongly related to special values of $L$-function via Beilinson's regulator map (resp. exponential or dual exponential map). Then try to show that these magical elements form integral basis of some canonically defined modules which typically requires an Iwasawa type argument.

Finally history seems to come back to a full circle when in Wiles' recent proof of Fermat's Last theorem, he reduced the problem to an estimation of the order of certain Selmer group which in turn follows from $p$ part of Tamagwa number conjecture for motive $sym^2(E)$ for a semi-stable modular elliptic curve $E$ at $p = 3$ or $p = 5$ !

As for the verification of the Tamagawa number conjecture, besides the complete proof of it in the case of Riemann Zeta function (i.e for motive $\mathbb{Q}(r)$) and partial results for motive $h^1(E)(2)$ for CM elliptic curve $E$ by Bloch-Kato, the only progress has been M.Harrison's ([6]) proof of $p$ part of Tamagawa number conjecture for certain special Hecke character of Gaussian field $\mathbb{Q}(\mathbf{i})$ in the case $p$ splits in $\mathbb{Q}(\mathbf{i})$.

The goal of this paper is to extend Harrison's result to more general Hecke characters of more general fields, verifying $p$ part of Tamagawa number conjecture for $p$ which is nonsplit as well as for $p$ which splits. Our method is basically different from Harrison's and works equally well for all imaginary quadratic fields. We recover his result as a special case of ours. For simplicity of presentation, we will just consider the case for $cl(K) = 1$. In general the only change is notationwise. By the same method we can check $p$ part of Tamagawa number conjecture for all Hecke characters of $K$ of type $(k, -j)$ with $0 \leq j$ and $k - j > 1$, for all odd prime $p \neq 3$ not dividing the number of roots of unity in Hilbert class field of $K$ (when $p$ is nonsplit, we also need to assume that the Hecke character is nontrivial when restricted to $\Delta$, see our Main theorem in 2.2.) The last restriction comes from Iwasawas main conjecture. The results we used of Kato (i.e explicit reciprocity law), Rubin (i.e Iwasawa main conjecture of quadratic number field) and results on arithmetic of Hecke characters and Eisenstein series all extends to general $K$, so are our extensions of some of these results we need for the proof.

Brieftly, here are the contents of this paper:

- In chapter 2, we will recall some background material on Tamagawa number conjecture and fix some notations in the special case of Hecke characters we will consider. We do some preliminary reduction of the Tamagawa number formula to be proved.

- In chapter 3 we write the local Tamagawa numbers as a product of 2 terms: one term will appear again in the formula for $p$ part of Shafarevich group; another term involves evaluation of dual exponential map on special element (in Galois cohomlogy group) constructed from elliptic units. We review dual exponential map and define these special elements. We also study carefully deep results of Kato on explicit reciprocity law for 1-dimesional Lubin-Tate formal groups which will be used for the calculation in chapter 4 and chapter 5.

- In chapter 4 we define carefully the elliptic units we are going to use and calculate the Coleman power series associated to it. The logarithmic derivative of this power series is essentially Eisenstein series. We review the properties of Eisenstein series and its link to special values of Hecke-$L$ function. We evaluate the dual exponential map in the case $j = 0$. We distinguish two cases according to whether $p$ is a good reduction prime or not. When $p$ is a bad reduction prime, since $E$ has CM, we go up a finite field extension and pass on to a twist of $E$ which has good reduction at $p$. We study the behavior of dual exponential map under field extension and twist.

- Chapter 5 is the most technical part of this thesis. Here we are trying to get information on the value of dual exponentail map evaluated on the special element in the case $j > 0$ from Kato's explicit reciprocity law (which require $j = 0$). Roughly we go up along a tower of fields (essentially unramified when $p$ splits but totally ramified when $p$ is nonsplit, which make it necessary to consider these 2 cases separately) and "untwist" our character. Interestingly the $p$-adic period of $E$ play a major role in these translation from $j > 0$ case to $j = 0$ case. In particular we need a congruence property between $p$-adic period and special values of weight 1 Eisenstein series evaluated on $p^n$ torsion points (which actually holds for general modular elliptic curves!). On the $L$-function side, we deduce realtion between $L(\bar{\psi}^{k+j}, k)$ and $L(\bar{\psi}^k, k)$ via their relation to values of Eisenstein series. It is known in the theory of elliptic functions how to relate $E_{j,k}$ with $E_k$. We arrive at our result by a congruence argument.

- In chapter 6 we deal with the Shafarevich-Tate group associated to our motive. Again we need to consider the case $p$ splits and the case $p$ is nonsplit separately. We need first reduce the Shafarevich group to more familiar objects studied in Iwasawa theory. Here we use crucially the Iwasawa Main conjecture proved by Rubin (both 1 and 2 variable case).

# CHAPTER 2

# MOTIVE FOR HECKE CHARACTER; PRELIMINARY RESULTS

## 2.1 Tamagawa number conjecture

Here we recall Bloch-Kato's conjecture in more detail.

**Definition 1.** Let $K$ be a finite extension of $\mathbf{Q}_p$. Let $V$ be a $\mathbf{Q}_l$ vector space endowed with a continous $Gal(\bar{K}/K)$ action. If $p = l$ we assume moreover that $V$ is a de Rham representation. Let $T$ be a Galois stable lattice of $V$. Then

$$H^1_f(K,V) \stackrel{\text{def}}{=} \begin{cases} Ker(H^1(K,V) \longrightarrow H^1(K^{nr},V)), & \text{if } l \neq p; \\ Ker(H^1(K,V) \longrightarrow H^1(K, B_{\text{crys}} \otimes V)), & \text{if } l = p; \end{cases}$$

and we define $H^1_f(K,T)$ as the full iverse image of $H^1_f(K,V)$ under the natural map $H^1(K,T) \to H^1(K,V)$. In particular $H^1_f(K,T)$ containes all the torsion of $H^1(K,T)$.

*Notation.* For a de Rham representation $V$ of $Gal(\bar{K}/K)$, denote $DR(V)$ to be $H^0(K, V \otimes B_{\text{DR}})$ and $Fil^i DR(V) = DR(V)^i$ to be $H^0(K, V \otimes Fil^i B_{\text{DR}})$.

We define the exponential map

$$DR(V)/DR(V)^0 \to H^1_f(K,V)$$

to be the connecting homomorphism from the exact squence of Galois modules obtained by tensoring $V$ with Fontaine's fundamental exact squence:

$$0 \to \mathbf{Q}_p \stackrel{\alpha}{\to} B_{\text{crys}} \oplus B^+_{\text{DR}} \stackrel{\beta}{\to} B_{\text{crys}} \oplus B_{\text{DR}} \to 0$$

where $\alpha(x) = (x,x)$ and $\beta(x,y) = (x - f(x), x - y)$, $f$ being the usual frobenuis.

Now define local and global points for a motivic pair $(V,D)$ of weights $\leqslant -1$. These are finite dimesional $\mathbf{Q}$ vector spaces with extra structures, satisfying some axioms. For details see [1] 5.5. One can just think of them as realizations of motives, with various comparison isomorphism comming from ($p$-adic) Hodge theory, $V$ given by Betti cohomology and $D$ coming from de Rham cohomology. Fix a $\mathbf{Z}$ lattice $M$ in $V$ such that $M \otimes \hat{\mathbf{Z}}$ is Galois stable.

**Definition 2.**

$$\mathbf{A}(\mathbf{Q}_p) \stackrel{\text{def}}{=} \begin{cases} H^1_f(\mathbf{Q}_p, M \otimes \hat{\mathbf{Z}}), & \text{if } p < \infty \\ ((D_\infty \otimes_{\mathbf{R}} \mathbf{C})/(D^0_\infty \otimes_{\mathbf{R}} \mathbf{C}) + M)^+, & \text{if } p = \infty \end{cases}$$

We regard $\mathbf{A}(\mathbf{Q}_p)$ for $p < \infty$ as a compact group with natural topology, and $\mathbf{A}(\mathbf{R})$ as a locally compact group. Fix an isomorphism

$$\omega : det_{\mathbf{Q}}(D/D^0) \cong \mathbf{Q}$$

which induces for each $p \leqslant \infty$

$$det_{\mathbf{Q}_p}(D_p/D^0_p) \cong \mathbf{Q}_p$$

This defines a measure on each $D_p/D^0_p$ and hence a Haar measure on $\mathbf{A}(\mathbf{Q}_p)$ via exponential map.

Bloch and Kato proved (see Theorem 4.1 of [1]) that there is a finite set of bad primes $S$ which consists of places at infinity, bad reduction primes and prime $p$ less than the length of filtration of $DR(V)$ and for $p \notin S$

$$\mu_{p,\omega}(\mathbf{A}(\mathbf{Q}_p)) = P_p(V, 1)$$

where for $l \neq p$, $P_p(V,1) \stackrel{\text{def}}{=} det(1 - frob_p : V_l^{I_p})$ is the $p$-Euler factor of $L(V,0)$. Now assume the weights are $\leqslant -3$, then the product

$$L_S(V,0)^{-1} = \prod_{p \notin S} \mu_{p,\omega}(\mathbf{A}(\mathbf{Q}_p))$$

converges. Bloch and Kato defined the Tamagawa measure $\mu$ for the motivic pair $(V, D)$ as $\mu = \prod_{p \leqslant \infty} \mu_{p,\omega}$ on $\prod_{p \leqslant \infty} \mathbf{A}(\mathbf{Q}_p)$. It is clear $\mu$ is independent of choice of $\omega$.

$\mathbf{A}(\mathbf{Q})$ is a finitely generated abelian group (generalized global Selmer group) such that

$$\mathbf{A}(\mathbf{Q}) \otimes \hat{\mathbf{Z}} = H^1_{f,Spec(\mathbf{Z})}(\mathbf{Q}, M \otimes \hat{\mathbf{Z}})$$

where $H^1_{f,Spec(\mathbf{Z})}$ means the set of elements in the global Galois cohomology group whose image under localization map at each prime $p$ lies in the $H^1_f$ part. There are natural homomorphism $\mathbf{A}(\mathbf{Q}) \to \mathbf{A}(\mathbf{Q}_p)$ for $p < \infty$, as well as $\mathbf{A}(\mathbf{Q}) \to \mathbf{A}(\mathbf{R})/\mathbf{A}(\mathbf{R})_{\text{cpt}} = D_\infty/(D^0_\infty + V^+_\infty)$. Now define

$$Tam(M) = \mu((\prod \mathbf{A}(\mathbf{Q}_p))/\mathbf{A}(\mathbf{Q}))$$

and define

$$\bigsqcup(M) \overset{\text{def}}{=} Ker\left(\frac{H^1(\mathbf{Q}, M \otimes \mathbf{Q}/\mathbf{Z})}{\mathbf{A}(\mathbf{Q}) \otimes \mathbf{Q}/\mathbf{Z}} \to \bigoplus_{p \leqslant \infty} \frac{H^1(\mathbf{Q}_p, M \otimes \mathbf{Q}/\mathbf{Z})}{\mathbf{A}(\mathbf{Q}_p) \otimes \mathbf{Q}/\mathbf{Z}}\right)$$

One knows that $p$ part of $\bigsqcup(M)$ is finite for each finite prime $p$, even though in general it is very hard to verify $\bigsqcup(M)$ itself is finite. Tamagawa number conjecture can also be written as

$$L_S(V,0) = \frac{\#(\bigsqcup(M))}{\#(H^0(\mathbf{Q}, M^* \otimes \mathbf{Q}/\mathbf{Z}(1)))}\mu_{\infty,\omega}(\mathbf{A}(\mathbf{R})/\mathbf{A}(\mathbf{Q})) \prod_{p\in S-\infty} \mu_{p,\omega}(\mathbf{A}(\mathbf{Q}_p)) \quad (2.1.1)$$

## 2.2 Hecke characters and associated motives

Let $K$ be an imaginary quadratic number field of class number 1 and discriminant $-d_K$, where $d_K > 0$. Let $\phi$ be a Grössencharacter of type $A_0$. We say $\phi$ has conductor $\mathfrak{f}$ and infinity type $(r,s)$if $\phi((\alpha)) = \alpha^r\bar{\alpha}^s$ whenever $\alpha \equiv 1 \pmod{\mathfrak{f}}$ , $\alpha \in K^*$ and $\mathfrak{f}$ is the smallest ideal of $\mathcal{O}_K$ such a relation holds. Fix f such that $(\mathrm{f}) = \mathfrak{f}$. Fix $\psi$ of type $(1,0)$ and conductor $\mathfrak{f}$ which satisfies $\psi(\bar{\boldsymbol{\alpha}}) = \overline{\psi(\boldsymbol{\alpha})}$ for all ideals $\boldsymbol{\alpha}$ of $\mathcal{O}_K$. By the theory of complex multiplication there exists an elliptic curve $E$ over $\mathbf{Q}$ such that $End(E) = \mathcal{O}_K$ and the Grössencharacter associated to $E$ in the sense of Deuring is precisely $\psi$. For each prime $\mathfrak{p}$ of $\mathcal{O}_K$, let $\psi_{\mathfrak{p}}$ be the Weil realization of $\psi$ at $\mathfrak{p}$.

$$\begin{array}{ccccc} \mathbf{A}_K^*/K^* & \xrightarrow{\psi} & K^* & \longrightarrow & K_{\mathfrak{p}}^* \\ & \searrow^{\text{artin map}} & & \nearrow_{\psi_{\mathfrak{p}}} & \\ & & Gal(K(\mathfrak{f}\mathfrak{p}^\infty)/K) & & \end{array}$$

*Notation.* We view all our global fields as subfields of $\bar{\mathbf{Q}}$ and choose once for all a place $\upsilon_p$ of $\bar{\mathbf{Q}}$ above $p$, and an embedding $i_p : \bar{\mathbf{Q}} \to \bar{\mathbf{Q}}_p$. We identify decomposion group $D_p$ (resp. inertia group $I_p$) at $\upsilon_p$ with $Gal(\bar{\mathbf{Q}_p}/\mathbf{Q}_p)$ (resp. the inertia subgroup of $Gal(\bar{\mathbf{Q}_p}/\mathbf{Q}_p)$) via $i_p$. When $p = \mathfrak{p}\mathfrak{p}^*$ we pick $\mathfrak{p}$ so that $\upsilon_p \mid \mathfrak{p}$. we identify $K_{\mathfrak{p}}$ with $\mathbf{Q}_p$ by $i_{\mathfrak{p}}$ which is the restriction of $i_p$ to K. we identify $K_{\mathfrak{p}^*}$ with $\mathbf{Q}_p$ by $i_{\mathfrak{p}} \circ \tau$ where $\tau$ is the complex conjugation. In particular, we have $\psi_{\mathfrak{p}} \mid_{I_p} = \chi_p$ the cyclotomic character and $\psi_{\mathfrak{p}^*} \mid_{I_p} = 1$. Let $K(\mathfrak{f})$ denote the ray class field of conductor $\mathfrak{f}$ and for any prime ideal $\mathfrak{p} \subseteq \mathcal{O}_K$ set $K(\mathfrak{f}\mathfrak{p}^\infty)$ as the union of all $K(\mathfrak{f}\mathfrak{p}^n)$. It is well known that $K(\mathfrak{f}\mathfrak{p}^n) = K(E_{\mathfrak{f}\mathfrak{p}^n})$. Fix a global minimal Weierstrass model for $E$ and fix $\omega$ as the Neron differential. Let $L$ be the period lattice and fix $\Omega_\infty \in \mathbb{R}$ an $\mathcal{O}_K$ generator of $L$: $\Omega_\infty = \int_{E(\mathbb{R})} \omega$. (sometimes called the real period of $E$). Let $\rho : \mathbf{C}/L \to E(\mathbf{C})$ be given by $\rho(z) = (\mathcal{P}(z), \mathcal{P}'(z))$ where $\mathcal{P}$ is the Weierstrass $\mathcal{P}$ function.

Let $\varphi = \psi^k \overline{\psi}^{-j}$. Define imprimitive $L$-series

$$L(\varphi^{-1}, s) = \sum_{\substack{\boldsymbol{\alpha} \subseteq \mathcal{O}_K \\ (\boldsymbol{\alpha}, \mathfrak{f}) = 1}} \varphi^{-1}(\boldsymbol{\alpha}) \mathbf{N}\boldsymbol{\alpha}^{-s}$$

The analytic continuation and functional equation for $L(\varphi^{-1}, s)$ are known. Since $\psi\bar{\psi} = \chi \overset{\text{def}}{=} \chi_{\text{cyclo}}$, $L(\bar{\psi}^j \psi^{-k}, 0) = L(\overline{\psi^{j+k}}, k)$. It is easy to see that if $k - j > 1$, then $L(\overline{\psi^{j+k}}, k) \neq 0$. A theorem of Damerell says that $(\frac{2\pi}{\sqrt{d_K}})^j \Omega_\infty^{-(k+j)} L(\overline{\psi^{j+k}}, k)$ is algebraic, belongs to $K$ if $0 \leqslant j < k$. From now on assume $0 \leqslant j < k$. Moreover since $\psi(\boldsymbol{\alpha}) = \overline{\psi(\boldsymbol{\alpha})}$, it follows that $L(\overline{\psi^{j+k}}, k) \in \mathbf{R}$. Hence in our case $(\frac{2\pi}{\sqrt{d_K}})^j \Omega_\infty^{-(k+j)} L(\overline{\psi^{j+k}}, k) \in \mathbf{Q}^*$.

On the other hand it is known how to associate motives to general Hecke character using abelian varieties with complex multiplication. For our $\varphi$ there exists a motive $\mathcal{M}_{k,j}$ such that $L(\mathcal{M}_{k,j},0) = L(\overline{\psi^{j+k}},k)$. (see [5] prop 2.1). Let $\Delta = Gal(K(E_{\mathfrak{p}})/K)$ which is a finite group of order prime to $p$ since by theory of complex multiplication we have $\Delta \hookrightarrow (\mathcal{O}_{\mathfrak{p}}/\mathfrak{p})^*$. Hence we can view $\Delta$ also as a subgroup of $Gal(K(\mathfrak{f}\mathfrak{p}^\infty)/K)$ and consider the restriction of $\psi_{\mathfrak{p}}$ on $\Delta$.

**Main Theorem.** *Fix integer $k,j$ such that $k-j>1$ and $j \geq 0$. Then the $p$ part of Tamagawa number conjecture for $\mathcal{M}_{k,j}$ is true for all prime $p \neq 2,3$ which splits in $K/\mathbb{Q}$ and all $p$ which is nonsplit in $K/\mathbb{Q}$ such that $\psi_p^k \bar{\psi}_p^{-j}$ is nontrivial when restricted to $\Delta$.*

$\mathcal{M}_{k,j}$ comes from $\underbrace{E \times \cdots \times E}_{k+j}$ twisted by $-j$ times the Tate motive, i.e it is a factor of $(h^1(E)(1))^{\otimes(k+j)}(-j) = h_1(E)^{\otimes(k+j)}(-j)$. Notice $wt(\mathcal{M}_{k,j}) = -k+j < -1$. $\mathcal{M}_{k,j}$ have been constructed in [6] section 2 ($j=0$) and [5] section 1 (for $j>0$). Harrison essentially proved the special case of our main theorem when $p$ splits in $K/\mathbf{Q}$ for $K=\mathbb{Q}(i)$ and $j=0$. Guo proved the special case of main theorem for $p$ good ordinary prime sufficently big which is unramified in $K/\mathbb{Q}$ (under this restriction he does not need to calculate local Tamagawa numbers). Notice both of these authors rely on (1 or 2 variable) $p$-adic $L$-function for the Hecke character of $K$ when $p$ splits, whose counterpart is still lacking when $p$ is inert. We get around this and our proof is different from theirs even in the case that $p$ splits.

Now brieftly recall results and notations for $\mathcal{M}_{k,j}$. For motive $h^1(E)(1)$, its realization is given by

- $M_B = H_B^{-1}(E) = \mathbf{Q}\Gamma_1 \oplus \mathbf{Q}\Gamma_2$, where $\Gamma_1$ is the image under $\rho$ of the line segement from the origin to $\Omega_\infty$ along the real axis and $\Gamma_2$ is the image under $\rho$ of the line segement from the origin to $\sqrt{d}\Omega_\infty$ along the imaginary axis.

- $M_{\mathrm{dR}} = H_{\mathrm{dR}}^{-1}(E) = \mathbf{Q}\hat{\omega} \oplus \mathbf{Q}\hat{\eta}$ where $\hat{\omega}, \hat{\eta}$ are dual basis of basis $\omega = dx/y$ and $\eta = xdx/y$. Since $Fil^1 H_{\mathrm{dR}}^1(E) = \mathbf{Q}\omega$ we see that

$$M_{\mathrm{dR}} = Fil^{-1} M_{\mathrm{dR}} \supseteq Fil^0 M_{\mathrm{dR}} = \mathbf{Q}\hat{\eta} \supseteq Fil^1 M_{\mathrm{dR}} = 0$$

- $$M_p = H_p^{-1}(E) = \begin{cases} (T_{\mathfrak{p}}(E) \oplus T_{\mathfrak{p}^*}(E)) \otimes \mathbf{Q}_p = \mathbf{Q}_p(\psi_{\mathfrak{p}}) \oplus \mathbf{Q}_p(\psi_{\mathfrak{p}}^*) & p \text{ splits} \\ T_{\mathfrak{p}}(E) \otimes \mathbf{Q}_p = K_{\mathfrak{p}}(\psi_{\mathfrak{p}}) & p \text{ is nonsplit} \end{cases} \tag{2.2.1}$$

  where $\tau$ acts by swapping the two factors when $p$ splits and by its natural action on $K_{\mathfrak{p}}$ otherwise.

Realizations for motive $\mathbb{G}_m$ (or $\mathbb{Q}(1)$) are

- $M_B = H_B^{-1}(\mathbb{G}_m) = \mathbf{Q}\Gamma$ where $\Gamma$ is the projection of the directed path from 0 to $2\pi i$ under $\mathbf{C} \to \mathbf{C}/2\pi i\mathbf{Z}$
- $M_{\mathrm{dR}} = H_{\mathrm{dR}}^{-1}(\mathbb{G}_m) = \mathbf{Q}\hat{\varepsilon}$ where $\hat{\varepsilon}$ is the dual basis of $\varepsilon = dx/x \in H_{\mathrm{dR}}^1(\mathbb{G}_m)$.
  $Fil^0 M_{\mathrm{dR}} = M_{\mathrm{dR}} \supseteq Fil^1 M_{\mathrm{dR}} = 0$
- $M_p = H_p^{-1}(\mathbb{G}_m) = \mathbf{Q}_p(\chi)$

It follows that the realizations for our motive $\mathcal{M}_{k,j}$ are given by

- $V \overset{\mathrm{def}}{=} M_B = (\mathbf{Q}e_1 \oplus \mathbf{Q}e_2) \otimes \mathbf{Q}\Gamma^{\otimes(-j)}$ where $e_1$ and $e_2$ are tensors in $\Gamma_1$ and $\Gamma_2$ with rational coefficents defined by

  $$d^{-\frac{k+j}{2}} \sum_{\Lambda \subseteq \{1,\dots,k+j\}} d^{\#(\frac{\Lambda}{2})} \Gamma_1^{\sigma} \otimes \cdots \otimes \Gamma_{k+j}^{\sigma} = e_1 + \frac{1}{\sqrt{d}} e_2$$

  where $\Gamma_i^{\sigma} = \Gamma_1$ if $i \in \Lambda$ and $\Gamma_i^{\sigma} = \Gamma_2$ otherwise.
- $D \overset{\mathrm{def}}{=} M_{\mathrm{dR}} = (\mathbf{Q}\hat{\omega}^{\otimes(k+j)} \oplus \mathbf{Q}\hat{\eta}^{\otimes(k+j)}) \otimes \mathbf{Q}\hat{\varepsilon}^{\otimes(-j)}$, and
  $Fil^0 = \mathbf{Q}\hat{\eta}^{\otimes(k+j)} \otimes \mathbf{Q}\hat{\varepsilon}^{\otimes(-j)}$
- $$V_p = \begin{cases} (T_{\mathfrak{p}}^{\otimes k} \otimes T_{\mathfrak{p}^*}^{\otimes(-j)} \oplus T_{\mathfrak{p}^*}^{\otimes k} \otimes T_{\mathfrak{p}}^{\otimes(-j)}) \otimes \mathbf{Q}_p \overset{\mathrm{def}}{=} \mathbf{Q}_p(\varphi_{\mathfrak{p}}) \oplus \mathbf{Q}_p(\varphi_{\mathfrak{p}^*}) & p \text{ splits} \\ T_{\mathfrak{p}}^{\otimes(k+j)}(-j) \otimes \mathbf{Q}_p \overset{\mathrm{def}}{=} K_{\mathfrak{p}}(\varphi_{\mathfrak{p}}) & p \text{ is nonsplit} \end{cases} \tag{2.2.2}$$

It is also known that $(V, D)$ gives rise to a motivic pair (see [5] section 1). Fix lattice $M = (\mathbf{Z}e_1 \oplus \mathbf{Z}e_2) \otimes \mathbf{Z}\Gamma^{\otimes(-j)}$ of $V$. It is shown in [5] that under comparision map $M$ is mapped onto the natural lattice $T_p$ where

$$T_p = \begin{cases} \mathbf{Z}_p(\varphi_{\mathfrak{p}}) \oplus \mathbf{Z}_p(\varphi_{\mathfrak{p}^*}), & \text{if } p \text{ splits} \\ \mathcal{O}_{\mathfrak{p}}(\varphi_{\mathfrak{p}}), & \text{if } p \text{ does not split} \end{cases} \tag{2.2.3}$$

It follows that $M \otimes \hat{\mathbf{Z}}$ is galois stable. Now fix a lattice $L = \mathbf{Z}\hat{\omega}^{\otimes(k+j)} \oplus \mathbf{Z}\hat{\eta}^{\otimes(k+j)} \otimes \mathbf{Z}\hat{\varepsilon}^{\otimes(-j)}$ of $D$. This gives a measure on $M_{\mathrm{dR}}/M^0_{\mathrm{dR}}$ by setting the measure of $L/L^0 = \mathbf{Z}\hat{\omega}^{\otimes(k+j)} \otimes \mathbf{Z}\hat{\varepsilon}^{\otimes(-j)}$ to be 1. We are going to verify our Main Theorem for this lattice $M$ with respect to Tamagawa measure induced by this chosen measure.

*Notation.* From now on, set $\varpi = \omega^{\otimes(k+j)} \otimes \varepsilon^{\otimes(-j)}$ and let $\hat{\varpi}$ denote $\hat{\omega}^{\otimes(k+j)} \otimes \hat{\varepsilon}^{\otimes(-j)}$.

Let us deal with $\mu_\infty(\mathbf{A}(\mathbf{R})/\mathbf{A}(\mathbf{Q}))$ first. $\mathbf{A}(\mathbf{Q}) \subseteq H^1_{f,spec\mathbf{Z}}(\mathbf{Q}, M \otimes \hat{\mathbf{Z}})$ and $\mathbf{A}(\mathbf{Q}) \otimes_{\mathbf{Z}} \mathbf{R} = D_\infty/(D^0_\infty + V^+_\infty) = 0$ in our case. We will see later in chapter6 (6.1.8) that $H^1_{f,spec\mathbf{Z}}(\mathbf{Q}, M \otimes \hat{\mathbf{Z}})$ is finite group. Hence $\mathbf{A}(\mathbf{Q})$ is finite torsion, and the $p$-part of $\#(\mathbf{A}(\mathbf{Q})) = \#(H^0(\mathbf{Q}, T_p \otimes \mathbf{Q}_p/\mathbf{Z}_p))$.

To find out $\mu_\infty(\mathbf{A}(\mathbf{R}))$ we need to compute a period integral. Notice under $\theta : H^1_{\mathrm{dR}}(E/\mathbf{Q}) \to H^1_{\mathrm{dR}}(E^{\mathrm{an}})$, we have (compare lemma 1.5 in [6])

$$\theta(\omega) = 2dz, \theta(\eta) = -\frac{2\pi}{\Omega^2_\infty} d\bar{z}$$

It is easy to compute that $\int_{\Gamma_1} dz = \Omega_\infty$, $\int_{\Gamma_2} dz = \sqrt{d}\Omega_\infty$, $\int_{\Gamma_1} d\bar{z} = \Omega_\infty$, $\int_{\Gamma_2} d\bar{z} = -\sqrt{d}\Omega_\infty$ hence under the map

$$\begin{aligned} &H^{-1}_{\mathrm{dR}}(E) \overset{(\theta^{\mathrm{dual}})^{-1}}{\longrightarrow} H^{-1}_{\mathrm{dR}}(E^{\mathrm{an}}) \to H^{-1}_B(E) \otimes \mathbf{C} \\ &\text{cycle} \mapsto \{\text{differential}\Delta \mapsto \int_{\text{cycle}} \Delta\} \\ &\hat{\omega} \mapsto \frac{1}{4\Omega_\infty}(\Gamma_1 + \frac{1}{\sqrt{d}}\Gamma_2) \text{ and } \hat{\eta} \mapsto \frac{-\Omega_\infty}{4\pi}(\Gamma_1 - \frac{1}{\sqrt{d}}\Gamma_2) \end{aligned} \tag{2.2.4}$$

Notice also under period isomorphism map for $\mathbb{G}_m$, $\Gamma^{\otimes}(-j) \mapsto (2\pi i)^{-j}\varepsilon^{-j}$. Hence by the above definition of $e_1$ amd $e_2$ , we have the following

**Proposition 1.** *Under our choice of Tamagawa measure and lattice $M$,*

$$\mu_\infty(\mathbf{A}(\mathbf{R})) = \begin{cases} (2\pi)^{-j}(4\Omega_\infty)^{k+j}, & \text{if } j \text{ is even} \\ (2\pi)^{-j}(4\Omega_\infty)^{k+j}\sqrt{-d}, & \text{if } j \text{ is odd} \end{cases}$$

For more details, see lemma 3.1 in [5].

**Lemma 1.**

$$\mu_p(\mathbf{A}(\mathbf{Q}_p)) = \frac{\mu_p(H^1_f(\mathbf{Q}_p, T_p))}{|P_p(V,1)|_p^{-1}} P_p(V,1) \text{ up to powers of 2 and 3} \tag{2.2.5}$$

*Proof.* Since $wt(\mathcal{M}_{k,j}) = -k \leqslant -2$, we have

$$\forall l \neq p, H^1_f(\mathbf{Q}_l, V_p) \cong H^1(\mathbf{Q}_l^{\mathrm{nr}}/\mathbf{Q}_l, V_p^{I_l}) \cong V_p^{I_l}/(1 - frob_l) = 0$$

Hence $\#(H^1_f(\mathbf{Q}_l, T_p)) = \#(H^0(\mathbf{Q}_l, V_p/T_p))$. It follows from the exact sequence of $I_l$ modules:

$$0 \to T_p \to V_p \to V_p/T_p \to 0$$

that

$$0 \to V_p^{I_l}/T_p^{I_l} \to (V_p/T_p)^{I_l} \to H^1(I_l, T_p)_{\mathrm{tors}} \to 0$$

Since $H^0(\mathbf{Q}_l, V_p/T_p) = H^0(\mathbf{Q}_l^{\mathrm{nr}}/\mathbf{Q}_l, (V_p/T_p)^{I_l})$, we get from the above exact sequence:

$$\#(H^0(\mathbf{Q}_l, V_p/T_p)) = \#(H^0(\mathbf{Q}_l^{\mathrm{nr}}/\mathbf{Q}_l, V_p^{I_l}/T_p^{I_l})) \times \#(H^0(\mathbf{Q}_l^{\mathrm{nr}}/\mathbf{Q}_l, H^1(I_l, T_p)_{\mathrm{tors}})$$

Obviously $\#(H^0(\mathbf{Q}_l^{\mathrm{nr}}/\mathbf{Q}_l, V_p^{I_l}/T_p^{I_l})) = |P_l(V,1)|_p^{-1}$. On the other hand, $\#(H^0(\mathbf{Q}_l^{\mathrm{nr}}/\mathbf{Q}_l,$ $H^1(I_l, T_p)_{\mathrm{tors}})) = 1$ if $I_l$ acts trivially on $T_p$ since then $H^1(I_l, T_p) \cong Hom(I_l, T_p)$ which has no torsions. In general since $l \neq p, H^1(I_l, T_p) \cong H^1(I_l/R, T_p) = T_p/1 - \sigma_l$ where $R$ is the wild inertia at $l$ and $\sigma_l$ is a topologically generator of $I_l/R$. But since our ellipic curve has CM, $I_l$ acts on $T_p$ through finite quotient. It is well known that eigenvalues of $\sigma_l$ are at most 24 th root of unity. Hence $\#(H^1(I_l, T_p)_{\mathrm{tors}})$ is divisible by at most powers of 2 and 3. ∎

Now the Tamagawa number conjecture for $\mathcal{M}_{k,j}$ is reduced to

$$\prod_{p\leqslant\infty}\frac{\mu_p(H^1_f(\mathbf{Q}_p,T_p))}{|P_p(V,1)|_p^{-1}}=\frac{L(\bar{\psi}^{k+j},k)}{(2\pi)^{-j}(4\Omega_\infty)^{k+j}(\sqrt{-d})^{\delta(j)}}\#(\mathbf{A}(\mathbf{Q}))\frac{\#H^0(\mathbf{Q},M^*\otimes\mathbf{Q}/\mathbf{Z}(1)}{\#(\bigsqcup(M))} \tag{2.2.6}$$

where $\delta(j)=1$ if $j$ is odd and 0 otherwise.

## 2.3 More preliminary results and notations

Let $\mathfrak{p}$ be a prime ideal of $\mathcal{O}_K$. Let $\hat{E}_\mathfrak{p}$ be the formal group of E around origin over $K_\mathfrak{p}$ which gives the kernel of reduction modulo $\mathfrak{p}$. Let $L$ be the period lattice of a global minimal Weierstrass model of $E$. Define

$$\varepsilon:\mathbb{C}/L\to E(\mathbb{C}),\varepsilon(z)=(x,y),x=\mathcal{P}(z),y=\mathcal{P}^{'}(z)$$

where $\mathcal{P}(z)$ is the Weierstrass function. We have $E_1(K_\mathfrak{p})\stackrel{\cong}{\to}\hat{E}_\mathfrak{p}$ under which $(x,y)\mapsto T=-2x/y$. Here $E_1(K_\mathfrak{p})$ is the kernel of reduction modulo $\mathfrak{p}$ on $E(K_\mathfrak{p})$. Fix $T$ as the variable on $\hat{E}_\mathfrak{p}$. Let $\omega_0$ be the normalized invariant differential on $\hat{E}_\mathfrak{p}$.

**Lemma 2.**

$$\omega=\omega_0$$

*Proof.* Let $\lambda$ be the inverse of $\varepsilon$, hence $z=\lambda(T)$. We have

$$\omega=\frac{dx}{y}=dz=\lambda^{'}(T)dT,\quad \omega_0=\lambda^{'}_{\hat{E}_\mathfrak{p}}(T)dT$$

where $\lambda_{\hat{E}_\mathfrak{p}}$ is the logarithm assoicated to $\omega$. Since both $\lambda$ and $\lambda_{\hat{E}_\mathfrak{p}}$ are group homomorphism : $\hat{E}_\mathfrak{p}\to\mathbb{G}_a$, we have $\lambda_{\hat{E}_\mathfrak{p}}=c\lambda$ for some constant c. But well known formulas, $x=T^{-2}a(T),y=T^{-3}a(T)$, where $a(T)$ has coefficents in $\mathcal{O}_\mathfrak{p}$ and constant term is 1. Hence $dx/y=(1+\cdots)dT$. It follows that $c=1$ as $\omega_0$ is normalized. ☐

*Notation.* From now on we fix a prime $\mathfrak{p}$ of $\mathcal{O}_K$ and the completion of $\mathcal{O}_K$ at $\mathfrak{p}$ by $\mathcal{O}_\mathfrak{p}$. When no confusion arise, we will sometimes omit the subscript $\mathfrak{p}$.

Now if $\mathfrak{p}$ is a good reduction prime for $E$, i.e $\mathfrak{p} \nmid \mathfrak{f}$, then $\hat{E}$ is a Lubin-Tate formal group with respect to the uniformizer $\pi = \psi(\mathfrak{p})$ of $\mathfrak{p}$. When $E$ has bad reduction at $\mathfrak{p}$, we will replace $E$ by a twist which has good reduction at $\mathfrak{p}$. This is possible since $E$ has CM and so $E$ must has potentially good reduction at $\mathfrak{p}$. More precisely fix character $\epsilon : Gal(\bar{K}/K) \to \mathcal{O}_K^*$ such that $\psi\epsilon$ is unramified at $\mathfrak{p}$, and such that cond$(\epsilon) = \mathfrak{g}$ the prime to $\mathfrak{p}$ part of $\mathfrak{f}$. Let $M$ denote the fixed field of Ker$(\epsilon)$. Let $E'$ denotes the twist of $E$ by $\epsilon$: $E'$ is an elliptic curve over $K$ which is isomorphic to $E$ over $M$. Fix $f : E \to E'$ over $M$. Let $\psi', \omega'$ and $\Omega'_\infty$ be defined the same way for $E'$ as for $E$. It is known that $\psi' = \psi\epsilon$. Set r such that $f^*\omega' = r\omega$. Then we also have $\Omega'_\infty = r\Omega_\infty$. Let $F$ denote the maximal extension of $K$ inside $M(E_{\mathfrak{p}^\infty}) = M(E'_{\mathfrak{p}^\infty})$ which is unramified at $\mathfrak{p}$. Let $F_n = F(E'_{\mathfrak{p}^n})$. It is known ([10] lemma 4.2 (ii)) that $F_n \supset MK_n, \forall n \geqslant m(\mathfrak{p})$ where $m(\mathfrak{p}) \overset{\text{def}}{=}$ the smallest integer such that $1 + \mathfrak{p}^{m(\mathfrak{p})}\mathcal{O} \subset (\mathcal{O}^*)^{w_K}$, $w_K = \#(\mathcal{O}_K^*)$.

*Notation.*

- when $p$ splits, $\forall n \in \mathbf{N}, \forall m \in \mathbf{N}$, set $K_{n,m} = K(E_{\mathfrak{p}^n\mathfrak{p}^{*m}})$ and let $K_{n,m,\mathfrak{P}}$ be the completion of $K_{n,m}$ with respect to the place of $K_{n,m}$ above $\mathfrak{p}$ induced by $\mathfrak{P}$ which we still denote by $\mathfrak{P}$. By convention $K_n$ means $K_{n,0}$ in this case.

- when $p$ is nonsplit, $\forall n \in \mathbf{Z}^+$, set $K_n = K(E_{\mathfrak{p}^n}) = K(E_{p^n})$ and set $K_{n,\mathfrak{P}}$ as the completion at $\mathfrak{P}$.

In all cases, set $G_\infty = Gal(K_\infty/K)$ and $G_{\infty,v} = Gal(K_{\infty,v}/K_\mathfrak{p})$, where

$$K_\infty \overset{\text{def}}{=} \begin{cases} K(E_{p^\infty}), & \text{if } j > 0 \\ K(E_{\mathfrak{p}^\infty}), & \text{if } j = 0 \end{cases}$$

Set $\Lambda = \mathbf{Z}_p[[G_\infty]]$ the usual Iwasawa algebra. In the case that $E$ has bad reduction at $\mathfrak{p}$ we aslo define the same way $F_{n,m}$ or $F_n$ and $F_\infty$ for base field $F$, but replacing $E$ by $E'$. For example, $F_n \overset{\text{def}}{=} F(E'_{\mathfrak{p}^n})$.

**Lemma 3.** *There are only finite number of places $\mathfrak{P}$ of $F_\infty$ and $K_\infty$ lying above $\mathfrak{p}$.*

*Proof.* We know the property of ramification of the fields under consideration from theory of complex multiplication. When $p$ is nonsplit, each $K_n$ is totally ramified over $K$ at $\mathfrak{p}$ hence there is a unique place of $K_\infty$ above $\mathfrak{p}$ , namely ( the one induced by ) $\upsilon$ and the set of places of $F_\infty$ lying above $\mathfrak{p}$ corresponds 1 - 1 with $\Sigma \stackrel{\text{def}}{=}$ a set of coset representative of $Gal(F/K)/Gal(F_\upsilon/K_\mathfrak{p}) : \forall \mathfrak{P} \mid \mathfrak{p}$, there exists a unique $\sigma \in \Sigma$ such that $\mathfrak{P} = \sigma\upsilon$ , then $\mathfrak{P}$ corresponds to $\sigma$.

When $p$ splits, let us first assume that $E$ has good reduction at $\mathfrak{p}$. Then $K_{n,m}$ is totally ramified over $K_{0,m}$ at places above $\mathfrak{p}$ and $K_{0,m}$ is unramified over $K$ at $\mathfrak{p}$. Moreover from the definition of Grössencharacter, it follows that the number of primes of $K_{0,m}$ above $\mathfrak{p}$ which we denote by $r_m$ is given by the index of the subgroup generated by $\pi$ in $(\mathcal{O}_\mathfrak{p}/\mathfrak{p}^{*m})^*$. Hence there is an interger $M$ such that $r_m = r_0 p^m$ for $m < M$ and $r_m = r_0 p^M$ for $m \geq M$. So there are only finite number of places of $K_\infty$ above $\mathfrak{p}$. When $\mathfrak{p}$ is a bad reduction prime for $E$, upon replacing $K$ by $F$ and $E$ by $E^{'}$ which has good reduction at $\mathfrak{p}$ , the above argument applies to $F_{n,m}$ and $F_\infty$. Recall $E$ and $E^{'}$ are isomorphic over $F_{m(\mathfrak{p}),0}$, hence $F_\infty$ is an overfield of $K_\infty$. From this we see there are only finite number of places of $K_\infty$ above $\mathfrak{p}$ too. In fact there is 1 - 1 correspondence between the set of places of $F_\infty$ above $\mathfrak{p}$ and a set of coset representative of $Gal(F_{0,M}/K)/Gal(F_{0,M,\upsilon}/K_\mathfrak{p})$ . $\square$

# CHAPTER 3

# LOCAL TAMAGAWA NUMBER; DUAL EXPONENTIAL MAP

## 3.1 Local Tamagawa number

In this section we will see how to reduce the computation of local Tamagawa number to an evaluation of dual exponential map $exp^*$ on a special element.

Recall in our case,

$$H^1(\mathbf{Q}_p, V_p) = \begin{cases} H^1(\mathbf{Q}_p, \mathbf{Q}_p(\varphi_{\mathfrak{p}})) \oplus H^1(\mathbf{Q}_p, \mathbf{Q}_p(\varphi_{\mathfrak{p}^*})), & \text{if } p \text{ splits} \\ H^1(\mathbf{Q}_p, K_{\mathfrak{p}}(\varphi_{\mathfrak{p}})), & \text{if } p \text{ does not split} \end{cases} \tag{3.1.1}$$

where $\varphi_{\mathfrak{p}} = \psi_{\mathfrak{p}}^k \bar{\psi}_{\mathfrak{p}}^{-j}$. Note that when $p$ splits, $H^1_f(\mathbf{Q}_p, V_p) = H^1(K_{\mathfrak{p}}, K_{\mathfrak{p}}(\varphi_{\mathfrak{p}}))$ since $DR^0(V_p) = DR(\mathbf{Q}_p(\varphi_{\mathfrak{p}^*}))$, from which it follows that $H^1_f(\mathbf{Q}_p, \mathbf{Q}_p(\varphi_{\mathfrak{p}^*})) = 0$. When $p$ is nonsplit, we have

$$H^1(\mathbf{Q}_p, V_p) \cong H^1(K_{\mathfrak{p}}, K_{\mathfrak{p}}(\varphi_{\mathfrak{p}}))^{Gal(K_{\mathfrak{p}}/\mathbf{Q}_p)}$$

$$DR(\mathbf{Q}_p, V_p) \otimes K_{\mathfrak{p}} = DR(K_{\mathfrak{p}}, K_{\mathfrak{p}}(\varphi_{\mathfrak{p}}))$$

We will first concentrate on $DR(K_{\mathfrak{p}}, K_{\mathfrak{p}}(\varphi_{\mathfrak{p}})) \stackrel{exp}{\to} H^1_f(K_{\mathfrak{p}}, K_{\mathfrak{p}}(\varphi_{\mathfrak{p}}))$ and then if $p$ is nonsplit, we will take into consideration the $Gal(K_{\mathfrak{p}}/\mathbf{Q}_p)$ action later.

Let $a \in K_{\mathfrak{p}}$ be the local Tamagawa measure of $H^1_f(K_{\mathfrak{p}}, \mathcal{O}_{\mathfrak{p}}(\varphi_{\mathfrak{p}}))$ with respect to lattice $L/L^0 \otimes \mathcal{O}_{\mathfrak{p}} \subseteq D/D^0 \otimes K_{\mathfrak{p}}$, i.e $a = a_1 \#(H^1(K_{\mathfrak{p}}, \mathcal{O}_{\mathfrak{p}}(\varphi_{\mathfrak{p}})))_{\text{tors}}$ where $a_1$ is the index of lattice $H^1_f(K_{\mathfrak{p}}, \mathcal{O}_{\mathfrak{p}}(\varphi_{\mathfrak{p}}))^{\text{tf}}$ and the lattice $exp(L/L^0 \otimes \mathcal{O}_{\mathfrak{p}})$ in $H^1_f(K_{\mathfrak{p}}, K_{\mathfrak{p}}(\varphi_{\mathfrak{p}}))$, i.e

$$exp(L/L^0 \otimes \mathcal{O}_{\mathfrak{p}}) = a_1 \times H^1_f(K_{\mathfrak{p}}, \mathcal{O}_{\mathfrak{p}}(\varphi_{\mathfrak{p}}))^{\text{tf}} \tag{3.1.2}$$

Note both sides are rank 1 $\mathcal{O}_\mathfrak{p}$ modules. Here "tf" means torsion free part. The restriction map gives an isomorphism:

$$\begin{aligned} H^1(K_\mathfrak{p}, K_\mathfrak{p}(\varphi_\mathfrak{p})) &\cong H^1(K_{\infty,v}, K_\mathfrak{p}(\varphi_\mathfrak{p}))^{G_{\infty,v}} \\ &\cong [\oplus_{\mathfrak{P}|\mathfrak{p}} H^1(K_{\infty,\mathfrak{P}}, K_\mathfrak{p}(\varphi_\mathfrak{p}))]^{G_\infty} \\ &\cong Hom_{G_\infty}(\Xi_\infty, K_\mathfrak{p}(\varphi_\mathfrak{p})) \quad \text{by local class field theory} \end{aligned} \tag{3.1.3}$$

where $\Xi_\infty = \oplus_{\mathfrak{P}|\mathfrak{p}} \mathbf{U}_{\infty,\mathfrak{P}}$ and

$$\mathbf{U}_{\infty,\mathfrak{P}} \overset{\text{def}}{=} \begin{cases} \lim_{\leftarrow n} K^*_{n,\mathfrak{P}} & \text{if } j = 0 \\ \lim_{\leftarrow n,m} K^*_{n,m,\mathfrak{P}} & \text{if } j > 0 \end{cases}$$

Let $\Xi'_\infty$ and $\mathbf{U}'_{\infty,\mathfrak{P}}$ be similarly defined for F. We denote $\mathfrak{A}_\infty$ to be inverse limit of global units similarly defined and we view $\mathfrak{A}_\infty$ as embeded inside $\Xi_\infty$ and denote the closure of the image as $\overline{\mathfrak{A}_\infty}$. Denote $\mathfrak{C}_\infty$ and $\overline{\mathfrak{C}_\infty}$ to be the subgroup defined for elliptic units.

*Remark 1.* The reseason we write the semi-local version of the above isomorphism is because of the link with Iwasawa Main conjecture.

Now $H^1_f(K_\mathfrak{p}, \mathcal{O}_\mathfrak{p}(\varphi_\mathfrak{p}))^{\text{tf}} \subseteq Hom_{G_\infty}(\Xi_\infty, K_\mathfrak{p}(\varphi_\mathfrak{p}))$ as a rank 1 free $\mathcal{O}_\mathfrak{p}$ module. There exists integer $e$ such that

$$H^1_f(K_\mathfrak{p}, \mathcal{O}_\mathfrak{p}(\varphi_\mathfrak{p}))^{\text{tf}} = \mathfrak{p}^e Hom_{G_\infty}(\Xi_\infty, \mathcal{O}_\mathfrak{p}(\varphi_\mathfrak{p})) \tag{3.1.4}$$

**Lemma 4.**

$$(\mathbf{N}\mathfrak{p})^e = \#(H^0(K_\mathfrak{p}, K_\mathfrak{p}/\mathcal{O}_\mathfrak{p}(\varphi_\mathfrak{p})))$$

*Proof.* A special case of this in the case $j = 0$ and $p$ splits is given in [6] prop 3.8, but the proof extends to our case. Nothing in that proof relies on $j = 0$ or $p$ being a split prime. ∎

Fix $\delta = exp(\hat{\varpi})$ a generator of $Hom_{G_\infty}(\Xi_\infty, K_\mathfrak{p}(\varphi_\mathfrak{p}))$. Suppose $\delta$ maps $\Xi_\infty$ to $b\mathcal{O}_\mathfrak{p}(\varphi_\mathfrak{p}) \subset K_\mathfrak{p}(\varphi_\mathfrak{p})$, i.e

$$exp(L/L^0 \otimes \mathcal{O}_\mathfrak{p}) = bHom_{G_\infty}(\Xi_\infty, \mathcal{O}_\mathfrak{p}(\varphi_\mathfrak{p})) \tag{3.1.5}$$

Comparing 3.1.2, 3.1.4 and 3.1.5, we get $|b|_{\mathfrak{p}} = |a_1|_{\mathfrak{p}} (\mathbf{N}\mathfrak{p})^e$. Since

$$\#(H^1_f(K_{\mathfrak{p}}, \mathcal{O}_{\mathfrak{p}}(\varphi_{\mathfrak{p}}))_{\text{tors}}) = \#(H^0(K_{\mathfrak{p}}, K_{\mathfrak{p}}/\mathcal{O}_{\mathfrak{p}}(\varphi_{\mathfrak{p}})))$$

we have $|a|_{\mathfrak{p}} = |b|_{\mathfrak{p}}$. Our approach to find $b$ is first calculate $\delta(u)$ for some $u$ which is a generator of rank 1 $\Lambda$ module $\mathfrak{C}_\infty$. In particular we see $\delta(u) \neq 0$. By well known results (e.g [10] Theorem5.1(ii)) in Iwasawa theory, $\Xi_\infty$ is of rank $[K_{\mathfrak{p}} : \mathbf{Q}_p]$ as $\Lambda$ module. So when $\mathfrak{p} \mid p$ and $p$ splits in $K/\mathbf{Q}$, then $\Xi_\infty/\overline{\mathfrak{C}_\infty}$ is a torsion $\Lambda$ module. Since $\delta$ is $G_\infty$ equivariant, we have

$$b = \delta(u) \times (\#([(\Xi_\infty/\overline{\mathfrak{C}_\infty})(\varphi_{\mathfrak{p}}^{-1})]_{G_\infty}))^{-1} \tag{3.1.6}$$

When $p$ is nonsplit, it is known ([10] lemma 11.9) there is a noncanonical direct sum decomposition $\Xi_\infty = \Xi^1_\infty \oplus \Xi^2_\infty$, where $\Xi^1_\infty$ and $\Xi^2_\infty$ are free rank 1 $\Lambda$ modules such that $\delta(\Xi^2_\infty) = 0$. Since $\delta(u) \neq 0$, it follows $\Xi_\infty/(\Xi^2_\infty, \overline{\mathfrak{C}_\infty})$ is a torsion $\Lambda$ module. Now $\delta$ viewed as a homorphism: $\Xi_\infty \to K_{\mathfrak{p}}(\varphi_{\mathfrak{p}})$ factors through $\Xi_\infty/\Xi^2_\infty$ and is also $G_\infty$ equivariant, so

$$b = \delta(u) \times \#(Hom_{G_\infty}(\Xi_\infty/(\Xi^2_\infty, \overline{\mathfrak{C}_\infty}), K_{\mathfrak{p}}/\mathcal{O}_{\mathfrak{p}}(\varphi_{\mathfrak{p}}))) \tag{3.1.7}$$

In the case $p$ splits, one can give a formula for $\#([(\Xi^1_\infty/\overline{\mathfrak{C}_\infty})(\varphi_{\mathfrak{p}}^{-1})]_{G_\infty})$ using (1 or 2 variable) $p$-adic L-function and Iwasawa Main conjecture. But this is not possible in the case $p$ is nonsplit. So instead of explicitly computing it, we relate it to $\#(\bigsqcup(M))$ and in verfying 2.2.6 this term get canceled. On the other hand, the claculation $\delta(u)$ turns out to be very complicated, especially when $j > 0$. First we will reduce it to an evaluation of $exp^*$ map.

## 3.2 Dual exponential map

**Definition 3.** $K$ will be a finite extension of $\mathbf{Q}_p$ in this section. Let $V$ be a $\mathbf{Q}_p$ vector space endowed with a continous action of $Gal(\bar{K}/K)$. Assume $V$ is a de Rham representation. There are 2 equivalent definition of $exp^*$ (for the proof of the equivalence, see [7] p124-125.)

1. $exp^* : H^1(K,V)/H^1_f(K,V) \to D^0_{\mathrm{dR}}(V)$ is the composite of the following maps induced by duality pairing:

$$\begin{aligned} H^1(K,V)/H^1_f(K,V) &\overset{\cong}{\to} Hom_{\mathbf{Q}_p}(H^1_f(K,V^*(1)),\mathbf{Q}_p) \qquad (3.2.1)\\ &\to Hom_{\mathbf{Q}_p}(D_{\mathrm{dR}}(V^*(1))/D_{\mathrm{dR}}(V^*(1))^0,\mathbf{Q}_p) \overset{\cong}{\to} D_{\mathrm{dR}}(V)^0 \end{aligned}$$

2. $exp^*$ is the composite map

$$H^1(K,V) \to H^1(K,V\otimes B^+_{\mathrm{dR}}) \overset{\cong}{\to} D^0_{\mathrm{dR}}(V)$$

where the second arrow is the inverse of isomorphism:

$$D^i_{\mathrm{dR}}(V) = H^0(K,V\otimes_{\mathbf{Q}_p} B^i_{\mathrm{dR}}) \overset{\cup(log(\chi_{\mathrm{cyclo}}))}{\longrightarrow} H^1(K,V\otimes_{\mathbf{Q}_p} B^i_{\mathrm{dR}}) \qquad (3.2.2)$$

in the case $i=0$. Here $\chi_{\mathrm{cyclo}} : Gal(\bar{K}/K) \to \mathbf{Z}_p^*$ is the character describing its action on $p^\infty$ th root of 1, and

$$log(\chi_{\mathrm{cyclo}}) \in H^1(K,\mathbf{Z}_p) = Hom_{\mathrm{cont}}(Gal(\bar{K}/K),\mathbf{Z}_p^*).$$

In our situation, we have

$$DR(K_{\mathfrak{p}},K_{\mathfrak{p}}(\varphi_{\mathfrak{p}})) \overset{\exp}{\to} H^1_f(K_{\mathfrak{p}},K_{\mathfrak{p}}(\varphi_{\mathfrak{p}})) \overset{\mathrm{res}}{\to} Hom(\Xi_\infty,K_{\mathfrak{p}}(\varphi_{\mathfrak{p}}))$$

or $DR(K_{\mathfrak{p}},V_{\mathfrak{p}}(E)^{\otimes(k+j)}(-j)) \overset{\exp}{\to} H^1_f(K_{\mathfrak{p}},V_{\mathfrak{p}}(E)^{\otimes(k+j)}(-j))$ and we want to evaluate

$$\delta = \oplus_{\mathfrak{P}|\mathfrak{p}} exp_{\mathfrak{P}}(\hat{\varpi}) : \Xi_\infty = \oplus_{\mathfrak{P}|\mathfrak{p}}\mathbf{U}_{\infty,\mathfrak{P}} \to K_{\mathfrak{p}}(\varphi_{\mathfrak{p}})$$

on some special element $u$ comming from global units, i.e $u = (\cdots, u_{\mathfrak{P}}, \cdots)_{\mathfrak{P}|\mathfrak{p}}$.

*Remark 2.* We are going to write down explicitly what $u$ is later on. Right now it suffices to know that when $p$ is a bad reduction prime, $u = Norm(u')$ for some $u' \in F^*_\infty$. We have the following commutative diagram:

$$\begin{array}{ccc} F^*_\infty & \longrightarrow & \oplus_{\mathfrak{P}|\mathfrak{p}}(\oplus_{\omega|\mathfrak{P}}F^*_{\infty,\omega}) \\ \mathrm{norm}\downarrow & & \downarrow \oplus_{\mathfrak{P}|\mathfrak{p}} \text{ local norm} \\ K^*_\infty & \longrightarrow & \oplus_{\mathfrak{P}|\mathfrak{p}}K^*_{\infty,\mathfrak{P}} \end{array}$$

Recall for a fixed place $\mathfrak{P} \mid \mathfrak{p}$ of $K_\infty$ , the set of places $\omega \mid \mathfrak{P}$ of $F_\infty$ corresponds 1 - 1 with a set of coset representative of $Gal(F/K)/Gal(F_v/K_{\mathfrak{p}}))$.

We need to consider

$$exp^* : H^1(K_{\mathfrak{p}}, V_{\mathfrak{p}}^{\otimes -(k+j)}(1+j)) \to DR^0(K_{\mathfrak{p}}, V_{\mathfrak{p}}^{\otimes -(k+j)}(1+j))$$

First we define special elements of $H^1(K_{\mathfrak{p}}, V_{\mathfrak{p}}^{\otimes(-k-j)}(1+j))$ by method of Soule:

*Notation.* Fix a basis $\epsilon$ of $\mathbf{Z}_p(1)$ once for all. Fix a basis $\xi = (\xi_n)_n$ for $T_{\mathfrak{p}}(E)$ as an $\mathcal{O}_{\mathfrak{p}}$ module. We will specify a particular choice later on in connection with Coleman power series. This also determines a basis $\zeta$ for $T_{\mathfrak{p}^*}$ when $\mathfrak{p} \mid p$ splits, by requiring the Weil pairing of $\xi$ and $\zeta$ to be $\epsilon$.

**Definition 4.** (1): if $j = 0$

(i): if $p$ is a good reduction prime of $E$, then $\forall n \geq 0$, denote the image of $u_{\mathfrak{P}}$ under the folowing composite map as $S_n(u_{\mathfrak{P}})$.

$$\begin{aligned} \varprojlim_m K^*_{m,\mathfrak{P}} \to \varprojlim_m H^1(K_{m,\mathfrak{P}}, \mathbf{Z}_p(1)) \overset{*}{\to} \varprojlim_m H^1(K_{m,\mathfrak{P}}, T_{\mathfrak{p}}^{\otimes(-k)}(1)/\pi^m) \\ \overset{\text{trace}}{\to} \varprojlim_m H^1(K_{n,\mathfrak{P}}, T_{\mathfrak{p}}^{\otimes(-k)}(1)/\pi^m) \cong H^1(K_{n,\mathfrak{P}}, T_{\mathfrak{p}}^{\otimes(-k)}(1)) \end{aligned} \tag{3.2.3}$$

where the map (*) is $\cup \xi_m^{\otimes(-k)}$ which is well defined because by the very definition of $K_m$, all $\mathfrak{p}^m$ torsion points of $E$ are rational over $K_m$.

(ii) if $p$ is a bad reduction prime, $\forall n \geq m(\mathfrak{p})$ we define for $\omega \mid \mathfrak{P} \mid \mathfrak{p}$, $S_n(u'_\omega)$ to be the image of $u'_\omega \in \lim_m F^*_{m,\omega}$ under analogous composite map as above defined for base field $K_{\mathfrak{p}}$. Define $S_n(u_{\mathfrak{P}}) = \sum_{\omega \mid \mathfrak{P}} Cores_{F_{n,\omega}/K_{n,\mathfrak{P}}}(S_n(u'_\omega))$ and define $S_0(u_{\mathfrak{P}}) = Cores_{K_{n,\mathfrak{P}}/K_{\mathfrak{p}}}(S_n(u_{\mathfrak{P}}))$. Of course it is clear that $S_0(u_{\mathfrak{P}})$ is independent of $n$. Let $S'_n(u'_\omega) \in H^1(F_{n,\omega}, T_{\mathfrak{p}}(E')^{\otimes(-k)}(1))$ has the same meaning, except now we use the basis $\xi'$ in the definition.

(2): if $j > 0$,

(i): if $p$ split is a good reduction prime, $\forall n > 0$, define $S_n(u_{\mathfrak{P}})$ to be the image of $u_{\mathfrak{P}}$ under the composite of the following maps:

$$\begin{aligned} \varprojlim_m K^*_{m,m,\mathfrak{P}} &\to \varprojlim_m H^1(K_{m,m,\mathfrak{P}}, \mathbf{Z}_p(1)) \\ &\overset{(*)}{\to} \varprojlim_m H^1(K_{m,m,\mathfrak{P}}, T_{\mathfrak{p}}^{\otimes(-k)} \otimes T_{\mathfrak{p}^*}^{\otimes j}(1)/\pi^m) \\ &\overset{\text{trace}}{\to} \varprojlim_m H^1(K_{1,n,\mathfrak{P}}, T_{\mathfrak{p}}^{\otimes(-k)} \otimes T_{\mathfrak{p}^*}^{\otimes j}(1)/\pi^m) \\ &\cong H^1(K_{1,n,\mathfrak{P}}, T_{\mathfrak{p}}^{\otimes(-k)} \otimes T_{\mathfrak{p}^*}^{\otimes j}(1)) \end{aligned} \tag{3.2.4}$$

where the map (*) is $\cup \xi_m^{\otimes(-k)} \otimes \zeta_m^{\otimes j}$which is well defined again by definition of $K_{m,m}$. We denote the image of $S_n(u_{\mathfrak{P}})$ in $H^1(K_{1,n,\mathfrak{P}}, T_{\mathfrak{p}}^{\otimes(-k)} \otimes T_{\mathfrak{p}^*}^{\otimes j}(1)/\pi^n)$ as $\overline{S_{n,n}}(u_{\mathfrak{P}})$. Define $S_0(u_{\mathfrak{P}}) = Cores_{K_{1,n,\mathfrak{P}}/K_{\mathfrak{p}}} S_n(u_{\mathfrak{P}})$ and $\overline{S_{0,n}}(u_{\mathfrak{P}}) = Cores_{K_{1,n,\mathfrak{P}}/K_{\mathfrak{p}}} \overline{S_{n,n}}(u_{\mathfrak{P}})$.

(ii): if $p$ is a nonsplit and good reduction prime, $\forall n \geq 0$, define $S_n(u_{\mathfrak{P}})$ as the image of $u_{\mathfrak{P}}$ under the composite of the following maps:

$$\begin{aligned} \varprojlim_m K^*_{m,\mathfrak{P}} \to \varprojlim_m H^1(K_{m,\mathfrak{P}}, \mathbf{Z}_p(1)) &\overset{(*)}{\to} \varprojlim_m H^1(K_{m,\mathfrak{P}}, T_{\mathfrak{p}}^{\otimes(-k-j)}(1+j)/\pi^m) \\ &\overset{\text{trace}}{\to} \varprojlim_m H^1(K_{n,\mathfrak{P}}, T_{\mathfrak{p}}^{\otimes(-k-j)}(1+j)/\pi^m) \\ &\cong H^1(K_{n,\mathfrak{P}}, T_{\mathfrak{p}}^{\otimes(-k-j)}(1+j)) \end{aligned} \tag{3.2.5}$$

where the map (*) is $\cup \xi_m^{\otimes(-k-j)} \otimes \epsilon^{\otimes j}$ which is well defined again by definition of $K_m$. We denote the image of $S_n(u_{\mathfrak{P}})$ in $H^1(K_{n,\mathfrak{P}}, T_{\mathfrak{p}}^{\otimes(-k-j)}(1+j)/\pi^n)$ as $\overline{S_{n,n}}(u_{\mathfrak{P}})$. Define $\overline{S_{0,n}}(u_{\mathfrak{P}}) = Cores_{K_{n,\mathfrak{P}}/K_{\mathfrak{p}}} \overline{S_{n,n}}(u_{\mathfrak{P}})$ which is just the image of $S_0(u_{\mathfrak{P}})$ in $H^1(K_{\mathfrak{p}}, T_{\mathfrak{p}}^{\otimes(-k-j)}$ $(1+j)/\mathfrak{p}^n)$.

(iii): when $p$ is a bad reduction prime, just as in the case $j = 0$, $\forall n \geq m(\mathfrak{p})$, define $S_n(u'_\omega)$ as the image of $u'_\omega$ in $H^1(F_{m(\mathfrak{p}),n,\omega}, T_{\mathfrak{p}}^{\otimes(-k-j)}(1+j))$ when $\mathfrak{p} \mid p$ , $p$ splits and $S_n(u'_\omega) \in H^1(F_{n,\omega}, T_{\mathfrak{p}}^{\otimes(-k-j)}(1+j))$ when $p$ nonsplits. Define $S_n(u_{\mathfrak{P}}) = \sum_{\omega|\mathfrak{P}} Cores(S_n(u'_\omega))$ and define $S_0(u_{\mathfrak{P}}) = Cores(S_n(u_{\mathfrak{P}}))$. Again this definition does not depend on the $n \geq m(\mathfrak{p})$ chosen. Define $\overline{S_{n,n}}(u_{\mathfrak{P}})$ the same way as when $\mathfrak{p}$ is a good reduction prime. Finally denote again $S'_n(u'_\omega)$ and $\overline{S'_{n,n}}(u'_\omega)$ similarly defined elements except we replace $T_{\mathfrak{p}}(E)$ by $T_{\mathfrak{p}}(E')$ and use $\xi'$ in the definition. It is clear by definition and our chosing $f(\xi) = \xi'$ that $fS_n(u'_\omega) = S'_n(u'_\omega)$, where $f : E \to E'$.

*Remark 3.* When $\mathfrak{p}$ is a bad reduction prime, we could have defined $S_n(u_{\mathfrak{P}})$ in another way, namely directly apply Soule's method to $u_{\mathfrak{P}}$ just as in the good reduction prime case. The following commutative diagram (in the case $j = 0$) tell us the two definitions are equivalent:

$$\begin{array}{ccc}
\oplus_{\omega|\mathfrak{P}}\mathbf{U}'_{\infty,\omega} & \xrightarrow{\oplus_{\omega|\mathfrak{P}}\text{local norm}} & \mathbf{U}_{\infty,\mathfrak{P}} \\
\downarrow & & \downarrow \\
\oplus_{\omega|\mathfrak{P}}H^1(F_{n,\omega}, T_{\mathfrak{p}}^{\otimes(-k)}(1)) & \xrightarrow{\sum_{\omega|\mathfrak{P}}\text{cores}} & H^1(K_{n,\mathfrak{P}}, T_{\mathfrak{p}}^{\otimes(-k)}(1))
\end{array}$$

The reason for the approach we defined them earlier is that in the process of calculating $exp^*(u)$ when $\mathfrak{p}$ is a bad reduction prime, we need to use $S_n(u')$ and $S'_n(u')$ anyway. Note the same remark applies to the definition of $S_n(u)$ for $\mathfrak{p}$ bad reduction prime and $j > 0$.

Now let us come back to the problem of relating $\delta(u)$ to $\sum_{\mathfrak{P}|\mathfrak{p}} exp^*(S_0(u_{\mathfrak{P}}))$. The crucial point is the following commutative diagram from Galois cohomology theory: (we have in mind now $j = 0$ case, $j > 0$ case follows the same way)

$$\begin{array}{ccccc}
H^1(K_{\mathfrak{p}}, K_{\mathfrak{p}}(\varphi_{\mathfrak{p}})) & \times & H^1(K_{\mathfrak{p}}, \mathcal{O}_{\mathfrak{p}}(\varphi_{\mathfrak{p}})^*(1)) & \longrightarrow & K_{\mathfrak{p}} \\
\downarrow{\scriptstyle \text{res}} & & \uparrow{\scriptstyle \sum_{\mathfrak{P}|\mathfrak{p}}\text{cores}} & & \| \\
\oplus_{\mathfrak{P}|\mathfrak{p}}H^1(K_{n,\mathfrak{P}}, K_{\mathfrak{p}}(\varphi_{\mathfrak{p}})) & \times & \oplus_{\mathfrak{P}|\mathfrak{p}}H^1(K_{n,\mathfrak{P}}, \mathcal{O}_{\mathfrak{p}}(\varphi_{\mathfrak{p}})^*(1)) & \longrightarrow & K_{\mathfrak{p}} \\
\downarrow{\scriptstyle \text{res}} & & \uparrow{\scriptstyle S_n} & & \| \\
\oplus_{\mathfrak{P}|\mathfrak{p}}Hom(\mathbf{U}_{\infty,\mathfrak{P}}, K_{\mathfrak{p}}(\varphi_{\mathfrak{p}})) & \times & \oplus_{\mathfrak{P}|\mathfrak{p}}\mathbf{U}_{\infty,\mathfrak{P}} & \longrightarrow & K_{\mathfrak{p}}
\end{array}$$

Here the upper diagram commutes by property of restriction and corestriction map while the lower diagram commutes because of the definition of Soule map $S_n$ we used and class field theory. Notice the composite of left vertical maps sends $exp(\hat{\varpi})$ to $\delta$ and the composite of right vertical maps sends $u$ to $\sum_{\mathfrak{P}|\mathfrak{p}} S_0(u_{\mathfrak{P}})$. Hence by this diagram, we have

$$\begin{aligned}
\delta(u) &= \sum_{\mathfrak{P}|\mathfrak{p}} (res(exp(\hat{\varpi})), S_n(u_{\mathfrak{P}})) && \text{by lower diagram} \\
&= \sum_{\mathfrak{P}|\mathfrak{p}} (exp(\hat{\varpi}), cores(S_n(u_{\mathfrak{P}}))) && \text{by upper diagram} \\
&= \sum_{\mathfrak{P}|\mathfrak{p}} (\hat{\varpi}, exp^*(S_0(u_{\mathfrak{P}})) && \text{by def of } exp^* \text{ and } S_0(u_{\mathfrak{P}}) \\
&= \sum_{\mathfrak{P}|\mathfrak{p}} exp^*(S_0(u_{\mathfrak{P}}))/\varpi && \text{by the choice of } \varpi \text{ and } \hat{\varpi}
\end{aligned}$$

We use and generalize results of Kato on $exp^*$ in the following section.

## 3.3 Explicit reciprocity law

In this section we recall Kato's explicit reciporcity law for 1-dimensional Lubin-Tate formal groups, generalizing that of Wiles and also his earlier results with Bloch on exponential map for $\mathbf{Q}_p(r)$. All refrences in this section are made to [7].

Let $F/K$ be an unramified extension of $p$-adic local fields. Fix $\pi$ a prime element of $K$. Let $\Gamma$ be the Lubin-Tate formal group over $\mathcal{O}_K$ corrsponding to uniformizer $\pi$ and let $G$ be the connected p-divisible group over $\mathcal{O}_F$ obtained from $\Gamma$ by extension of scalars. Denote the action of $a \in \mathcal{O}_K$ on $G$ by [a]. Let $T$ be the Tate module for $G$ and fix a basis $\xi$ for T. Let $F_n$ denotes the field obtained by adjoining $\pi^n$ torsion point of $G$. For a norm compatible system of units $u = (u_n)_n \in \lim F_n^*$, define $S_n(u) \in H^1(F_n, T^{\otimes(-r)}(1))$, $\forall n \geq 0, r \geq 1$, exactly as we defined in definition 4 in the case $j = 0$ and $\mathfrak{p}$ is good reduction prime. Associated to such $u$ is its Coleman power series $g_{u,\xi}$ characterized by:

$$(\phi^{-n}(g_{u,\xi}))(\xi_n) = u_n, \forall n \geq 1$$

where $\phi$ is the Frobenius element in $Gal(F/K)$ acting on the coefficents of $g_{u,\xi}$.

**Theorem 1 (The explicit reciprocity law,[7] 2.1.7).** *$\forall n, r \geq 1$,*

$$exp^* : H^1(F_n, T^{\otimes(-r)}(1)) \to DR^0(F_n, V^{\otimes(-r)}(1)) \cong colie(G)^{\otimes(r)} \otimes_{\mathcal{O}_F} F_n$$

*sends $S_n(u)$ to*

$$\frac{1}{(r-1)!}\pi^{-nr}\omega^{\otimes r} \otimes \{(\frac{d}{\omega})^r \log(\phi^{-n}(g_{u,\xi}))\}(\xi_n) \tag{3.3.1}$$

Here $\omega$ is any $\mathcal{O}_F$ basis of $colie(G) = Hom_{\mathcal{O}_F}(Lie(G), \mathcal{O}_F)$.

*Remark 4.* 1. When $j = 0$ we are going to apply theorem 1 to $\Gamma = \hat{E}_{\mathfrak{p}}$

- If $\mathfrak{p}$ is a good reduction prime, then take $F = K = K_{\mathfrak{p}}$, $\Gamma = \hat{E}_{\mathfrak{p}}$, and $\pi = \psi(\mathfrak{p})$
- If $\mathfrak{p}$ is a bad reduction prime, take $K = K_{\mathfrak{p}}, F = F_{\mathfrak{P}}, \Gamma = \hat{E}'_{\mathfrak{p}}$ and $\pi = \psi'(\mathfrak{p})$

2. When $j > 0$ we can not directly apply theorem 1 as the proof given in [7] really require the twist to be 1. We will have to dig much deeper into the whole proof and extract a "$mod\pi^n$" version of $exp^*$ and find ways to pass from twist by $j+1$ to twist by 1. Finally we need a congruence argument to nail down $\delta(u)$.

3. Because of the nature of the proof, theorem 1 just describes explictly $exp^*$ on Soule-type elements $S_n(u)$ We will brieftly recall the details of the argument. Also the condition $n \geq 1$ is necessary, and the argument does not apply to $n = 0$.

For this purpose recall now some important steps for the proof of theorem 1. Readers who are just interested in verifying Tmamagawa number conjecture for $j = 0$ can go directly to the next section.

Let $\Omega^1(G)$ be the space of differential forms on $G$ and $\mathcal{O}(G) \cong \mathcal{O}_F[[t]]$ be the ring of functions on $G$. We have

$$\Omega^1(G) = \mathcal{O}(G) \otimes_{\mathcal{O}_F} colie(G)$$

For $n \geq 0$, let $S_n$ be the topological ring $p-$adically complete and separated, characterized by the property: $\forall i \geq 1$, $Spec(S_n/p^iS_n)$ is the $PD$ envelope of $Spec(\mathcal{O}_{F_n}/p^i\mathcal{O}_{F_n})$ in $Spec(\mathcal{O}(G)/p^i\mathcal{O}(G))$ with repsect to the embedding induced by $\xi_n : Spec(\mathcal{O}_{F_n}) \to G$. $\forall r \geq 1$ define $J_{S_n}^{[r]} \subset S_n$ to be the inverse limit over $i$ of the $r$-th divided power of $Ker(S_n/p^iS_n \to \mathcal{O}_{F_n}/p^i\mathcal{O}_{F_n})$. Since $G$ is formally smooth over $\mathcal{O}_F$, the fundamental theorem in crystalline cohomology theory shows

$$\begin{aligned} &R\Gamma(Spec(\mathcal{O}_{F_n}/\pi^i)/Spec(\mathcal{O}_F/\pi^i)_{\mathrm{crys}}, J^{[r]}) \\ &\cong [J_{S_n}^{[r]}/\pi^i \xrightarrow{d} J_{S_{n-1}}^{[r]}/\pi^i \otimes_{\mathcal{O}_G} \Omega^1(G)] \quad (i \geq 1) \end{aligned} \tag{3.3.2}$$

where $J = Ker(\mathcal{O}_{\mathrm{crys}} \to \mathcal{O}_{\mathrm{zar}})$. Following Kato, define a pairing

$$colie(G) \otimes_{\mathcal{O}_K} T \to J_\infty : \omega \otimes \xi \mapsto (l_{\omega,\xi,n})_n \tag{3.3.3}$$

Here for $\omega \in colie(G)$ , let $l_\omega$ be the logarithm of $G$ associated to $\omega$, i.e $dl_\omega = \omega$. Define $l_{\omega,\xi,n} = [\pi^n]^* l_\omega \in J_{S_n}$. Then $dl_{\omega,\xi,n} = \pi^n\omega$. In particular, $dl_{\omega,\xi,n} \equiv 0 \pmod{\pi^n}$. Hence by 3.3.2

$$l_{\omega,\xi,n} \bmod \pi^n \in H^0(Spec(\mathcal{O}_{F_n}/\pi^n)/Spec(\mathcal{O}_F/\pi^n)_{\mathrm{crys}}, J)$$

We denote the image of it in

$$H^0(Spec(\mathcal{O}_{\bar{F}}/\pi^n)/Spec(\mathcal{O}_F/\pi^n)_{\mathrm{crys}}, J) \cong J_n = J_\infty/\pi^n$$

also as $l_{\omega,\xi,n}$. This gives the above pairing which indues an isomorphism:

$$colie(G)^{\otimes r} \otimes_{\mathcal{O}_F} F \cong DR^0(F, V^{\otimes(-r)}(1) \otimes_{\mathbf{Q}_p} B_{\mathrm{dR}}^+) \tag{3.3.4}$$

We also need an integral version of this, as in [7] 2.2.3, that

$$colie(G)^{\otimes r} \otimes_{\mathcal{O}_K} \mathcal{O}_L \hookrightarrow H^0(L, T^{\otimes(-r)} \otimes_{\mathcal{O}_K} J_\infty^{[r]}/J_\infty^{[r+1]}) \tag{3.3.5}$$

where L is any finite extension of F.

*Remark 5.* We will generalize these statements when $j > 0$.

We introduce more notation and results from [7]. Define $\mu_{\xi,n} : S_n \otimes_{\mathcal{O}(G)} \Omega^1(G) \to H^1(F_n, J_\infty)$ to be the composite map coming from isomorphism 3.3.2 and spectral sequence

$$\begin{aligned} & R\Gamma(Spec(\mathcal{O}_{F_n}/\pi^i)/Spec(\mathcal{O}_F/\pi^i)_{\mathrm{crys}}, J) \\ & \to R\Gamma(Gal(\bar{F}/F_n), R\Gamma(Spec(\mathcal{O}_{\bar{F}}/\pi^i)/Spec(\mathcal{O}_F/\pi^i)_{\mathrm{crys}}, J) \\ & = R\Gamma(F_n, J_\infty/\pi^i) \end{aligned} \tag{3.3.6}$$

By taking $\cup\xi_n^{\otimes(-r)}$ this induces a map

$$S_n \otimes_{\mathcal{O}(G)} \Omega^1(G) \xrightarrow{\mu_{r,\xi,n}} H^1(F_n, (T^{\otimes(-r)} \otimes_{\mathcal{O}_K} J_\infty/J_\infty^{[r+1]})/\pi^n) \tag{3.3.7}$$

Let $\Theta = \lim(\ldots \to \Omega^1(G) \to \Omega^1(G) \to \Omega^1(G))$ where the arrows are trace map $Tr_\pi$ associated to $[\pi]$. Define

$$\theta_{r,\xi,n} : \Theta \to H^1(F_n, T^{\otimes(-r)} \otimes_{\mathcal{O}_K} J_\infty/J_\infty^{[r+1]}) \text{ to be the composite of}$$

$$\begin{aligned} \Theta = \lim_m \Omega^1(G) & \overset{(\mu_{r,\xi,m})_m}{\longrightarrow} \lim_m H^1(F_m, (T^{\otimes(-r)} \otimes_{\mathcal{O}_K} J_\infty/J_\infty^{[r+1]})/\pi^m) \\ & \overset{\mathrm{trace}}{\to} \lim_m H^1(F_n, (T^{\otimes(-r)} \otimes_{\mathcal{O}_K} J_\infty/J_\infty^{[r+1]})/\pi^m) \\ & \cong H^1(F_n, T^{\otimes(-r)} \otimes_{\mathcal{O}_K} J_\infty/J_\infty^{[r+1]}) \end{aligned} \tag{3.3.8}$$

The following commuative diagram is a key point in Kato's proof:

$$\begin{array}{ccccccc} \mathcal{O}(G)^* & \xrightarrow{\xi_n} & \mathcal{O}_{F_n}^* & \longrightarrow & H^1(F_n, \mathbf{Z}_p(1)) & \longrightarrow & H^1(F_n, T^{\otimes(-r)}(1)/\pi^n) \\ \downarrow{\scriptstyle d\log} & & & & \downarrow & & \downarrow \\ \Omega^1(G) & & \xrightarrow{-\mu_{\xi,n}} & & H^1(F_n, J_\infty) & \longrightarrow & H^1(F_n, (T^{\otimes(-r)} \otimes_{\mathcal{O}_F} J_\infty)/\pi^n) \end{array} \tag{3.3.9}$$

which reduces the proof of theorem 1 to the following

**Theorem 2 ([7] 2.2.7).** *$\forall n, r \geq 1$, the composite map*

$$\begin{aligned} \Theta & \overset{\theta_{r,\xi,n}}{\longrightarrow} H^1(F_n, T^{\otimes(-r)} \otimes_{\mathcal{O}_K} J_\infty/J_\infty^{[r+1]}) \\ & \to H^1(F_n, V^{\otimes(-r)} \otimes_K B_{dR}^1/B_{dR}^{r+1}) \cong colie(G)^{\otimes r} \otimes_{\mathcal{O}_F} F_n \end{aligned} \tag{3.3.10}$$

*sends $(\eta_m)_m \in \Theta$ to*

$$-\frac{1}{(r-1)!}\pi^{-nr}\omega^{\otimes r} \otimes \{(\frac{d}{\omega})^{r-1}(\frac{\eta_n}{\omega})\}(\xi_n)$$

Actually Kato deduced this from "finite version" of it as follows:

**Proposition 2 ([7] 2.3.2).** *For fixed $n \geq 1, \exists$ integer $c = c(F,n) \neq 0$, such that $\forall m \geq n$, let $\alpha_m$ be defined as follows:*

$$\begin{array}{ccc} S_m \otimes \Omega^1(G) & \xrightarrow{\mu_{r,\xi,m}} & H^1(F_m, (T^{\otimes(-r)} \otimes_{\mathcal{O}_K} J_\infty/J_\infty^{[r+1]})/\pi^m) \\ & \searrow^{\alpha_m} & \downarrow trace \\ & & H^1(F_n, (T^{\otimes(-r)} \otimes_{\mathcal{O}_K} J_\infty/J_\infty^{[r+1]})/\pi^m) \end{array}$$

*and set $\alpha_{c,m} = c\alpha_m$. Define $\beta_{c,m}$ to be the map :*

$$\beta_{c,m}(\eta) = ((r-1)!^{-1}\pi^{-nr}c)\omega^{\otimes r} \otimes \{(\frac{d}{\omega})^{r-1}(\frac{1}{\omega}Trace_{m,n}(\eta))\}(\xi_n) \cup \log(\chi_G)$$

*Then $\alpha_{c,m} + \beta_{c,m} = 0$.*

Here $c(F,n)$ is independent of $m \geq n$ and $\chi_G$ describes the Galois action on $T$, but notice $\log(\chi(G))$ and $\log(\chi_{\mathrm{cyclo}})$have the same image in $H^1(F, \mathbb{C}_p)$. We do not have to worry about the definition of "Trace" map here, as in our application, we take $(\eta_m)_m \in \Omega^1(G)$ the element given by $(d\log\phi^{-m}(g_{u,\xi}))_m$ and it follows from the norm compatibility of $u = (u_m)_m$ that after taking $d\log$ this particular element is trace compatible.

**Lemma 5.** *We have the following commutative diagram*

$$\begin{array}{ccc} H^1(F_n, T^{\otimes(-r)}(1)/\pi^m) & \xrightarrow{c(n)\iota_1} & H^1(F_n, (T^{\otimes(-r)} \otimes_{\mathcal{O}_K} J_\infty/J_\infty^{[r+1]})/\pi^m) \\ \vdots\lambda & \nearrow^{\rho} & \uparrow \log\chi_G \\ colie(G)^{\otimes r} \otimes_{\mathcal{O}_K} \mathcal{O}_{F_n}/\pi^m & \xrightarrow{\iota_2} & H^0(F_n, (T^{\otimes(-r)} \otimes_{\mathcal{O}_K} J_\infty^{[r]}/J_\infty^{[r+1]})/\pi^m) \end{array}$$

*where $\iota_1$ on the first row is induced by inclusion $\mathbf{Z}_p(1) \to J_\infty$ and $\iota_2$ is induced from 3.3.5.*

*Proof.* This follows from proposition 2 and 3.3.9. The point is even though in general one might not be able to define $\lambda$ on the left side of the diagram, but we do know that in $H^1(F_n, (T^{\otimes(-r)} \otimes_{\mathcal{O}_K} J_\infty/J_\infty^{[r+1]})/\pi^m)$,

$$\rho(\frac{c(n)}{(r-1)!}\pi^{-nr}\omega^{\otimes r} \otimes \{(\frac{d}{\omega})^r \log(\phi^{-n}(g_{u,\xi}))\}(\xi_n)) = c(n)\iota_1(\overline{(S_{n,m}}(u)), \qquad (3.3.11)$$

□

*Remark 6.* We are going to apply proposition2 and 3.3.11 in the following situation: $j > 0$ (see also remark 4)

1. when $p$ splits,

$$(F,E) \stackrel{\text{def}}{=} \begin{cases} (K,E) \\ (F,E') \end{cases} \qquad n_0 \stackrel{\text{def}}{=} \begin{cases} 1, \text{ if } \mathfrak{p} \text{ is a good reduction prime} \\ m(\mathfrak{p}), \text{ if } \mathfrak{p} \text{ is a bad reduction prime} \end{cases} \tag{3.3.12}$$

   apply with base field $F = F(E_{\mathfrak{p}^{*t}})_{\mathfrak{P}}$, $\mathfrak{P} \mid \mathfrak{p}, t \to \infty$ and $n = n_0$ which is fixed. Denote $c(t)$ to be the constant needed when applying proposition 2 to the field $F(E_{\mathfrak{p}^{n_0}\mathfrak{p}^{*t}})_{\mathfrak{P}}$.

2. when $p$ is nonsplit, define $(F,E)$ just as above. Apply with $F = F_{\mathfrak{P}}$ and $t$ any positive integer. Denote $c(t)$ to be the constant needed when applying proposition 2 to the field $F(E_{\mathfrak{p}^t})_{\mathfrak{P}}$.

Note in all cases we have a fixed Lubin-Tate formal group when $t$ changes.

**Theorem 3.** *In the case p splits, the constant $c(t)$ in proposition 2 can be chosen independent of $t$.*

*Proof.* We need to study more carefully how the constant $c(F,n)$ in proposition 2 depends on $F,n$. Note in our situation, $n = n_0$ is fixed while $t$ changes (hence base field changes). Recall Kato has put several restriction on $c$:

(1) In lemma 2.3.3, he required $(r-1)!^{-1}\pi^{-nr}c \in \mathcal{O}_F$, and $\alpha_{c,m} + \beta_{c,m}$ kills $l^{r-1}_{\omega,\xi,m}S_m \otimes_{\mathcal{O}(G)} \Omega^1(G), \forall m \geq n$. He deduced this from a explicit description of $\mu_{r,\xi,m}$ on $l^{r-1}_{\omega,\xi,m}S_m \otimes \Omega^1(G)$. In 2.3.8 we see all he need is that c is big enough so that

$$c\pi^{-m}\log(\chi_G(\sigma)) \equiv c\pi^{-m}(\chi_G(\sigma)-1) \mod \pi^m, \forall m \geq n, \sigma \in Gal(\overline{\mathbf{Q}_p}/F_m)$$

It is clear one can choose such an integer c independent of m and even n (take c corresponding to $n = 0$). Then by Kato's calculation,

$$\alpha_{c,m}(l^{r-1}_{\omega,\xi,m}h \otimes \omega) = -\omega^{\otimes r} \otimes (\pi^{n-m}Tr_{m,n}(h(\xi_m)))(c\pi^{-n}\log(\chi_G))$$

$\forall h \in S_m$. The points is $Tr_{m,n} : \mathcal{O}_{F_m} \to \mathcal{O}_{F_n}$ has image contained in $\pi^{m-n}\mathcal{O}_{F_n}$. (This one of the place where the argument fails if one take $n = 0$.) On the other hand, a direct computation shows

$$\{(\frac{d}{\omega})^{r-1}(\frac{1}{\omega}Tr_{m,n}(l^{r-1}_{\omega,\xi,m}h \otimes \omega))\}(\xi_n) = (r-1)!\pi^{n(r-1)}\pi^{n-m}Tr_{m,n}(h(\xi_m))$$

This shows

$$\beta_{c,m}(l^{r-1}_{\omega,\xi,m}h \otimes \omega) = -\omega^{\otimes r} \otimes (\pi^{n-m}Tr_{m,n}(h(\xi_m)))(c\pi^{-n}\log(\chi_G))$$

and in particular we do not need to require that $(r-1)!^{-1}\pi^{-nr}c \in \mathcal{O}_F$ here.

(2) In lemma 2.3.10, Kato required c to kill $J^{[r-1]}_{S_m}/(J^{[r]}_{S_m} + l^{r-1}_{\omega,\xi,m}S_m), \forall m \geq n$ which he remarked by lemma 2.3.11 it is enough to ask c annihilates $J_{S_m}/(J^{[2]}_{S_m} + l_{\omega,\xi,m}S_m)$. Kato was able to find such an integer c which is independent of m and n. see the end of Page 141 of [7].

(3) In lemma 2.3.9, he wanted an integer c so that $\alpha_{c,m}+\beta_{c,m}$ kills $J^{[r-1]}_{S_m}\otimes_{\mathcal{O}(G)} \Omega^1(G)\forall m \geq n$. An interger c which satisfies both of the above condition will automatically satisfies this, since it follows from the definition that $\alpha_{c,m} + \beta_{c,m}$ vanishes on $J^{[r]}_{S_m} \otimes_{\mathcal{O}(G)} \Omega^1(G)$ .

(4) Finally a condtion on c in lemma 2.3.13 was reduced to condition that c annihilates group $H^q(F_n, \hat{\mathcal{O}}_{\bar{F}}(-r))$ $(q = 0, 1)$ in lemma 2.3.16.

Hence to prove theorem 3, we just need to prove:

**Proposition 3.** *Fix integer $r \neq 0$. Let $K$ be a finite extension of $\mathbf{Q}_p$ in (a)-(b).*
*(a): $H^0(K, \mathcal{O}_{\mathbb{C}_p}(r)) = 0$.*
*(b): $H^1(K, \mathcal{O}_{\mathbb{C}_p}(r))$ is killed by a power of $p$ that depends on $r$ and $K$.*
*(c): When $\mathfrak{p} \mid p$ splits in $K/\mathbf{Q}$, fix $n_0 \geq 1$, then*

$$p^{d_r}H^1(K_{n_0,m,\upsilon}, \mathcal{O}_{\mathbb{C}_p}(r))) = 0$$

*for some constant $d_r$ that is independent of $m$, but depends on $n_0, r$. Here $K$ is our imaginary quadratic field.*

We use and generalize some results of Tate [13] and Sen [11]. Note we are dealing with continuous Galois cohomology, but the coefficent module $\mathcal{O}_{\mathbb{C}_p}(r)$ is not a "discrete" Galois module, hence usual method to reduce it to computation of cohomology of finite groups fails.

*Proof.* First recall some results of Tate in [13]: Let $K$ be the quotient field of a complete DVR of char 0, let $C$ be the completion of algebraic closure of $K$; let $K_\infty$ be any infinite totally ramified extension of $K$ such that $Gal(K_\infty/K) \cong \mathbf{Z}_p$; let $X$ be the completion of $K_\infty$; let $\chi$ be a continuous character of $\Gamma = Gal(K_\infty/K)$ into the group of units of $K$. Then

1. Prop 8(b) [13]: if $\chi(\Gamma)$ is infinite, then $H^0(\Gamma, X(\chi))$ and $H^1(\Gamma, X(\chi))$ are 0.

2. Prop 10 [13] $H^0(Gal(\bar{K}/K_\infty), C) = X$, $H^r(Gal(\bar{K}/K_\infty), C) = 0, \forall r > 0$

Using these tow proposition, we get

$$H^0(K, \mathcal{O}_{\mathbb{C}_p}(r)) \subseteq H^0(K, \mathbb{C}_p(r)) = H^0(\Gamma, H^0(Gal(\bar{\mathbf{Q}}_p/K_\infty, \mathbb{C}_p(r)))$$
$$= H^0(\Gamma, X(r)) = 0$$

where (1) follows from prop 10 and (2) follows from prop 8(b) of [13]. This proves claim (a).

(b) is stated without proof in [7](2.2.4). We have by inflation-restriction exact sequence:

$$0 \to H^1(\Gamma, \mathcal{O}_X(r)) \to H^1(K, \mathcal{O}_{\mathbb{C}_p}(r)) \to H^1(K_\infty, \mathcal{O}_{\mathbb{C}_p}(r)) \tag{3.3.13}$$

First we prove $p$ annihilates $H^1(K_\infty, \mathcal{O}_{\mathbb{C}_p}(r)) \cong H^1(K_\infty, \mathcal{O}_{\mathbb{C}_p})$. Since for every continuous cochain $f$ on $Gal(\bar{\mathbf{Q}}_p/K_\infty)$ with values in $\mathcal{O}_{\mathbb{C}_p}$(taken with valuation topology), we can find a sequence of continous cochain $f_n$ with values in $\mathcal{O}_{\bar{\mathbf{Q}}_p}$(taken with discrete topology) such that $|f - f_n| \to 0$(see [13] 3.2 prop10), we just need to prove that $p$ kills $H^1(K_\infty, \mathcal{O}_{\bar{\mathbf{Q}}_p})$. This can also be seen in the following way: we have commutative diagram

$$\begin{array}{ccccc} \mathcal{O}_{\bar{\mathbf{Q}}_p} & \longrightarrow & \bar{\mathbf{Q}}_p & \longrightarrow & \bar{\mathbf{Q}}_p/\mathcal{O}_{\bar{\mathbf{Q}}_p} \\ \downarrow & & \downarrow & & \| \\ \mathcal{O}_{\mathbb{C}_p} & \longrightarrow & \mathbb{C}_p & \longrightarrow & \mathbb{C}_p/\mathcal{O}_{\mathbb{C}_p} \end{array}$$

where we put discrete topology on the first row and valuation topology on the second row and the vertical maps are continuous. Since by [13] prop 10, $H^1(K_\infty, \mathcal{O}_{\bar{\mathbf{Q}}_p})$ is torsion group, we have another commutative diagram where horizontal maps are epimorphisms:

$$\begin{array}{ccc} H^0(K_\infty, \bar{\mathbf{Q}}_p/\mathcal{O}_{\bar{\mathbf{Q}}_p}) & \longrightarrow & H^1(K_\infty, \mathcal{O}_{\bar{\mathbf{Q}}_p}) \\ \| & & \downarrow \\ H^0(K_\infty, \mathbb{C}_p/\mathcal{O}_{\mathbb{C}_p}(r)) & \longrightarrow & H^1(K_\infty, \mathcal{O}_{\mathbb{C}_p}) \end{array}$$

By the way exactly the same argument can be applied to show that if some power of $p$ kills $H^1(\Gamma, \mathcal{O}_{K_\infty}(r))$, then the same power of $p$ kills $H^1(\Gamma, \mathcal{O}_X(r))$ too.

By a well known result of Serre ([12] prop 8) that every continuous cochain in $Gal(\bar{\mathbf{Q}}_p/K_\infty)$ with values in $\mathcal{O}_{\bar{\mathbf{Q}}_p}$ comes by inflation from some finite Galois extension $M/K_\infty$, we see $pH^1(K_\infty, \mathcal{O}_{\mathbb{C}_p}) = 0$ follows immediately from the following generalization of [13] 3.2, cor1:

**Lemma 6.** *Let $K_\infty$ be a deeply ramified field. Let $M/K_\infty$ be a finite Galois extension with group $G$. Let $f$ be an $n-$cochain of $G$ with coefficents in $M(\chi)$ where $\chi$ is any character of $G$, $n \geq 0$. Let $c > 1$. Then there exists an $(n-1)-$cochain $g$ of $G$ in $M(\chi)$ such that*

$$|f - \delta g| \leq c|\delta f|, \textit{and} \quad |g| \leq c|f|$$

*In particular it follows $H^q(G, \mathcal{O}_M(\chi))$ is killed by $p$, $\forall q > 0$. Here $|f|$ denotes the maximum of the absolute values of the coefficents of $f$ and by a $(-1)-$cochain we mean an element $y \in M(\chi)$. The coboundary $\delta y$ of such a $y$ is the $0-$cochain $Tr_{M/K_\infty}$.*

*Proof.* For more details about deeply ramified extensions see [2]. They are generalization of totally ramified $\mathbf{Z}_p$ extension studied by Tate[13]. There are 5 equivalent characterization of it one of which says if $K_\infty$ is deeply ramified, then

$$Tr_{M/K_\infty}(\mathfrak{m}_M) = \mathfrak{m}_\infty, \quad \forall M \text{ finite extension of } K_\infty$$

where $\mathfrak{m}_M$ and $\mathfrak{m}_\infty$ are maxiaml ideal of $M$ and $L_\infty$ respectively. Hence there exists a $(-1)-$cochain $y \in M(\chi)$ such that $|y| \leq 1$ and $|\delta y| \geq c^{-1}$. Define an $(n-1)-$cochain $y \cup f$ by the formulas:

$$y \cup f = yf, \quad \text{if } n = 0$$

$$(y \cup f)(s_1, \ldots, s_{n-1}) = (-1)^n \sum_{s_n \in G} s_1 s_2 \cdots s_n y f(s_1, \ldots, s_n), \quad \text{if } n > 0$$

One checks easily identity $(\delta y)f - \delta(y \cup f) = y \cup (\delta f)$ which holds formally and really has nothing to do with what $\chi$ is. Now divide last equality by $x = \delta y = Tr_{M/K_\infty} y \in K_\infty$ we find

$$f - \delta g = x^{-1}(y \cup \delta f), \quad \text{with } g = x^{-1}(y \cup f)$$

Since $|x^{-1}| \leq c$, and $|y| \leq 1$, we conclude the first part of the lemma.
Now if $f \in H^q(G, \mathcal{O}_M(\chi))$, then we can find an $(q-1)-$cochain $g$ of $G$ in $M(\chi)$, such that

$$f = \delta g \quad \text{and } |g| \leq c$$

We can pick e.g $c = 2$. Then we see $pf = \delta(pg)$ with $pg \in H^{q-1}(G, \mathcal{O}_M(\chi))$. Hence $pf = 0$ in $H^q(G, \mathcal{O}_M(\chi))$. $\square$

Secondly we prove $H^1(\Gamma, \mathcal{O}_{K_\infty}(r))$ is killed by a power of $p$. By 3.3.13, this will give (b) of proposition 3. Denote $K_n$ as the cyclic subextension of $K$ inside $K_\infty$ of degree $p^n$. Let $\sigma$ be a generator of $\Gamma = Gal(K_\infty/K)$ and so $\Gamma_n \stackrel{\text{def}}{=} Gal(K_\infty/K_n)$ is generated by $\sigma^{p^n}$. Let $\lambda = \chi_{\text{cyclo}}^{-r}(\sigma)$ which is $\equiv 1 \pmod{\pi}$ where $\pi$ is a uniformizer of $K$. Tate has defined a continuous operator $t$ (an idempotent) on $X$(see [13] prop 5,6,7) such that $X$ is a direct sum of $K$ and $X_0 = ker(t)$ and $\sigma - 1$ is bijective with a continuous inverse $\rho$ on $X_0$. Moreover he proved

$$|\rho y| \leq d|y|, \quad \forall y \in X_0$$

where $d$ is a constant depending only on the $\mathbf{Z}_p$ extension in question(see Tate's remark after prop 6 in [13]). In particular if we replace $K$ by $K_n$ as ground field, we get $\sigma^{p^n} - 1$ is bijective with a continuous inverse $\rho_n$ on $X_0$ and

$$|\rho_n y| \leq d|y|, \quad \forall y \in X_0$$

with the same $d$.

**Lemma 7.** *There exists integer $a(K)$ which depends on $K$ such that $t(\mathcal{O}_{K_\infty}) \subseteq p^{a(K)}\mathcal{O}_K$ and we have the following exact sequence*

$$0 \to \mathcal{O}_{K_{\infty,0}} \to \mathcal{O}_{K_\infty} \xrightarrow{t} p^{a(K)}\mathcal{O}_K \tag{3.3.14}$$

*where $\mathcal{O}_{K_{\infty,0}} = \mathcal{O}_{K_\infty} \cap X_0$.*

*Proof.* We need a known result which follows easily from the definition of different:

**Lemma 8.** *Let $M/F$ be extensions of local fields. Then*

$$Trace(\mathcal{O}_M) = \pi_F^{[ord_F(\delta_{M/F})]}\mathcal{O}_F$$

*where $\delta_{M/F}$ is the relative different; $\pi_F$ is a uniformizer of $F$ and $ord_F$ is the valuation on $F$ normalized so that $ord_F(\pi_F) = 1$.*

Applying this lemma to $K_n/K$, we get

$$tr(\mathcal{O}_{K_n}) = p^{-n}Trace(\mathcal{O}_{K_n}) \overset{(1)}{=} p^{-n}\pi_K^{[e_K n + c_K]}\mathcal{O}_K \overset{(2)}{=} \pi_K^{[c_K]}\mathcal{O}_K \tag{3.3.15}$$

where (1) follows from [13] proposition 5 which use the classical formula expressing relative different in terms of higher ramification groups. $e_K$ is the absolute ramification index of $K$, i.e $ord_K(p) = e_K$. $c_K$ is a constant. Let $a(K) = [c_K/e_K]$, this gives the claim in the lemma 7. ∎

We have exact sequence:

$$0 \to H^1(\Gamma, \mathcal{O}_{K_{\infty,0}}(r)) \to H^1(\Gamma, \mathcal{O}_{K_\infty}(r)) \to p^{a(K)} H^1(\Gamma, \mathcal{O}_K(r)) \tag{3.3.16}$$

We have by inflation-restriction sequence,

$$0 \to H^1(\Gamma/\Gamma_n, \mathcal{O}_{K_{\infty,0}}(r)^{\Gamma_n}) \to H^1(\Gamma, \mathcal{O}_{K_{\infty,0}}(r)) \hookrightarrow H^1(\Gamma_n, \mathcal{O}_{K_{\infty,0}}(r)) \tag{3.3.17}$$

This is because by [13] prop 8(b), $\mathcal{O}_{K_{\infty,0}}(r)^{\Gamma_n} \subset X(r)^{\Gamma_n} = 0$.

$$H^1(\Gamma_n, \mathcal{O}_{K_{\infty,0}}(r)) = \mathcal{O}_{K_{\infty,0}}(r)/1 - \sigma^{p^n} \cong \mathcal{O}_{K_{\infty,0}}/\lambda - \sigma^{p^n} \tag{3.3.18}$$

We have

$$\sigma^{p^n} - \lambda^{p^n} = (\sigma^{p^n} - 1) - (\lambda^{p^n} - 1) = (\sigma^{p^n} - 1)(1 - (\lambda^{p^n} - 1)\rho_n) \tag{3.3.19}$$

On the other hand, by Sen's result, $H^1(K_n, \mathcal{O}_{K_\infty})$ is killed by $p$ for all $n$, which implies

$$\forall x \in \mathcal{O}_{K_{\infty,0}}, \exists\, y \in \mathcal{O}_{K_\infty}, \text{ such that } px = (\sigma^{p^n} - 1)y$$

Take any $n$ big enough so that $|(\lambda^{p^n} - 1)d| < 1$ and consequently $1 - (\lambda^{p^n} - 1)\rho_n$ is an automorphism on $X_0$. Let $z = \{1 - (\lambda^{p^n} - 1)\rho_n\}^{-1} y$, then by 3.3.19 we have

$$px = (\sigma^{p^n} - \lambda^{p^n})z \tag{3.3.20}$$

Notice $|z| = |z - (\lambda^{p^n} - 1)\rho_n z| = |y| \leq 1$, hence $z \in \mathcal{O}_{K_\infty}$. It is also clear that $z \in X_0$ as $1 - (\lambda^{p^n} - 1)\rho_n$ is an automorphism on $X_0$. Hence $z \in \mathcal{O}_{K_{\infty,0}}$. We conclude from 3.3.18 and 3.3.20 that $pH^1(\Gamma_n, \mathcal{O}_{K_{\infty,0}}(r)) = 0$ for all $n$ sufficently big and hence by 3.3.17, $pH^1(\Gamma, \mathcal{O}_{K_{\infty,0}}(r)) = 0$. Notice $H^1(\Gamma, \mathcal{O}_K(r))$ is killed by $p^{d(K)+ord_p(r)}$ where $d(K)$ is the largest integer $n$ such that $K$ contains all $p^n$ th roots of unity. Now by 3.3.16, $p^b H^1(\Gamma, \mathcal{O}_{K_\infty}(r)) = 0$, for some $b$ that depends on $r$ and $K$. This gives (b) in the proposition 3.

Finally to prove (c), we just remark the constant $a(K)$ in lemma 7 depends only on the higher ramification groups of $Gal(K_n/K), \forall n$. If one consider $L/K$ finite extension disjoint from cyclotomic extension of $K$ (in particular this is true if $L/K$ is unramified), then it is easy to see $a(K) = a(L)$. On the other hand, $d(K)$ also stays the same as $K$ runs through the fields $K_{n_0,m,v}, \forall m$, which has fixed ramification. Since the only dependence of the power of $p$ needed to kill $H^1(K, \mathcal{O}_{\mathbb{C}_p}(r))$ in (c) is in $a(K)$ and $d(K)$, we see (c) follows. □

*Remark 7.* We had believed earlier that $H^1(K, \mathcal{O}_{\mathbb{C}_p}(r))$ is killed by a power of $p$ which only depends on $r$, by the influence of result of Tate and Sen. In order to prove this, it is clear one just need to elaborate on the dependence of $a(K)$ in lemma 7 on $K$, especially when one replace $K$ by $K_n$. The crucial point is even though $X = X_0 \oplus K$, but $t$ does not respect integral structure and $H^1(\Gamma, \mathcal{O}_K(r))$ (which cause us all the trouble) is not a direct summand of $H^1(\Gamma, \mathcal{O}_X(r))$. One wants to prove the image of $H^1(\Gamma, \mathcal{O}_K(r)) \to H^1(\Gamma, \mathcal{O}_X(r))$ is small. We have trouble proving this now. Hence we do not have the corresponding result of (c) in proposition 3 for $p$ is nonsplit case yet. This will force us to modify somewhat the calculation of $exp^*$ in the following when $p$ is nonsplit.

A classical result of Sen([11], theorem 3) says $H^1(K, \mathcal{O}_{\bar{\mathbf{Q}}_p})$ is killed by $p$, for all local field $K$. Recall Sen proved this in [11] by a careful study of wildly ramified automorphism of local fields, in particular

**Lemma 9 ([11] Theorem 2 ).** *Let $M/F$ be a Galois extension of local fields whose residue field has char $p$ and $G = Gal(M/F)$ is a cyclic group of order $p^n$ generated by $\sigma$. Let $\pi_F$(resp. $\pi_M$) be a uniformizer of $F$(resp.$M$). Then*

$$H^1(Gal(M/F), \mathcal{O}_M) = \oplus_{\mu=1}^{\mu=p^n-1} \mathcal{O}_F / \pi_F^{(\mu)} \mathcal{O}_F$$

*where $(\mu)$ is the greatest integer less than $(\mu + i(\mu))/p^n$ and $i(\mu) = ord_M((\sigma^\mu - 1)\pi_M/\pi_M)$, $ord_M$ is the discrete valuation on $M$ such that $ord_M(\pi_M) = 1$. In particular $H^1(Gal(M/F), \mathcal{O}_M)$ is killed by $p$.*

We will need the following 2 variants of it later:

**Lemma 10.** *Let $M/F$ be an abelian $p$ extension of local fields whose residue field extension has characteristic $p$. Then $H^1(Gal(M/F), \mathcal{O}_M)$ is killed by $p^e$, where the constant $e$ is the number of cyclic components in the decomposition of $Gal(M/F)$.*

*Proof.* By easy induction on the number of cyclic components, the lemma is reduced to the special case $Gal(M/F) = C_1 \times C_2$ where $C_1$ and $C_2$ are cyclic $p$ groups. Let $N \subset M$ be the field corresponding to $C_2$, i.e $Gal(M/N) = C_2$. Applying the spectral sequence associated to

$$0 \to C_2 \to Gal(M/F) \to C_1 \to 0$$

we get the following exact sequence:

$$H^1(C_1, \mathcal{O}_M^{C_1}) \to H^1(Gal(M/F), \mathcal{O}_M) \to H^1(C_2, \mathcal{O}_M)^{C_1}$$

Since $\mathcal{O}_M^{C_1} = \mathcal{O}_N$, and by lemma 9 (applying to $M/N$ and $N/F$ respectively), we get

$$H^1(C_1, \mathcal{O}_M^{C_1}) = 0, \quad H^1(C_2, \mathcal{O}_M) = 0$$

hence the lemma follows. □

**Lemma 11.** *With the same assumption as in lemma 9, let $\chi$ be a nontrivial character of $G = Gal(M/F)$. Then $H^1(G, \mathcal{O}_M/p^n(\chi))$ is killed by a power of $p$ which depends only on $F$ and $\chi$.*

*Proof.* We follow some argument of Sen. Recall that $H^1(G, \mathcal{O}_M(\chi)/p^n) = A/(\sigma - 1)\mathcal{O}_M(\chi)/p^n$ where $\sigma \in G$ is a generator and

$$A \stackrel{\text{def}}{=} ker(\sum_{i=0}^{p^n-1} \sigma^i : \mathcal{O}_M(\chi)/p^n \to \mathcal{O}_M(\chi)/p^n)$$

Sen has found $x_\mu \in M, 1 \leq \mu \leq p^n - 1$, such that the $x_\mu' s$ and 1 together form an $\mathcal{O}_F$ basis of $\mathcal{O}_M$. Let $y_\mu = (\sigma - 1)x_\mu$. It is known ([11] theorem 1) that $ord_M(y_\mu)$ are distinct mod $p^n$. It follows that $\{y_\mu, \mu = 1, \dots, p^n - 1\}$ are $\mathcal{O}_F$ basis for $(\sigma - 1)\mathcal{O}_M$. For any character $\chi$ of $G$, it is clear that $\bar{1}, \overline{y_\mu/\pi_F^{(\mu)}}$ form an $\mathcal{O}_F/p^n$ basis of A. Hence the claim follows from lemma 9 upon observing that $\mathcal{O}_F/p^n/(1 - \chi(\sigma))$ (for $\chi$ nontrivial) is killed by a power depending only on $F$ and $\chi$. □

The difference between lemma 11 and lemma 9 is now there is contribution from $\mathcal{O}_F$.

# CHAPTER 4

# COMPUTATION OF DUAL EXPONENTIAL MAP WHEN $J = 0$

## 4.1 Elliptic units and Eisenstein series

Here we will define exactly the elliptic unit $u \in \mathfrak{A}_\infty$ we are going to use. The bridge between its Coleman power series and Hecke $L$-series in provided by Eisenstein series. We will recall some important results we need for our calculation. Of special interest is a congruence property between Eisenstein numbers (which can be thought of as special values of Eisenstein series at CM points of modular curves) and congruences with $p$-adic periods of an elliptic curve. We need these congruences to link $L(\bar{\psi^k}, k)$ with $L(\bar{\psi^{k+j}}, k)$ and relate $exp^*$ map for $\mathcal{M}_{k,j}$ with that for $\mathcal{M}_{k,0}$ in the next section. Our main refrence for this section is [3].

Let $\theta$ be the fundamental theta function (see [3] p48 (7))

$$\theta(z, L) = \Delta(L)e^{-6\eta(z,L)z}\sigma(z, L)^{12}$$

and for $\boldsymbol{\alpha}$ an integral ideal of $K$, define $\Theta(z, L, \boldsymbol{\alpha}) = \frac{\theta(z,L)^{N\boldsymbol{\alpha}}}{\theta(z,\boldsymbol{\alpha}^{-1}L)}$ which is an elliptic function with respect to the second variable $L$.

$$\Theta(z, L, \boldsymbol{\alpha}) = \frac{\Delta(L)}{\Delta(\boldsymbol{\alpha}^{-1}L)} \prod_{u\in\boldsymbol{\alpha}^{-1}L/L}{}' \frac{\Delta(L)}{(\mathcal{P}(z, L) - \mathcal{P}(u, L))^6}$$

Now let us assume moreover that $L$ has complex multiplication. Here are some important results we will use: (see[3] prop 2.3)

1. Let $\mathfrak{m}$ be a non-trivial integral ideal of $K$, $v$ a primitive $\mathfrak{m}$ division point of $L$. Assume $(\boldsymbol{\alpha}, \mathfrak{m}) = 1$. Then

   - $\Theta(v, L, \boldsymbol{\alpha}) \in K(\mathfrak{m})$ : it is a unit if $\mathfrak{m}$ is not a prime power

- $\Theta(\upsilon, L, \boldsymbol{\alpha})^{\sigma_{\mathfrak{c}}} = \Theta(\upsilon, \mathfrak{c}^{-1}L, \boldsymbol{\alpha}) = \Theta(\psi(\mathfrak{c})\upsilon, L, \boldsymbol{\alpha})$ where $\mathfrak{c}$ is any integral ideal prime to $\mathfrak{m}$,$\sigma_{\mathfrak{c}}$ is the image of $\mathfrak{c}$ under artin map and $\psi = \psi_E$ for elliptic curve $E$ corresponding to $L$.

2. Let $\mathfrak{a}$ and $\mathfrak{b}$ be two integrals of $K$, prime to each other. then there is distribution relation:

$$\prod_{\upsilon \in \mathfrak{b}^{-1}L/L} \Theta(z + \upsilon, L, \boldsymbol{\alpha}) = \Theta(z, \mathfrak{b}^{-1}L, \boldsymbol{\alpha}) \tag{4.1.1}$$

Define elliptic unit $u \in \mathfrak{A}_\infty$ as follows:

1 (a): if $p$ is a good split reduction prime, set $\Omega = \Omega_\infty/f$ where $(f) = \mathfrak{f}$ is the conductor of $E$(or $\psi$), where $L = \Omega_\infty \mathcal{O}_K$ the period lattice of our CM elliptic curve $E$. Define

$$e_{n,m} = \Theta(\Omega, \mathfrak{p}^n \mathfrak{p}^{*m} L, \boldsymbol{\alpha}) = \Theta(\frac{\Omega}{\pi^n \pi^{*m}}, L, \boldsymbol{\alpha}) \in K(\mathfrak{f}\mathfrak{p}^n\mathfrak{p}^{*m})^*$$

Define

$$u_{n,m} = Norm_{K(\mathfrak{f}\mathfrak{p}^n\mathfrak{p}^{*m})/K_{n,m}}(e_{n,m})$$

From properties above we see they are units in $K_{n,m}$.

**Lemma 12.** *$(u_{n,m})_{n,m}$ are norm compatible.*

*Remark 8.* A special case of this for $e_{n,0}$ is stated without proof in [3].

*Proof.* $\forall n_1 \geq n_2, m_1 \geq m_2$,

$$Norm_{K(\mathfrak{f}\mathfrak{p}^{n_1}\mathfrak{p}^{*m_1})/K(\mathfrak{f}\mathfrak{p}^{n_2}\mathfrak{p}^{*m_2})}(e_{n_1,m_1}) = \prod_{\mathfrak{c} \in B} \Theta(\psi(\mathfrak{c})\frac{\Omega}{\pi^{n_1}\pi^{*m_1}}, L, \boldsymbol{\alpha}) \tag{4.1.2}$$

by property of Galois action, where $B$ is a set of integral ideals of K such that $\{(\mathfrak{b}, K(\mathfrak{f}\mathfrak{p}^{n_1}\mathfrak{p}^{*m_1})/K) : \mathfrak{b} \in B\}$ describes precisely $Gal(K(\mathfrak{f}\mathfrak{p}^{n_1}\mathfrak{p}^{*m_1})/K(\mathfrak{f}\mathfrak{p}^{n_2}\mathfrak{p}^{*m_2}))$. Since conductor of $(\psi) = \mathfrak{f}$, by definition of set $B$, we have

$$\psi(\mathfrak{c}) \equiv 1 \mod \mathfrak{f}\mathfrak{p}^{n_2}\mathfrak{p}^{*m_2}$$

Hence as $\mathfrak{c}$ runs through set $B$, $(\psi(\mathfrak{c})-1)\Omega/\pi^{n_1}\pi^{*m_1}$ exactly runs through $\mathfrak{p}^{n_1-n_2}\mathfrak{p}^{*(m_1-m_2)}$ primitive torsion points of $L$. Hence by distribution preperty we have right hand side of 4.1.2 is also equal to

$$\Theta(\frac{\Omega}{\pi^{n_1}\pi^{*m_1}}, \mathfrak{p}^{n_2-n_1}\mathfrak{p}^{*(m_2-m_1)}L, \boldsymbol{\alpha}) = \Theta(\Omega, \mathfrak{p}^{n_2}\mathfrak{p}^{*m_2}L, \boldsymbol{\alpha}) = e_{n_2,m_2}$$

Now it is easy to see $u_{n,m}$ are also norm compatible. □

1 (b): if p is a split bad reduction prime, define

$$e_{n,m} = \Theta(\Omega^{'}, \mathfrak{p}^n\mathfrak{p}^{*m}L^{'}, \boldsymbol{\alpha}) \in K(\mathfrak{g}\mathfrak{p}^n\mathfrak{p}^{*m})^*$$

where $L^{'}, \Omega^{'}$ and $g$ is defined the same way as above, but for $E^{'}$. Define

$$u^{'}_{n,m} = Norm_{K(\mathfrak{g}\mathfrak{p}^n\mathfrak{p}^{*m})/F_{n,m})}(e_{n,m})$$

and set $u_{n,m} = Norm_{F_{n,m}/K_{n,m}}(u^{'}_{n,m})$ for $n \geq m(\mathfrak{p})$ It follows similary as above that $u_{n,m}$ are norm compatible units.

2 (a): when $\mathfrak{p}$ is a good nonsplit prime, define $e_n = \Theta(\Omega, \mathfrak{p}^nL, \boldsymbol{\alpha}) \in K(\mathfrak{f}\mathfrak{p}^n)^*$ and set $u_n = Norm_{K(\mathfrak{f}\mathfrak{p}^n)/K_n}(e_n)$ which by the same proof above can be seen to be norm compatible units, hence indeed $u \in \mathfrak{A}_\infty$.

2(b): when $p$ is bad nonsplit prime, define $e_n = \Theta(\Omega^{'}, \mathfrak{p}^nL^{'}, \boldsymbol{\alpha}) \in K(\mathfrak{g}\mathfrak{p}^n)$ and set $u^{'}_n = Norm_{K(\mathfrak{g}\mathfrak{p}^n)/F_n}(e_n)$ and $u_n = Norm_{F_n/K_n}(u^{'}_n)$ for $n \geq m(\mathfrak{p})$

Now for each prime ideal $\mathfrak{p}$ of $\mathcal{O}_K$ which is a good reduction prime for $E$, fix once for all a basis for $T_\mathfrak{p}(E)$ as follows: Let $w_n$ be the unique $\mathfrak{f}$ division point of $\mathfrak{p}^nL$ such that

$$\phi^n(\varepsilon(\pi^{-n}w_n, L)) = \varepsilon(\Omega, L)$$

where $(\varepsilon, L) : \mathbb{C}/L \to E(\mathbb{C})$ is defined in the begining of section 2.3 (note this map is not Galois equivariant) and $\phi = frob_\mathfrak{p}$. It follows that $w_n \equiv w_{n-1} \pmod{\mathfrak{p}^{n-1}L}$. In particular, $w_n \equiv w_0 \equiv \Omega \pmod{L}$. Define $v_n = w_n - \Omega \bmod \mathfrak{p}^nL$ and let $\xi_n = \mathcal{P}(\pi^{-n}v_n)/\mathcal{P}^{'}(\pi^{-n}v_n)$. It is clear $\xi_n$ is a primitive $\pi^n$ torsion point of $\phi^{-n}\hat{E}_\mathfrak{p} = \hat{E}_\mathfrak{p}$ since $E$ is defined over $K$. Also by definition $[\pi]\xi_n = \xi_{n-1}$, hence $\{\xi_n\}_n \in T_\mathfrak{p}(E)$ is a basis. Note $\lambda(\xi_n) = \pi^{-n}v_n$, where $\lambda$ is the inverse of $\varepsilon$ which (as we have seen in the proof of lemma 2.2) can be viewed as the canonical logarithm: $\hat{E}_\mathfrak{p} \to \mathbb{G}_a$.

**Lemma 13.** *The Coleman power series of* $\{u_{n,m}\}_n$ *with respect to the chosen basis above is:*

$$\prod_{\sigma \in Gal(K(\mathfrak{f})/K)} P(\pi^{*(-m)}\lambda(t))^{\sigma} \tag{4.1.3}$$

*where* $P\ (z)$ *is the Taylor series expansion of* $\Theta(\pi^{*(-m)}\Omega - z, L, \boldsymbol{\alpha})$ *and* $\sigma$ *acts on coefficents. Recall* $t$ *is the variable on formal group* $\hat{E}_{\mathfrak{p}}$*. and* $\pi^* = \psi(\mathfrak{p}^*)$*. When* $p$ *is nonsplit, there is no* $m$ *and the statement about coleman power series holds without* $\pi^{*(-m)}$ *term.*

*Proof.* First notice for any $\mathfrak{m}$ torison point $\rho_{\mathfrak{m}}$ and integral ideal $\mathfrak{c}$ prime to $\mathfrak{m}$, we have

$$\Theta(z + \rho_{\mathfrak{m}}, L, \boldsymbol{\alpha})^{\sigma_{\mathfrak{c}}} = \Theta(z + \psi(\mathfrak{c})\rho_{\mathfrak{m}}, L, \boldsymbol{\alpha}) \tag{4.1.4}$$

This follows because $\Theta(z, L, \boldsymbol{\alpha})$ is a rational function of $\mathcal{P}(z)$ with coefficents in $K$, so by addition theorem $\Theta(z + \rho_{\mathfrak{m}}, L, \boldsymbol{\alpha})$ is a rational function of $\mathcal{P}(z)$ and $\mathcal{P}^{'}(z)$ with coefficent in $K(\mathfrak{fm})$. Since

$$\varepsilon(\rho_{\mathfrak{m}})^{(\mathfrak{c}, K(\mathfrak{fm})/K)} = \varepsilon(\psi(\mathfrak{c})\rho_{\mathfrak{m}})$$

we get 4.1.4 on applying $(\mathfrak{c}, K(\mathfrak{fm})/K)$ to the coefficents of $\Theta(z + \rho_{\mathfrak{m}}, L, \boldsymbol{\alpha})$.

By definition of Coleman power series, we want to show that

$$P^{\phi^{-n}}(\pi^{*(-m)}\lambda(t)) \mid_{t=\xi_n} = e_{n,m}$$

By the careful choice of $w_n$ it follows that

$$P^{\phi^{-n}}(z) = \Theta(w_n/\pi^n\pi^{*(-m)} - z, L, \boldsymbol{\alpha})$$

Hence

$$P^{\phi^{-n}}(\pi^{*(-m)}\lambda(t)) \mid_{t=\xi_n} = \Theta(\frac{w_n}{\pi^n\pi^{*m}} - \frac{v_n}{\pi^n\pi^{*m}}, L, \boldsymbol{\alpha}) = e_{n,m}$$

the last equality holds because of the choice of $v_n$ and that

$$\Theta(\frac{\Omega}{\pi^n\pi^{*m}}, L, \boldsymbol{\alpha}) = \Theta(\Omega, \mathfrak{p}^n\mathfrak{p}^m L, \boldsymbol{\alpha})$$

Finally it is clear that the Coleman power series of $\{u_{n,m}\}_n$ is "norm" from $K(\mathfrak{f})$ to K of the Coleman power series for $\{e_{n,m}\}_n$. $\square$

When $\mathfrak{p}$ is a bad reduction prime for $E$, as usual we replace $E$ by $E^{'}$and argue just the same to conclude that the Coleman power series for $\{u^{'}_{n,m}\}_n$ is :

$$\prod_{\sigma \in Gal(K(\mathfrak{g})/F)} P(\pi^{*(-m)}\lambda(t))^{\sigma} \tag{4.1.5}$$

which is to be interpretated the same way as above when $p$ is nonsplit prime.

Recall we want to evaluate $exp^*$ on $S_0(u)$. To apply Kato's result we have to calculate logarithmic derivatives of Coleman power series above which are essentially Eisenstein series. We brieftly recall as follows: ([3] 3.1) For integers $0 \leq -i < k$ (we are going to apply these results with $i = -j$)

(a):

$$E_{i,k}(z, L) \overset{\text{def}}{=} (k-1)! A(L)^i \sum_{\omega \in L} (z + \omega)^{-k} (\bar{z} + \bar{\omega})^{-i},\ k + i \geq 3 \tag{4.1.6}$$

where $A(L) = \pi^{-1} Area(\mathbb{C}/L)$, $\pi = 3.14\cdots$. Define

$$E_1(z, L) = \frac{1}{12}\frac{\partial}{\partial z} \log \Theta(z, L) - \frac{1}{2}\bar{z}A(L)^{-1}$$

Define $E_{i,k}(z, L, \boldsymbol{\alpha}) = N\boldsymbol{\alpha} E_{i,k}(z, L) - E_{i,k}(z, \boldsymbol{\alpha}^{-1}L)$ and denote $E_k(z, L) = E_{0,k}(z, L)$, $E_k(z, L, \boldsymbol{\alpha}) = E_{0,k}(z, L, \boldsymbol{\alpha})$.

(b): (see [3] p57-58)

$$(\frac{d}{dz})^k \log \Theta(z, L, \boldsymbol{\alpha}) = -12 E_k(z, L, \boldsymbol{\alpha}),\ \forall k \geq 1 \tag{4.1.7}$$

Here we can kill the constant 12 if we had used 12th root of $\Theta(z, L, \boldsymbol{\alpha})$, but since we are verifying Tamagawa number conjecture up to power of 2 and 3, we do not have to do that.

(c): Let $\mathfrak{m}$ be any nontrivial ideal and $\rho_{\mathfrak{m}}$ be a primitive $\mathfrak{m}$ division point of $L =$ period lattice of CM elliptic curve $E$ of conductor $\mathfrak{f}$. By 4.1.7 it is not surprising that we have following properties similar to that of theta functions: ([3] 3.3)

- $E_{i,k}(cz, cL) = c^{i-k} E_{i,k}(z, L)$
- $E_{i,k}(\rho_{\mathfrak{m}}) \in K(l.c.m(\mathfrak{f}, \mathfrak{m}))$
- $E_k(\rho_{\mathfrak{m}})$ is $\mathfrak{p}$ integral if $\mathfrak{m}$ is not a pure power of $\mathfrak{p}$, $\forall k \geq 2$ (see lemma 22)
- If $\mathfrak{c}$ is integral and $(\mathfrak{c}, \mathfrak{f}) = 1$, then

$$E_{i,k}(\rho_{\mathfrak{m}}, L)^{\sigma_{\mathfrak{c}}} = \psi(\mathfrak{c})^{i-k} E_{i,k}(\rho_{\mathfrak{m}}, \mathfrak{c}^{-1} L) \tag{4.1.8}$$

(d):([3] 3.2) There exists a unique polynomial $\Phi_{i,k}$ in $\mathbf{Z}[X_1, \cdots, X_{k-i}]$, of degree $1-i$, isobaric of weight $k-i$ ($X_k$ is given weight k) such that

$$E_{i,k} = 2^j \Phi_{i,k}(E_1, \cdots, E_{k-i})$$

Furthermore

$$\Phi_{i,k} = (-2X_1)^{-i} X_k + \text{ terms in which } X_1 \text{ appears to degree} < -i \tag{4.1.9}$$

(e): We have the following link between Eisenstein numbers and special values of partial Hecke L-series:([3] 3.5) let $\mathfrak{c}$ be an integral ideal prime to $\mathfrak{f}$ and $\mathfrak{f} \mid \mathfrak{m}$, $\forall \Omega \in \mathbb{C}$,

$$\mathbf{N}\mathfrak{m}^{-i} E_{i,k}(\Omega, \mathfrak{c}^{-1}\mathfrak{m}\Omega) = (k-1)!(\frac{\sqrt{d_K}}{2\pi})^i \Omega^{i-k} \psi(\mathfrak{c})^{k-i} L(\bar{\psi}^{k-i}, (\frac{K(\mathfrak{m})/K}{\mathfrak{c}}), k) \tag{4.1.10}$$

Here the partial $L$ value means summing over all integral ideals of $K$ prime to $\mathfrak{m}$ and whose artin symbol in $Gal(K(\mathfrak{m})/K)$ is equal to that of $\mathfrak{c}$. Combined with 4.1.9 this gives a link between $L(\bar{\psi}^{k-i}, k)$ with $L(\bar{\psi}^k, k)$. We delay this study till we need it later.

## 4.2 Calculation in the case $j = 0$

**Lemma 14.** *Let $L/K$ be a finite extension of $p$-adic local fields. Let $V$ be a de Rham representation of $Gal(\bar{K}/K)$. Then (1) the following diagram commutes:*

$$\begin{array}{ccc} H^1(L,V) & \xrightarrow{exp_L^*} & DR^0(K,V)\otimes_K L \\ \uparrow{\scriptstyle res} & & \uparrow{\scriptstyle inclu} \\ H^1(K,V) & \xrightarrow{exp_K^*} & DR^0(K,V) \end{array}$$

*(2) $exp^*$ commutes with Gal(L/K) action, where $\forall \phi \in H^1(L,V), \forall g \in Gal(\bar{K}/L)$*

$$\sigma\phi(g) \overset{def}{=} \sigma(\phi(\bar{\sigma}^{-1}g\bar{\sigma})), \forall \sigma \in Gal(L/K), \bar{\sigma} \text{ is a lift of } \sigma \text{ to } Gal(\bar{K}/K)$$

*Proof.* The commutativity of the diagram follows from the functorial property of either definition of $exp^*$. Let $exp^*(\phi) = a$, i.e by second definition, $\exists m \in V \otimes B_{\mathrm{dR}}$, such that

$$\phi(g) - \log \chi_{\mathrm{cyclo}}(g) \cup a = gm - m, \ \forall g \in Gal(\bar{K}/L)$$

in particular we have

$$\phi(\bar{\sigma}^{-1}g\bar{\sigma}) - \log \chi(\bar{\sigma}^{-1}g\bar{\sigma}) \cup a = (\bar{\sigma}^{-1}g\bar{\sigma})m - m$$

Let $\sigma$ acts on both side of the above equation, noting that $\chi(\bar{\sigma}^{-1}g\bar{\sigma}) = \chi(g)$, we get

$$(\sigma\phi)(g) - \log \chi(g) \cup (\sigma a) = g(\sigma m) - \sigma m$$

hence by definition $exp^*(\sigma\phi) = \sigma exp^*(\phi)$. □

*Remark 9.* Since $V$ is a vector space, the above restriction map is injective. We can view the bottom row as the $Gal(L/K)$ fixed part of the upper row.

**Lemma 15.** *Let $K/\mathbf{Q}_p$ be finite extension. Let $V_i$ be de Rham representation of $Gal(\overline{\mathbf{Q}_p}/K)$ for $i = 1, 2$. Suppose given a $Gal(\overline{\mathbf{Q}_p}/K)$ equivariant map $f : V_1 \to V_2$, then the following diagram commutes:*

$$\begin{array}{ccc} H^1(K,V_1) & \xrightarrow{f} & H^1(K,V_2) \\ \downarrow{\scriptstyle exp^*} & & \downarrow{\scriptstyle exp^*} \\ DR^0(V_1) & \xrightarrow{f} & DR^0(V_2) \end{array}$$

*Proof.* $\forall \phi \in H^1(K, V_1), \forall g \in Gal(\overline{\mathbf{Q}_p}/K)$, define $(f\phi)(g) = f(\phi(g))$. Check that $f\phi$ such defined is a cocycle:

$$\begin{aligned}
(f\phi)(g_1g_2) &= f(\phi(g_1g_2)) && \text{by definition} \\
&= f(\phi(g_1) + g_1\phi(g_2)) && \text{since } \phi \text{ is a cocycle} \\
&= f(\phi(g_1)) + g_1 f(\phi(g_2)) && \text{since f is } Gal(\overline{\mathbf{Q}_p}/K) \text{ equivariant}
\end{aligned}$$

Let $exp^*(\phi) = a$. Applying $f$ to the defining equation for $a$, we get

$$f(\phi(g)) - \log\chi(g) \cup (fa) = f(gm) - fm = g(fm) - fm$$

hence $exp^*(f\phi) = fexp^*(\phi)$. □

### 4.2.1 $j = 0$ and $p$ is a good reduction prime

In this case $K_\infty/K$ is totally ramified at $\mathfrak{p}$ and $v$ is the unique place of $K_\infty$ above $\mathfrak{p}$. Observe by definition (see definition 4 in section 3) of $S_n(u_v) \in H^1(K_{n,v}, T_{\mathfrak{p}}(E)^{\otimes(-k)}(1))$ that

$$S_0(u_v) \mapsto \sum_{\tau \in Gal(K_{n,v}/K_{\mathfrak{p}})} \tau S_n(u_v)$$

under the restriction map

$$H^1(K_{\mathfrak{p}}, T_{\mathfrak{p}}(E)^{\otimes(-k)}(1)) \to H^1(K_{n,v}, T_{\mathfrak{p}}(E)^{\otimes(-k)}(1))$$

using Kato's results (theorem 1), we have

$$\begin{aligned}
\frac{exp^*(S_n(u_v)}{\varpi} &= \frac{1}{(k-1)!}\pi^{-nk}\{(\frac{d}{\omega})^k \log(\phi^{-n} g_{u,\xi})\}(\xi_n) \\
&\overset{*}{=} \frac{1}{(k-1)!}\pi^{-nk}\{(\frac{d}{dz})^k \log(\prod_{\sigma\in Gal(K(\mathfrak{f})/K)} (\sigma P)^{\phi^{-n}}(z))\}(v_n) \\
&= \frac{1}{(k-1)!}\pi^{-nk}\{(\frac{d}{dz})^k \sum_{\sigma\in Gal(K(\mathfrak{f})/K)} \log\Theta^\sigma(-z+\frac{w_n}{\pi^n}, L, \boldsymbol{\alpha})\}(\frac{v_n}{\pi^n}) \\
&= -\frac{12}{(k-1)!}\pi^{-nk}\{\sum_{\sigma\in Gal(K(\mathfrak{f})/K)} E_k^\sigma(-z+\frac{w_n}{\pi^n}, L, \boldsymbol{\alpha})\}((\frac{v_n}{\pi^n}) \\
&= -\frac{12}{(k-1)!}\pi^{-nk}\sum_{\sigma\in Gal(K(\mathfrak{f})/K)} E_k^\sigma(\frac{\Omega}{\pi^n}, L, \boldsymbol{\alpha})
\end{aligned} \tag{4.2.1}$$

where equality $(*)$ follows from lemma 12 and by definition of $\xi_n$. Here $\omega = \lambda^{'}(T)dT$ and $\varpi = \omega^{\otimes(k+j)}\epsilon^{\otimes(-j)}$. Now we relate the last equality to special values of Hecke $L$-series as follows:

$$\begin{aligned}
&\sum_{\sigma\in Gal(K(\mathfrak{f})/K)} E_k^\sigma(\frac{\Omega}{\pi^n}, L, \boldsymbol{\alpha}) \\
&\overset{(1)}{=} \sum_{\mathfrak{b}\in B} \mathbf{N}\boldsymbol{\alpha} E_k(\psi(\mathfrak{b})\frac{\Omega}{\pi^n}, L) - E_k(\psi(\mathfrak{b})\frac{\Omega}{\pi^n}, \boldsymbol{\alpha}^{-1}L) \\
&\overset{(2)}{=} \sum_{\mathfrak{b}\in B} \mathbf{N}\boldsymbol{\alpha} E_k(\psi(\mathfrak{b})\frac{\Omega}{\pi^n}, \frac{\Omega}{\pi^n}\mathfrak{f}\mathfrak{p}^n) - E_k(\psi(\mathfrak{b})\frac{\Omega}{\pi^n}, \psi(\boldsymbol{\alpha}^{-1})\frac{\Omega}{\pi^n}\mathfrak{f}\mathfrak{p}^n) \qquad (4.2.2) \\
&= \sum_{\mathfrak{b}\in B} \left(\frac{\Omega}{\pi^n}\right)^{-k} \{\mathbf{N}\boldsymbol{\alpha} E_k(\psi(\mathfrak{b}), \mathfrak{f}\mathfrak{p}^n) - \psi^k(\boldsymbol{\alpha}) E_k(\psi(\boldsymbol{\alpha}\mathfrak{b}), \mathfrak{f}\mathfrak{p}^n)\} \\
&= \left(\frac{\Omega}{\pi^n}\right)^{-k} \{\mathbf{N}\boldsymbol{\alpha} L_{K_n}(\psi^{-k}, 1, 0) - \psi^k(\boldsymbol{\alpha}) L_{K_n}(\psi^{-k}, \sigma_{\boldsymbol{\alpha}}, 0)\}
\end{aligned}$$

where $B$ is a set of integral ideals of $K$ whose artin symbol in $Gal(K(\mathfrak{f}\mathfrak{p}^n)/K)$ describe precisely $Gal(K(\mathfrak{f}\mathfrak{p}^n)/K_n)$ and in the last equation,

$$L_{K_n}(\psi_k, \sigma_{\boldsymbol{\alpha}}, s) = \sum_{\substack{\mathfrak{c}\in\mathcal{O}_K, (\mathfrak{c},\mathfrak{f}\mathfrak{p})=1 \\ (\mathfrak{c}, K_n/K)=\sigma_{\boldsymbol{\alpha}}}} \psi^{-k}(\mathfrak{c})\mathbf{N}\mathfrak{c}^{-s}$$

(1) holds by definition and (2) follows since $\psi(\boldsymbol{\alpha})$ is a generator of $\boldsymbol{\alpha}$. Compare with

$$E_k(\psi(\mathfrak{b}), \mathfrak{f}\mathfrak{p}^n) = \sum_{\omega \in \mathfrak{f}\mathfrak{p}^n} \frac{(k-1)!}{(\psi(\mathfrak{b}) + \omega)^k}, \, \forall k \geq 2, \forall \mathfrak{b} \in B$$

Notice $\{(\psi(\mathfrak{b}) + \omega) : \mathfrak{b} \in B, \omega \in \mathfrak{f}\mathfrak{p}^n\}$ is precisely the set of integral ideals of $K$ prime to $\mathfrak{f}\mathfrak{p}^n$ whose artin symbol in $Gal(K_n/K)$ is $\sigma_{\boldsymbol{\alpha}}$ by choice of $B$. Also since $\psi$ has conductor $\mathfrak{f}$ , we have

$$\psi((\psi(\boldsymbol{\alpha}\mathfrak{b}) + \omega)) = \psi(\boldsymbol{\alpha}\mathfrak{b}) + \omega, \, \forall \omega \in \mathfrak{f}\mathfrak{p}^n$$

It follows from the definition of $E_k(z, L)$ that

$$E_k(\psi(\boldsymbol{\alpha}\mathfrak{b}), \mathfrak{f}\mathfrak{p}^n) = (k-1)! L_{K(\mathfrak{f}\mathfrak{p}^n)}(\psi^{-k}, \sigma_{\boldsymbol{\alpha}\mathfrak{b}}, 0), \, \forall (\boldsymbol{\alpha}, \mathfrak{f}\mathfrak{p}) = 1$$

and also

$$\sum_{\mathfrak{b} \in B} E_k(\psi(\boldsymbol{\alpha}\mathfrak{b}), \mathfrak{f}\mathfrak{p}^n) = (k-1)! L_{K_n}(\psi^{-k}, \sigma_{\boldsymbol{\alpha}}, 0) \tag{4.2.3}$$

Combining these equations, we get:

$$\begin{aligned} \frac{exp^*(S_0(u_v))}{\varpi} &= \sum_{\tau \in Gal(K_{n,v}/K_{\mathfrak{p}})} \tau exp^*(S_n(u_v)) \\ &\stackrel{*}{=} -12\omega^{-k} \sum_{\substack{\mathfrak{b} \in B \\ \mathfrak{c} \in C}} \mathbf{N}\boldsymbol{\alpha} E_k(\psi(\mathfrak{b}\mathfrak{c}), \mathfrak{f}\mathfrak{p}^n) - \psi^k(\boldsymbol{\alpha}) E_k(\psi(\boldsymbol{\alpha}\mathfrak{b}\mathfrak{c}), \mathfrak{f}\mathfrak{p}^n) \qquad (4.2.4) \\ &= -12 \left(\frac{\Omega_\infty}{f}\right)^{-k} (\mathbf{N}\boldsymbol{\alpha} - \psi^k(\boldsymbol{\alpha})) L_{\mathfrak{f}\mathfrak{p}}(\psi^{-k}, 0) \end{aligned}$$

where $C$ is a set of integral ideals of $K$ whose artin symbol in $Gal(K(\mathfrak{f}\mathfrak{p}^n)/K)$ describe $Gal(K_{n,v}/K_{\mathfrak{p}}) = Gal(K_n/K)$ and $L_{\mathfrak{f}\mathfrak{p}}(\psi^{-k}, s)$ means the incomplete complex $L$-series with Euler factor at $\mathfrak{p}$ and places dividing $\mathfrak{f}$ deleted.

*Remark 10.* 1. equality $(*)$ follows either from lemma 14 (2) or one can take note of the fact that $\tau S_n(u_v) = S_n((u^\tau)_v)$ by definition and that the Coleman power series for $u^\tau$ is simply $\tau$ acting on the coleman power series for $u$. Then we can repeat the above argument exactly the same way for $u^\tau$. Later in the case when $j > 0$, we do not have analogue of lemma 14 as we do not quite have the exact definition of finite version of $exp^*$. Then we will apply this alternative method.

2. note in the above calculation we take any $n \geq 1$ and get the same results (as we should). One can just take $n = 1$ in this case.

### *4.2.2 $j = 0$ and $p$ is a bad reduction prime*

In this case, $K_\infty/K$ is still totally ramified at $\mathfrak{p}$, $\upsilon$ being the unique place lying above $\mathfrak{p}$ in $K_\infty$. The places of $F_\infty$ lying above $\mathfrak{p}$ are $\{\upsilon^\sigma; \sigma \in D\}$ where $D$ is a set of coset representatives of $Gal(F/K)/Gal(F_\upsilon/K_\mathfrak{p})$. We want to relate $exp^*(S_0(u_\upsilon))$ with $exp^*(S_n'(u_\upsilon'))$ which we can calculate by the same method from last section applied to $(F, E')$. The following argument works for any $n \geq m(\mathfrak{p})$. In particular one can just take $n = m(\mathfrak{p})$. By lemma 14 and 15, we have the following commutative diagram:

$$\begin{array}{ccccc} H^1(K_\mathfrak{p}, T_\mathfrak{p}(E)^{\otimes(-k)}(1)) & \xrightarrow{\text{res}} & \oplus_{\mathfrak{P}|\mathfrak{p}} H^1(F_{n,\mathfrak{P}}, T_\mathfrak{p}(E)^{\otimes(-k)}(1)) & \xrightarrow{f} & \oplus_{\mathfrak{P}|\mathfrak{p}} H^1(F_{n,\mathfrak{P}}, T_\mathfrak{p}(E')^{\otimes(-k)}(1)) \\ \downarrow{\scriptstyle exp^*} & & \downarrow{\scriptstyle exp^*} & & \downarrow{\scriptstyle exp^*} \\ D^0(K_\mathfrak{p}, V_\mathfrak{p}(E)^{\otimes(-k)}(1)) & \rightarrow & D^0(K_\mathfrak{p}, V_\mathfrak{p}(E)^{\otimes(-k)}(1)) \otimes F_{n,\mathfrak{P}} & \xrightarrow{f} & \oplus_{\mathfrak{P}|\mathfrak{p}} colie(G')^{\otimes k} \otimes F_{n,\mathfrak{P}} \end{array}$$

Here the maps (induced by) $f$ are isomorphisms. It follows

$$\begin{aligned} exp^*(S_0(u_\upsilon)) &= f^{-1}\{exp^* \circ f \circ (\oplus_{\mathfrak{P}|\mathfrak{p}}\text{res})(S_0(u_\upsilon)\} \\ &= f^{-1}\{exp^* \circ f(\sum_{\tau \in Gal(F_{n,\upsilon}/K_\mathfrak{p}) \times D} \tau S_n(u_\upsilon'))\} \\ &= f^{-1}\{exp^*(\sum_{\tau \in Gal(F_n/K)} \epsilon^k(\tau)\tau S_n'(u_\upsilon'))\} \end{aligned} \tag{4.2.5}$$

The last equality holds because:

**Lemma 16.** *Let $L/K_\mathfrak{p}$ be finite extension. Let $T_i, i = 1, 2$ be rank 1 $\mathcal{O}_\mathfrak{p}$ module with continuous $Gal(\overline{K_\mathfrak{p}}/K_\mathfrak{p})$ action. Suppose given $f : T_1 \to T_2$ a $Gal(\overline{K_\mathfrak{p}}/L)$ isomorphism under which $T_2 = T_1(\eta)$ where $\eta$ is a finite character of $Gal(L/K_\mathfrak{p})$. Then*

$$f(\sigma\phi) = \eta^{-1}(\sigma)\sigma(f\phi), \forall \sigma \in Gal(L/K_\mathfrak{p}), \phi \in H^1(L, T_1)$$

*Same conclusion holds for finite version of this, i.e replace $T_i$ by $T_i/\mathfrak{p}^n$.*

*Proof.* Under map $f$ we have $\sigma(fa) = \eta(\sigma)f(\sigma a)$, $\forall a \in T_1, \sigma \in Gal(L/K_{\mathfrak{p}})$, i.e $\sigma f = \eta(\sigma)f\sigma$. Hence $\forall g \in Gal(\overline{K_{\mathfrak{p}}}/L)$,

$$[\sigma(f\phi)](g) = (\sigma(f\phi))(\bar{\sigma}^{-1}g\bar{\sigma}) = \eta(\sigma)f(\sigma\phi(\bar{\sigma}^{-1}g\bar{\sigma})) = \eta(\sigma)(f(\sigma\phi))(g)$$

so $\sigma(f\phi) = \eta(\sigma)f(\sigma\phi)$. □

We apply this lemma here with $T_1 = T_{\mathfrak{p}}(E)^{\otimes(-k)}(1)$ and $T_2 = T_{\mathfrak{p}}(E')^{\otimes(-k)}(1)$. Recall $\psi' = \psi\epsilon$, so in applying the above lemma, we take $\eta = \epsilon^{-k}$. Now we calculate $exp^*(S_n'(u_v'))$ using Kato's theorem:

$$\frac{exp^*(\tau S_n'(u_v'))}{\varpi'} = -12(\Omega')^{-k}\{\mathbf{N}\boldsymbol{\alpha} L_{F_n}(\psi'^{(-k)}, \tau, 0) - \psi'^{(-k)}(\boldsymbol{\alpha})L_{F_n}(\psi'^{(-k)}, \sigma_{\boldsymbol{\alpha}}\tau, 0)\} \tag{4.2.6}$$

where $\tau \in Gal(F_n/K)$. Recall when $j = 0, \varpi = \omega^{\otimes k}$. Also $f^{-1}(\Omega_\infty') = r\Omega_\infty$ and $f^*\omega' = r\omega$. Combining 4.2.5 and 4.2.6 we get

$$\begin{aligned}
\frac{exp^*(\tau S_0(u_v))}{\varpi} &= r^k\{exp^*(\sum_{\tau \in Gal(F_n/K)} \epsilon^k(\tau)\tau S_n'(u_v'))\}/\varpi' \\
&= -12r^k \sum_{\tau \in Gal(F_n/K)} \epsilon^k(\tau)\{(\Omega')^{-k}(\mathbf{N}\boldsymbol{\alpha} L_{F_n}(\psi'^{(-k)}, \tau, 0) \\
&\quad - \psi'^{(-k)}(\boldsymbol{\alpha})L_{F_n}(\psi'^{(-k)}, \sigma_{\boldsymbol{\alpha}}\tau, 0))\} \\
&= -12\left(\frac{\Omega_\infty}{g}\right)^{-k}\{\mathbf{N}\boldsymbol{\alpha} L_{\mathfrak{g}\mathfrak{p}}((\psi'\epsilon^{-1})^{-k}, 0) \qquad (4.2.7) \\
&\quad - \sum_{\tau \in Gal(F_n/K)} \psi'^k(\boldsymbol{\alpha})\epsilon^{-k}(\boldsymbol{\alpha})\epsilon^k(\tau\sigma_{\boldsymbol{\alpha}})L_{F_n}(\psi'^{(-k)}, \tau\sigma_{\boldsymbol{\alpha}}, 0)\} \\
&= -12\left(\frac{\Omega_\infty}{g}\right)^{-k}\{\mathbf{N}\boldsymbol{\alpha} L_{\mathfrak{g}\mathfrak{p}}(\psi^{-k}, 0) - \psi^k(\boldsymbol{\alpha})L_{\mathfrak{g}\mathfrak{p}}(\psi^{-k}, 0)\} \\
&= -12\left(\frac{\Omega_\infty}{g}\right)^{-k}(\mathbf{N}\boldsymbol{\alpha} - \psi^k(\boldsymbol{\alpha}))L_{\mathfrak{g}\mathfrak{p}}(\psi^{-k}, 0)
\end{aligned}$$

# CHAPTER 5

# DUAL EXPONENTIAL MAP WHEN $J > 0$

Here we will deal with $p$ splits and $p$ is nonsplit case seperately. Fix integers $k, j$ such that $0 < j < k$ and $k - j > 1$. Consider a general setting: $E$ is an CM elliptic curve over $K$, $\mathfrak{p}$ a prime of $K$ such that $E$ has good reduction at $\mathfrak{p}$. Define

$$V_p = \begin{cases} V_{\mathfrak{p}}(E)^{\otimes(k+j)}(-j) \oplus V_{\mathfrak{p}^*}(E)^{\otimes(k+j)}(-j) & \text{if } p \text{ splits} \\ V_p(E)^{\otimes(k+j)}(-j) & \text{if } p \text{ nonsplit} \end{cases} \tag{5.0.1}$$

Equivalently,

$$V_p = \begin{cases} \mathbf{Q}_p(\psi_{\mathfrak{p}}^{\otimes(k+j)}\chi^{-j}) \oplus \mathbf{Q}_p(\psi_{\mathfrak{p}^*}^{\otimes(k+j)}\chi^{-j}) & \text{if } p \text{ splits} \\ K_{\mathfrak{p}}(\psi_{\mathfrak{p}}^{\otimes(k+j)}) \otimes_{\mathbf{Q}_p} \mathbf{Q}_p(-j) & \text{if } p \text{ nonsplit} \end{cases}$$

Let $W_p = V_p(j)$. We have

$$DR(W_p) = Fil^{-(k+j)}DR(W_p) \supset \cdots Fil^0 DR(W_p) \supset Fil^1 DR(W_p) = 0$$

In particular when $p$ splits, since $\psi_{\mathfrak{p}^*}$ is unramified at $\mathfrak{p}$, we get

$$Fil^0 DR(W_p) = DR(\mathbf{Q}_p(\psi_{\mathfrak{p}^*}^{\otimes(k+j)})), \quad gr^0 DR(W_p) = DR(\mathbf{Q}_p(\psi_{\mathfrak{p}}^{\otimes(k+j)}))$$

Since filtration for $DR(V_p)$ is just a shift of the above filtration to the right by $j$ we conclude that $gr^{-k}DR(V_p)$ and $gr^j DR(V_p)$ are 1 dimensional $\mathbf{Q}_p$ vector space and $gr^n DR(V_p) = 0, \forall n \neq -k, j$. When $p$ splits, we still have $DR(V_p)/DR^0(V_p) = DR(\mathbf{Q}_p(\psi_{\mathfrak{p}}^{\otimes(k+j)}\chi^{-j}))$.

With calculation of $exp^*$ in mind, let us now consider $DR(L, V_{\mathfrak{p}}^*(1))$ where $L/K_{\mathfrak{p}}$ is a finite extension and for $\mathfrak{p} \mid p$

$$V_{\mathfrak{p}} \overset{\text{def}}{=} V_{\mathfrak{p}}(E)^{\otimes(k+j)}(-j)$$

Note $DR^0(L, V_{\mathfrak{p}}^*(1))$ is 1 dimensional over $L$. Let $T_{\mathfrak{p}}$ be defined similarly using $T_{\mathfrak{p}}(E)$.

**Lemma 17.**

$$DR^0(L, V_{\mathfrak{p}}^*(1)) \cong colie(G)^{\otimes(k+j)} \otimes colie(\mathbb{G}_m)^{\otimes(-j)} \otimes L$$

*where $G$ is the $p$ divisible group associated to $\hat{E}_{\mathfrak{p}}$ and $\mathbb{G}_m$ is the formal multiplicative group.*

*Proof.* the key point is the following commutative diagram:

$$\begin{array}{ccc} colie(G)^{\otimes(k+j)} \otimes B_{\mathrm{dR}}^+ & \xrightarrow[\psi_1]{\cong} & V_{\mathfrak{p}}(E)^{\otimes(-k-j)} \otimes B_{\mathrm{dR}}^{k+j} \\ \Big\downarrow{\scriptstyle (ds/1+s)^{\otimes(-j)}} & & \Big\downarrow{\scriptstyle \theta} \\ colie(G)^{\otimes(k+j)} \otimes colie(\mathbb{G}_m)^{\otimes(-j)} \otimes B_{\mathrm{dR}}^+ & \xrightarrow[\psi_2]{} & V_{\mathfrak{p}}(E)^{\otimes(-k-j)}(j) \otimes B_{\mathrm{dR}}^k \end{array}$$

where map $\theta$: $\forall v \otimes x \in V_{\mathfrak{p}}(E)^{\otimes(-k-j)} \otimes B_{\mathrm{dR}}^{k+j}, \theta(v \otimes x) = (v \otimes \epsilon^j) \otimes (xt^{-j})$. Here $\epsilon$ is our fixed basis of $\mathbf{Z}_p(1)$ and define $t = \log[\epsilon]$ as usual in the theory of $B_{\mathrm{crys}}$ which gives us a Galois equivariant embedding:

$$\mathbf{Z}_p(1) \hookrightarrow J_\infty : \epsilon \mapsto t$$

The above diagram commutes because the horizontal maps come from pairing

$$colie(G) \otimes T_{\mathfrak{p}}(E) \to J_\infty, colie(\mathbb{G}_m) \otimes \mathbf{Z}_p(1) \to J_\infty$$

Using definition of this pairing, we see $(ds/1+s, \epsilon) \mapsto t$ as the logarithm of $\mathbb{G}_m$ associated to $ds/1+s$ is just the usual function log. This says the diagram is commutative. Now notice the vertical maps are Galois equivariant isomorphisms. Kato proved $\psi_1$ is a Galois equivariant isomorphism, which follows from the pairing above. Hence we see $\psi_2$ is also an isomorphism. Take $H^0(L, -)$ of it, we get:

$$\begin{aligned} colie(G)^{\otimes(k+j)} \otimes colie(\mathbb{G}_m)^{\otimes(-j)} \otimes L &\cong H^0(L, V_{\mathfrak{p}}(E)^{\otimes(-k-j)}(j) \otimes B_{\mathrm{dR}}^k) \\ &\overset{(*)}{=} H^0(L, V_{\mathfrak{p}}(E)^{\otimes(-k-j)}(j) \otimes B_{\mathrm{dR}}^1) \\ &= H^0(L, V_{\mathfrak{p}}(E)^{\otimes(-k-j)}(1+j) \otimes B_{\mathrm{dR}}^+) \\ &= DR^0(L, V_{\mathfrak{p}}) \end{aligned} \tag{5.0.2}$$

where $(*)$ follows from $p$-adic Hodge-Tate decomposition of $V_\mathfrak{p}$ as we remarked above that $V_\mathfrak{p}(E)^{\otimes(-k-j)}(j)\otimes B^i_{\mathrm{dR}}/B^{i+1}_{\mathrm{dR}} \cong \mathbb{C}_p(-k+i)$ and the fact that $H^0(L,\mathbb{C}_p(-k+i)) = 0$, since when $1 \leq i \leq k-1$, $-k+i \neq 0$ and $H^0(L,\mathbb{C}_p(m)) = 0, \forall m \neq 0$. $\square$

We also need the integral version of above lemma, namely there exists integer $N$ such that under $\psi_2$ in lemma 17,

$$colie(G)^{\otimes(k+j)} \otimes colie(\mathbb{G}_m)^{\otimes(-j)} \otimes B_\infty \cong p^N T_\mathfrak{p}(E)^{\otimes(-k-j)}(j) \otimes J^k_\infty \tag{5.0.3}$$

We also have the following analogue of 3.3.5: $\forall L$ finite extension of $K_\mathfrak{p}$,

$$\begin{aligned} colie(G)^{\otimes(k+j)} \otimes colie(\mathbb{G}_m)^{\otimes(-j)} \otimes \mathcal{O}_L &\hookrightarrow H^0(L, T_\mathfrak{p}(E)^{\otimes(-k-j)}(j) \otimes J^k_\infty/J^{[k+1]}_\infty) \\ &\to H^0(L, T_\mathfrak{p}(E)^{\otimes(-k-j)}(j) \otimes J_\infty/J^{[k+1]}_\infty) \end{aligned} \tag{5.0.4}$$

with cokernel killed by $p^N$, where $N$ is independent of $L$.

There is a unique homomorphism of formal groups over $\mathcal{O}_{\mathbb{C}_p}$

$$\eta : \hat{E}_\mathfrak{p} \to \mathbb{G}_m, \eta(\xi) = \epsilon$$

and $\exists \Omega_p \in \mathcal{O}_{\mathbb{C}_p}$ making the following diagram commutes:

$$\begin{array}{ccc} \hat{E}_\mathfrak{p} & \xrightarrow{\eta} & \mathbb{G}_m \\ & {\scriptstyle \Omega_p\lambda}\searrow & \uparrow{\scriptstyle exp} \\ & & \mathbb{G}_a \end{array}$$

This is well known in the case $p$ splits, since then $\hat{E}_\mathfrak{p}$ and $\mathbb{G}_m$ are two 1 dimensional Lubin-Tate formal group over $\mathbf{Q}_p$ which correspond to uniformizer $\pi = \psi(\mathfrak{p})$ and $p$ respectively. Moreover they become isomorphic over $\mathbf{Q}_p^{\mathrm{unr}}$. In general a result of Katz [14](section 6) says that $Hom_{\mathcal{O}_{\mathbb{C}_p}}(\hat{E}_\mathfrak{p}, \mathbb{G}_m)$ is a free $\mathcal{O}_\mathfrak{p}$ module of rank 1. Since $\rho(T) \stackrel{\text{def}}{=} \log(1+\eta(T))$ satisfies $\rho(T_1[+]T_2) = \rho(T_1) + \rho(T_2)$, it follows from the uniqueness of normalized logarithm of $\hat{E}_\mathfrak{p}$ that $\exists \Omega_p \in \mathcal{O}_{\mathbb{C}_p}$, such that $\rho = \Omega_p\lambda$ ($\Omega_p$ is just the coefficent of T in $\rho(T)$). From this we have

$$\eta^*(\frac{ds}{1+s}) = \Omega_p\lambda^{'}(T)dT = \Omega_p\omega, \text{ by lemma2}$$

Since $[\pi] = [p\alpha^{-1}]$ for some $\alpha \in \mathcal{O}^*_{K_\mathfrak{p}}$(even in $\mathcal{O}^*_K$ when $p$ is nonsplit), again following the definition of the above paring, we get:

$$\omega \otimes \xi \mapsto \Omega_p^{-1}\alpha t$$

## 5.1 $p$ splits

### 5.1.1 $p$ is a good reduction prime

In this case we have seen in lemma 3 there is an integer $M$ such that $\forall m \geq M$, the set of places in $K_{0,m}$ above $\mathfrak{p}$ is the same as the set of places of $K_\infty$ lying above $\mathfrak{p}$ and they are in 1-1 correspondence with $B$ which is a set of coset representatives of $Gal(K_{0,M}/K)/Gal(K_{o,M,v}/K_\mathfrak{p})$. In this section, we denote $V'_\mathfrak{p} = V_\mathfrak{p}(E)^{\otimes k}$ and $T'_\mathfrak{p} = T_\mathfrak{p}(E)^{\otimes k}$ (we work exactly with $V'_\mathfrak{p}$ when $j = 0$).

**Theorem 4.** *(1) $\exists$ integer $c, e$ both independent of $m$ and $A_m \in \mathcal{O}_{K_\mathfrak{p}}/\mathfrak{p}^m$, such that $\forall m > M$ we have equality in $H^1(K_\mathfrak{p}, (T^*_\mathfrak{p} \otimes J_\infty/J_\infty^{[k+1]})/\mathfrak{p}^m)$ :*

$$p^e c \sum_{\mathfrak{P}|\mathfrak{p}} \overline{S_{0,m}}(u_\mathfrak{P}) = p^e A_m \varpi \cup \log \chi_G \tag{5.1.1}$$

*(2) $\exists l \in \{1, \cdots p-1\}, \exists A \in K_\mathfrak{p}$, such that $cA \in \mathcal{O}_{K_\mathfrak{p}}$, and for all $m$ big enough we have*

$$cA \equiv A_{ml} \pmod{\mathfrak{p}^m}$$

*(3)*

$$\sum_{\mathfrak{P}|\mathfrak{p}} \frac{exp^*(S_0(u_\mathfrak{P}))}{\varpi} = -12\alpha^{-j}\left(\frac{\sqrt{d_K}}{2\pi}\right)^j \Omega^{-j-k}(\mathbf{N}\boldsymbol{\alpha} - \psi^k\bar{\psi}^{-j}(\boldsymbol{\alpha}))L_{\mathfrak{f}\mathfrak{p}\mathfrak{p}^*}(\bar{\psi}^{k+j}, k) \tag{5.1.2}$$

*Remark 11.* • It follows from (1) of the theorem and the definition of $exp^*$ combined with the fact that $\log \chi_G$ and $\log \chi_{\text{cyclo}}$ have the same image in $H^1(K_\mathfrak{p}, \mathbb{C}_p)$ that:

$$c \sum_{\mathfrak{P}|\mathfrak{p}} \frac{exp^*(S_0(u_\mathfrak{P}))}{\varpi} \equiv A_m \pmod{\mathfrak{p}^{m-e}}, \forall m > M + e$$

By (2) and the fact that c is independent of $m$, we have

$$\sum_{\mathfrak{P}|\mathfrak{p}} exp^*(S_0(u_{\mathfrak{P}})) = A\varpi \tag{5.1.3}$$

In the process of proving this theorem, we will give an explicit formula for $A_m$ and $A$ (i.e just the righ hand side of 5.1.2). In particular they are closely related to $L(\mathcal{M}_{k,j}, 0)$ as we want.

- This theorem which is an analogue of Kato's theorem when $j = 0$ (see theorem 1) seems to apply only to special element $S_0(u)$ where $u$ is our chosen generator of elliptic units. In the proof we need to be able to compute explicitly the logarithmic derivative of coleman power series for $u$ (which for our $u$ is essentially Eisenstein series) and to pass from $j = 0$ case to $j > 0$ case. For our choice of $u$ we rely on some special property of Eisenstein numbers. However it seems impossible to describe explicitly the image of $exp^*$ on arbitrary element of $H^1(K_{\mathfrak{p}}, V_{\mathfrak{p}}^*(1))$.

*Proof.* The key point is the following 2 commutative diagram: $\forall \mathfrak{P} \mid \mathfrak{p}$

$$\begin{array}{ccc}
H^1(K_{\mathfrak{p}}, T_{\mathfrak{p}}^*(1)/\mathfrak{p}^m) & \xrightarrow{\rho_1} & H^1(K_{1,m,\mathfrak{P}}, T_{\mathfrak{p}}^*(1)/\mathfrak{p}^m) \\
\downarrow \theta_1 & & \downarrow \theta_2 \\
H^1(K_{\mathfrak{p}}, (T_{\mathfrak{p}}^* \otimes \frac{J_\infty}{J_\infty^{[k+1]}})/\mathfrak{p}^m) & \xrightarrow{\rho_3} & H^1(K_{1,m,\mathfrak{P}}, (T_{\mathfrak{p}}^* \otimes \frac{J_\infty}{J_\infty^{[k+1]}})/\mathfrak{p}^m) \\
\cup \log \chi_G \uparrow & & \cup \log \chi_G \uparrow \\
H^0(K_{\mathfrak{p}}, (T_{\mathfrak{p}}^* \otimes \frac{J_\infty}{J_\infty^{[k+1]}})/\mathfrak{p}^m) & \xrightarrow{\rho_5} & H^0(K_{1,m,\mathfrak{P}}, (T_{\mathfrak{p}}^* \otimes \frac{J_\infty}{J_\infty^{[k+1]}})/\mathfrak{p}^m) \\
\theta_4 \uparrow & & \theta_5 \uparrow \\
colie(G)^{\otimes(k+j)} \otimes colie(\mathbb{G}_m)^{\otimes(-j)}/\mathfrak{p}^m & \xrightarrow{\rho_7} & colie(G)^{\otimes(k+j)} \otimes colie(\mathbb{G}_m)^{\otimes(-j)} \otimes \mathcal{O}_{K_{1,m,\mathfrak{P}}}/\mathfrak{p}^m
\end{array}$$

$$\begin{array}{ccc}
H^1(K_{1,m,\mathfrak{P}}, T^*_{\mathfrak{p}}(1)/\mathfrak{p}^m) & \xrightarrow{\rho_2} & H^1(K_{1,m,\mathfrak{P}}, T'^*_{\mathfrak{p}}(1)/\mathfrak{p}^m) \\
\downarrow \theta_2 & & \downarrow \theta_3 \\
H^1(K_{1,m,\mathfrak{P}}, (T^*_{\mathfrak{p}} \otimes \frac{J_\infty}{J_\infty^{[k+1]}})/\mathfrak{p}^m) & \xrightarrow{\rho_4} & H^1(K_{1,m,\mathfrak{P}}, (T'^*_{\mathfrak{p}} \otimes \frac{J_\infty}{J_\infty^{[k+1]}})/\mathfrak{p}^m) \\
\cup \log \chi_G \uparrow & & \cup \log \chi_G \uparrow \\
H^0(K_{1,m,\mathfrak{P}}, (T^*_{\mathfrak{p}} \otimes \frac{J_\infty}{J_\infty^{[k+1]}})/\mathfrak{p}^m) & \xrightarrow{\rho_6} & H^0(K_{1,m,\mathfrak{P}}, (T'^*_{\mathfrak{p}} \otimes \frac{J_\infty}{J_\infty^{[k+1]}})/\mathfrak{p}^m) \\
\theta_5 \uparrow & & \theta_6 \uparrow \\
colie(G)^{\otimes(k+j)} \otimes colie(\mathbb{G}_m)^{\otimes(-j)} \otimes \mathcal{O}_{K_{1,m,\mathfrak{P}}}/\mathfrak{p}^m & \xrightarrow{\rho_8} & colie(G)^{\otimes k} \otimes \mathcal{O}_{K_{1,m,\mathfrak{P}}}/\mathfrak{p}^m
\end{array}$$

Here are the maps used in the diagram: $\rho_1$ and $\rho_3$ are restriction maps; $\rho_5, \rho_7$ are inclusions; $\theta_1$ and $\theta_2$ is what we had called $\iota_1$ in lemma 5; $\theta_3$ is composition of $\iota_1$ with the map induced by multiplication by $\Omega_p^j \alpha^{-j} \mod \mathfrak{p}^m$ on $J_\infty/J_\infty^{[k+1]}/\mathfrak{p}^m$; $\rho_2$ is $\cup \zeta_m^{\otimes(-j)}$ and $\rho_4 = \rho_6$ are composition of map $\cup \zeta_m^{\otimes(-j)}$ with the map induced by multiplication by $\Omega_p^{\otimes j} \alpha^{-j} \mod \mathfrak{p}^m$; $\rho_8$ is tensoring with $\omega^{\otimes(-j)} \otimes (ds/1+s)^{\otimes j} \mod \mathfrak{p}^m$; $\theta_6$ is what we called $\iota_2$ in lemma 5 and induced from 3.3.5; similarly $\theta_4$ and $\theta_5$ are induced from 5.0.4. Note that $\rho_2$ is well defined since $\zeta_m$ (which is a fixed basis of $E_{\mathfrak{p}^{*m}}$) is rational over $K_{1,m}$ by definition; $\rho_4$ is well defined since

$$\sigma \Omega_p = (\psi(\sigma)/\chi(\sigma))\Omega_p, \forall \sigma \in Gal(\overline{\mathbf{Q}_p}/K_{\mathfrak{p}})$$

which can be deduced from the definition of $\Omega_p$ (see [6] 2.7 when $p$ splits and [14] theorem 6.1 when $p$ is nonsplit). Note our definition of $\Omega_p$ coincides with Yager's but it is actually inverse to that of de Shalit and Harrison. When $p$ splits, $\psi/\chi$ exactly describes the Galois action on $\zeta$. By the way it follows from these that $\Omega_p$ lies inside $\hat{\mathbf{Q}}_p^{\text{unr}}$ when $p$ splits and when $p$ does not split, $\Omega_p$ lies in the completion of $K_{\mathfrak{p}}(E_{p^\infty})$.

*Remark 12.* When $p$ splits we can easily makes sense of the map

$$\text{mult by } \Omega_p^j \alpha^{-j} \mod \mathfrak{p}^m : J_\infty/J_\infty^{[k+1]}/\mathfrak{p}^m \to J_\infty/J_\infty^{[k+1]}/\mathfrak{p}^m$$

since $\hat{\mathbf{Q}}_p^{\text{unr}} \subset B_{\text{crys}}$ and $\Omega_p$ is integral (even a $p$ unit).

**Lemma 18.** *The above diagram commutes.*

*Proof.* By definitions of various maps used in the diagram, the only nontrivial part of the commutativity is the lower right corner. This can be seen from Kato's pairing $colie(G) \otimes T_G \to J_\infty$ (take $G = \mathbb{G}_m$, or the $p$ divisible group associated to $\hat{E}_\mathfrak{p}$) under which

$$\frac{ds}{1+s} \otimes \epsilon \mapsto t, \quad \omega \otimes \xi \mapsto \Omega_p^{-1}\alpha t$$

**Lemma 19.** *(1): the maps $\rho_{2i}, i = 1, 2, 3, 4$ are isomorphisms;*
*(2): $ker(\rho_3)$ is killed by $p^e$ where $e$ is an integer independent of $m$ and for all $m$ big enough;*
*(3): the image of $\rho_7$ is the $Gal(K_{1,m,\mathfrak{P}}/K_\mathfrak{p})$ invariant part of the target space.*

*Proof.* The claim (1) for $\rho_{2i}, i = 1, \ldots, 4$ is clear, as $\Omega_p, \alpha$ are $\mathfrak{p}$-units.
(3) follows from exact sequence

$$0 \to \mathcal{O}_{K_{1,m,\mathfrak{P}}} \xrightarrow{\pi^m} \mathcal{O}_{K_{1,m,\mathfrak{P}}} \to \mathcal{O}_{K_{1,m,\mathfrak{P}}}/\pi^m \to 0$$

we have

$$0 \to \mathcal{O}_{K_\mathfrak{p}}/\pi^m \to (\mathcal{O}_{K_{1,m,\mathfrak{P}}}/\pi^m)^{Gal(K_{1,m,\mathfrak{P}}/K_\mathfrak{p})} \to H^1(Gal(K_{1,m,\mathfrak{P}}/K_\mathfrak{p}), \mathcal{O}_{K_{1,m,\mathfrak{P}}}) \tag{5.1.4}$$

Apply lemma 9, recall $K_{1,m,\mathfrak{P}}/K_\mathfrak{p}$ is not wildly ramified which implies that $i(\mu)$ defined in lemma 9 is 0 in our case. Hence

$$H^1(Gal(K_{1,m,\mathfrak{P}}/K_\mathfrak{p}), \mathcal{O}_{K_{1,m,\mathfrak{P}}}) = 0$$

Now we prove $ker(\rho_3)$ is killed by $p^e$ for some integer $e$ independent of $m$. To show this, we will use (a),(b),and (c) below:
(a): $ker(\rho_3) = H^1(Gal(K_{1,m,\mathfrak{P}}/K_\mathfrak{p}), \left((T_\mathfrak{p}^* \otimes J_\infty/J_\infty^{[k+1]})/\mathfrak{p}^m\right)^{Gal(\bar{\mathbf{Q}}_p/K_{1,m,\mathfrak{P}})})$.
(b): we have exact sequence

$$\mathcal{O}_{K_{1,m,\mathfrak{P}}}/\mathfrak{p}^m \to \left(B_\infty/J_\infty^{[k+1]}/\mathfrak{p}^m\right)^{Gal(\bar{\mathbf{Q}}_p/K_{1,m,\mathfrak{P}})} \to H^1(K_{1,m,\mathfrak{P}}, B_\infty/J_\infty^{[k+1]})$$

By obvious exact sequence and proposition 3(c), $H^1(K_{1,m,\mathfrak{P}}, B_\infty/J_\infty^{[k+1]})$ is killed by some $p^d$ which depends only on $k$. Hence we see

$$p^d(\mathcal{O}_{K_{1,m,\mathfrak{P}}}/\mathfrak{p}^m) = p^d \left(B_\infty/J_\infty^{[k+1]}/\mathfrak{p}^m\right)^{Gal(\bar{\mathbf{Q}}_p/K_{1,m,\mathfrak{P}})} \tag{5.1.5}$$

(c): by 5.1.5 and 5.0.4, we have

$$p^{N+d}ker(\rho_3) \cong H^1(Gal(K_{1,m,\mathfrak{P}}/K_{\mathfrak{p}}), \mathcal{O}_{K_{1,m,\mathfrak{P}}}/\mathfrak{p}^{m-d-N})$$

where $N$ is independent of $m$. Consider the following exact sequence:

$$H^1(G, \mathcal{O}_{K_{1,m,\mathfrak{P}}}) \to H^1(G, \mathcal{O}_{K_{1,m,\mathfrak{P}}}/\mathfrak{p}^{m-d-N}) \to H^2(G, \mathcal{O}_{K_{1,m,\mathfrak{P}}}) \tag{5.1.6}$$

where $G = Gal(K_{1,m,\mathfrak{P}}/K_{\mathfrak{p}})$. By lemma 9, $H^1(G, \mathcal{O}_{K_{1,m,\mathfrak{P}}}) = 0$. Since as $m \to \infty$, the fields $K_{1,m,\mathfrak{P}}$ has fixed ramification, we see

$$H^2(Gal(K_{1,m,\mathfrak{P}}/K_{\mathfrak{p}}), \mathcal{O}_{K_{1,m,\mathfrak{P}}}) = \mathcal{O}_{K_{\mathfrak{p}}}/Tr(\mathcal{O}_{K_{1,m,\mathfrak{P}}})$$

is killed by $p^a$ for some $a$ independent of $m$. Hence

$$p^a H^1(Gal(K_{1,m,\mathfrak{P}}/K_{\mathfrak{p}}), \mathcal{O}_{K_{1,m,\mathfrak{P}}}/\mathfrak{p}^{m-d-N}) = 0$$

Now take $e = N + d + a$ which is independent of $m$, we get

$$p^e ker(\rho_3) = 0$$

□

From now on we make the convention that $\forall \sigma \in Gal(K_{1,m}/K), \psi(\sigma)$ means $\psi(\bar{\sigma}) \mod \mathfrak{p}^m$ which is well defined for any $\bar{\sigma}$ which is a lift of $\sigma$ to $Gal(\bar{K}/K)$. When we talk about congruence, we meant for $\mathbb{C}_p$ and not necessarily for integers.

Now chasing diagram begining with $\overline{S_{0,m}}(u_{\mathfrak{P}}), \forall \mathfrak{P} = \upsilon^\eta$, $\eta \in B$. First

$$\begin{aligned} \rho_2 \circ \rho_1(\overline{S_{0,m}}(u_{\mathfrak{P}})) &= \rho_2(\sum_{\sigma \in Gal(K_{1,m,\mathfrak{P}}/K_{\mathfrak{p}})} \sigma \overline{S_m}(u_{\mathfrak{P}})) \\ &= \sum_{\sigma \in Gal(K_{1,m,\upsilon}/K_{\mathfrak{p}})} \bar{\psi}^j(\sigma\eta)\sigma\eta \hat{S_m}(u_\upsilon) \end{aligned} \tag{5.1.7}$$

where $\bar{\psi} = \psi/\chi$ and $\bar{\psi}^{-1}$ describes the Galois action on $T_{\mathfrak{p}^*}(E)$. Here $\hat{S_m}(u_{\mathfrak{P}})$ denotes $\overline{S_{0,m}}(u_{\mathfrak{P}}) \cup \zeta_m^{\otimes(-j)}$ which is nothing but the image of norm compatible units $\{u_{n,m}\}_n$ in $H^1(K_{1,m,\mathfrak{P}}, T_{\mathfrak{p}}(E)^{\otimes(-k)}(1)/\mathfrak{p}^m)$ defined just as we did in definition 4 when $j = 0$ and $\mathfrak{p}$ is a good reduction prime, except instead of using base field $K_{\mathfrak{p}}$, we make an unramified field extension and use $K_{0,m,\mathfrak{P}}$ as base field. Coleman power series for $\{u_{n,m}\}_n$ has been computed in lemma 13 and we can now apply Kato's "finite" version of $exp^*$ proposition 2 and in particular 3.3.11, which gives us equality in $H^1(K_{1,m,\mathfrak{P}}, (T'^*_{\mathfrak{p}} \otimes J_\infty/J_\infty^{[k+1]})/\mathfrak{p}^m)$:

$$\begin{aligned} c\eta\hat{S_m}(u_v) &= \sum_{\tau\in C} \rho(\frac{c}{(k-1)!}\pi^{-k}\omega^{\otimes k} \otimes \{(\frac{d}{\omega})^k \log(\phi^{-1}P^{\tau\eta}(\pi^{*(-m)}\lambda(t)))\}(\xi_1) \\ &= \sum_{\tau\in C} \rho(\frac{c}{(k-1)!}\pi^{-k}\omega^{\otimes k} \otimes \{(\frac{d}{dz})^k \log(\phi^{-1}\Theta^{\tau\eta}(\Omega/\pi^{*m} - z, L, \boldsymbol{\alpha}))\} \mid_{z=v_1/\pi}) \\ &= -12\sum_{\tau\in C} \rho(\frac{c}{(k-1)!}\pi^{-k}\pi^{*(-mk)}E_k^{\tau\eta}(\frac{\Omega}{\pi\pi^{*m}}, L, \boldsymbol{\alpha})\omega^{\otimes k}) \end{aligned} \tag{5.1.8}$$

where we use $C$ to denote $Gal(K(\mathfrak{f}\mathfrak{p}^{*m})/K(E_{\mathfrak{p}^{*m}}))$. We can argue completely the same way and get a formula for $c\sigma\hat{S_m}(u_v), \forall\sigma \in Gal(K_{1,m,v}/K_{\mathfrak{p}})$ and the only change in the corresponding result is to replace $E_k^{\tau\eta}$ in the above equation by $E_k^{\sigma\tau\eta}$. Note map $\rho$ in 3.3.11 is just $\theta_6 \circ \cup(\log\chi_G)$ here. Combining 5.1.7 and 5.1.8 and chasing the diagram, we have equality in $H^1(K_{1,m,\mathfrak{P}}, (T'^*_{\mathfrak{p}} \otimes J_\infty/J_\infty^{[k+1]})/\mathfrak{p}^m)$:

$$\begin{aligned} &c\rho_4 \circ \rho_3 \circ \theta_1(\overline{S_{0,m}}(u_{\mathfrak{P}})) \\ &\quad = c\theta_3 \circ \rho_2 \circ \rho_1(\overline{S_{0,m}}(u_{\mathfrak{P}})) \\ &\quad = -12 \sum_{\substack{\sigma\in Gal(K_{1,m,v}/K_{\mathfrak{p}}) \\ \tau\in C}} \rho(\Omega_p^j\alpha^{-j}\frac{c}{(k-1)!}\pi^{-k}\pi^{*(-mk)}\bar{\psi}^j(\sigma\eta)E_k^{\sigma\tau\eta}(\frac{\Omega}{\pi\pi^{*m}}, L, \boldsymbol{\alpha})\omega^{\otimes k}) \\ &\quad = -12 \sum_{\sigma\in Gal(K_{1,m,v}/K_{\mathfrak{p}})\times C} \rho(\Omega_p^j\alpha^{-j}\frac{c}{(k-1)!}\pi^{-k}\pi^{*(-mk)}\bar{\psi}^j(\sigma\eta)E_k^{\sigma\eta}(\frac{\Omega}{\pi\pi^{*m}}, L, \boldsymbol{\alpha})\omega^{\otimes k}) \end{aligned} \tag{5.1.9}$$

where the last equality holds becasue

$$\bar{\psi}(\tau) \equiv 1 \pmod{\mathfrak{p}^m}, \forall\tau \in Gal(K(\mathfrak{f}\mathfrak{p}^{*m})/K(E_{\mathfrak{p}^{*m}}))$$

It follows from this that in equality in $H^1(K_{1,m,\mathfrak{P}}, (T^*_{\mathfrak{p}} \otimes J_\infty/J_\infty^{[k+1]})/\mathfrak{p}^m)$

$$\begin{aligned} &c\rho_3 \circ \theta_1(\overline{S_{0,m}}(u_{\mathfrak{P}})) \\ &= \rho\{-12\Omega_p^j \alpha^{-j} \frac{c}{(k-1)!} \pi^{-k} \pi^{*(-mk)} \sum_{\sigma \in Gal(K_{1,m,v}/K_{\mathfrak{p}}) \times C} \bar{\psi}^j(\sigma\eta) E_k^{\sigma\eta}(\frac{\Omega}{\pi\pi^{*m}}, L, \boldsymbol{\alpha}) \varpi\} \\ &\stackrel{\text{def}}{=} \rho(A_{m,\eta}\varpi) \end{aligned} \tag{5.1.10}$$

Now observe $A_{m,\eta} \in \mathcal{O}_{K_{1,m,\mathfrak{P}}}/\pi^m$ is fixed by $Gal(K_{1,m,\mathfrak{P}}/K_{\mathfrak{p}})$ , since

$$\begin{aligned} &\tau\{\Omega_p^j \sum_{\sigma \in Gal(K_{1,m,v}/K_{\mathfrak{p}}) \times C} \bar{\psi}^j(\sigma\eta) E_k^{\sigma\eta}(\frac{\Omega}{\pi\pi^{*m}}, L, \boldsymbol{\alpha})\} \\ &\equiv \bar{\psi}(\tau)^j \Omega_p^j \bar{\psi}(\tau)^{-j} \sum_{\sigma \in Gal(K_{1,m,v}/K_{\mathfrak{p}}) \times C} \bar{\psi}^j(\sigma\eta\tau) E_k^{\sigma\tau\eta}(\frac{\Omega}{\pi\pi^{*m}}, L, \boldsymbol{\alpha})\} \mod \mathfrak{p}^m \quad (5.1.11) \\ &\equiv \{\Omega_p^j \sum_{\sigma \in Gal(K_{1,m,v}/K_{\mathfrak{p}}) \times C} \bar{\psi}^j(\sigma\eta) E_k^{\sigma\eta}(\frac{\Omega}{\pi\pi^{*m}}, L, \boldsymbol{\alpha})\} \mod \mathfrak{p}^m \end{aligned}$$

So by lemma 19 (3), $cA_{m,\eta}\varpi$ comes from $colie(G)^{\otimes(k+j)} \otimes colie(\mathbb{G}_m)^{\otimes(-j)}/\mathfrak{p}^m$ via $\rho_7$ and by the commutative diagram, together with lemma 19 (2), we have the following equality in $H^1(K_{\mathfrak{p}}, (T^*_{\mathfrak{p}} \otimes J_\infty/J_\infty^{[k+1]})/\mathfrak{p}^m)$:

$$p^e c\theta_1(\overline{S_{0,m}}(u_{\mathfrak{P}})) = p^e A_{m,\eta}\varpi \cup \log \chi_G \tag{5.1.12}$$

Define $A_m \in \mathcal{O}_{\mathfrak{p}}/\mathfrak{p}^m$ as follows:

$$\begin{aligned} A_m &= \sum_{\eta \in B} A_{m,\eta} \\ &= -12\Omega_p^j \alpha^{-j} \frac{c}{(k-1)!} \pi^{-k} \pi^{*(-mk)} \sum_{\sigma \in Gal(K(\mathfrak{f}\mathfrak{p}\mathfrak{p}^{*m})/K)} \bar{\psi}^j(\sigma) E_k^{\sigma}(\frac{\Omega}{\pi\pi^{*m}}, L, \boldsymbol{\alpha}) \quad (5.1.13) \\ &= -12\Omega_p^j \alpha^{-j} \frac{c}{(k-1)!} \pi^{*(-mk)} \sum_{\sigma \in Gal(K(\mathfrak{f}\mathfrak{p}\mathfrak{p}^{*m})/K)} \bar{\psi}^j(\sigma) E_k^{\sigma}(\frac{\Omega}{\pi^{*m}}, \mathfrak{p}L, \boldsymbol{\alpha}) \end{aligned}$$

We have by 5.1.12 equality in $H^1(K_{\mathfrak{p}}, (T^*_{\mathfrak{p}} \otimes J_\infty/J_\infty^{[k+1]})/\mathfrak{p}^m)$:

$$p^e c \sum_{\mathfrak{P}|\mathfrak{p}} \theta_1(\overline{S_{0,m}}(u_{\mathfrak{P}})) = p^e A_m\varpi \cup \log \chi_G \tag{5.1.14}$$

This settles part (1) of the theorem 4. Now our objective is to rewrite $A_m$ and to study its dependence on $m$. This gives us naturally the definition of $A$. Here we use crucially a congruence between $p$-adic period $\Omega_p$ with values of weight 1 Eisenstein series which in the case $p$ splits can be found in [3] ( 4.14 (41) ): (but be careful his notation $\mathfrak{f}, \Omega, \Omega_p$ are different from ours, $m$ has the same meaning though)

$$\Omega_p \equiv -\mathbf{N}(\mathfrak{f}\mathfrak{p}\mathfrak{p}^{*m}) E_1(\frac{\Omega}{\pi^{*m}}, \mathfrak{p}L) \pmod{\mathfrak{p}^{m-m_0}} \tag{5.1.15}$$

where $m_0$ is a constant independent of $m$.

We also need a congruence between different Eisenstein numbers as follows: (see [3] 3.4 (12), in the notation of [3], take $c = \mathbf{N}(\mathfrak{f}\mathfrak{p}\mathfrak{p}^{*m}), m = \mathfrak{f}\mathfrak{p}\mathfrak{p}^{*m}, \mu = \Omega/\pi\pi^{*m}$) $\forall\, 0 \leq j < k$,

$$(-2cE_1(\mu, L))^j E_k(\mu, L) \equiv (2c)^j E_{-j,k}(\mu, L) \pmod{\bar{f}\mathfrak{p}^m\mathfrak{p}^*} \tag{5.1.16}$$

or equivalently

$$(2c)^j E_{-j,k}(\frac{\Omega}{\pi^{*m}}, \mathfrak{p}L) \equiv (-2cE_1(\frac{\Omega}{\pi^{*m}}, \mathfrak{p}L))^j E_k(\frac{\Omega}{\pi^{*m}}, \mathfrak{p}L) \mod \mathfrak{p}^{m-j-k} \tag{5.1.17}$$

Using 5.1.15 and 5.1.17, we have:

$$\begin{aligned}
&c\Omega_p^j \bar{\psi}^j(\sigma) E_k^\sigma(\frac{\Omega}{\pi^{*m}}, \mathfrak{p}L) \\
&\quad\equiv c(\Omega_p^\sigma)^j E_k^\sigma(\frac{\Omega}{\pi^{*m}}, \mathfrak{p}L) \mod \mathfrak{p}^m \\
&\quad\equiv c(-\mathbf{N}(\mathfrak{f}\mathfrak{p}\mathfrak{p}^{*m}))^j (E_1^\sigma(\frac{\Omega}{\pi^{*m}}, \mathfrak{p}L))^j E_k^\sigma(\frac{\Omega}{\pi^{*m}}, \mathfrak{p}L) \mod \mathfrak{p}^{m-m_0-j-k} \\
&\quad\equiv c(-\mathbf{N}(\mathfrak{f}\mathfrak{p}\mathfrak{p}^{*m}))^j E_{-j,k}^\sigma(\frac{\Omega}{\pi^{*m}}, \mathfrak{p}L) \mod \mathfrak{p}^{m-m_0-j-k} \\
&\quad\equiv c(k-1)! \left(\frac{\sqrt{d_K}}{2\pi}\right)^j \left(\frac{\Omega}{\pi^{*m}}\right)^{-j-k} L(\bar{\psi}^{k+j}, \left(\frac{K(\mathfrak{f}\mathfrak{p}\mathfrak{p}^{*m})/K}{\mathfrak{c}}\right), k) \mod \mathfrak{p}^{m-m_0-j-k}
\end{aligned} \tag{5.1.18}$$

where $\mathfrak{c}$ is an integral ideal of $K$ such that $\sigma = \sigma_\mathfrak{c} \in Gal(K(\mathfrak{f}\mathfrak{p}\mathfrak{p}^{*m})/K)$ and the last equality follows from :

$$E_{-j,k}^{\sigma_\mathfrak{c}}(\frac{\Omega}{\pi^{*m}}, \mathfrak{p}L) = \psi(\mathfrak{c})^{-j-k} E_{-j,k}(\frac{\Omega}{\pi^{*m}}, \mathfrak{c}^{-1}\mathfrak{p}L)$$

and we apply 4.1.10 with $\mathfrak{m} = \mathfrak{f}\mathfrak{p}\mathfrak{p}^{*m}$ and $\Omega = \Omega/\pi^{*m} = \Omega_\infty/f\pi^{*m}$. Note also $\pi$ in $\sqrt{d_K}/2\pi$ in 5.1.18 is 3.14...

*Remark 13.* Recall $E_k(z, L, \boldsymbol{\alpha}) = \mathbf{N}\boldsymbol{\alpha} E_k(z, L) - E_k(z, \boldsymbol{\alpha}^{-1}L)$, one can apply similar arguments to $E_k(z, \psi(\boldsymbol{\alpha})\boldsymbol{\alpha}^{-1}L)$. The point is the integral and congruence property we used before for $E_{-j,k}(z, L)$ really requires that $L$ is taken to be the period lattice of some CM ellptic curve with good reduction at $\mathfrak{p}$. Hence we have to apply arguments above to $E_k(z, \psi(\boldsymbol{\alpha})\boldsymbol{\alpha}^{-1}L)$ and not $E_k(z, \boldsymbol{\alpha}^{-1}L)$, since $\psi(\boldsymbol{\alpha})\boldsymbol{\alpha}^{-1}L$ is the period lattice of $E^{\sigma_{\boldsymbol{\alpha}}}$, not $\boldsymbol{\alpha}^{-1}L$. This is more nontrivial when we consider $K$ with class number bigger than 1. Also when $E$ has bad reduction at $\mathfrak{p}$ as we are going to deal with in the next section, we replace $E$ by $E'$ which has good reduction at $\mathfrak{p}$ and apply above arguments to $E_k(z, L')$.

So we get similarly to 5.1.18

$$\begin{aligned}
&c\Omega_p^j \bar{\psi}^j(\sigma) E_k^\sigma(\frac{\Omega}{\pi^{*m}}, \boldsymbol{\alpha}^{-1}\mathfrak{p}L) \\
&\quad \equiv c\psi^k(\boldsymbol{\alpha})\bar{\psi}^j(\sigma)(\Omega_p)^j E_k^\sigma(\psi(\boldsymbol{\alpha})\frac{\Omega}{\pi^{*m}}, \mathfrak{p}L) \mod \mathfrak{p}^m \\
&\quad \equiv c\psi^k\bar{\psi}^{-j}(\boldsymbol{\alpha})\bar{\psi}^j(\sigma\sigma_{\boldsymbol{\alpha}})(\Omega_p)^j\bar{\psi}^j(\sigma) E_k^{\sigma\sigma_{\boldsymbol{\alpha}}}(\frac{\Omega}{\pi^{*m}}, \mathfrak{p}L) \mod \mathfrak{p}^m \\
&\quad \equiv c\psi^k\bar{\psi}^{-j}(\boldsymbol{\alpha})(k-1)!\left(\frac{\sqrt{d_K}}{2\pi}\right)^j \left(\frac{\Omega}{\pi^{*m}}\right)^{-j-k} L(\bar{\psi}^{k+j}, \left(\frac{K(\mathfrak{f}\mathfrak{p}\mathfrak{p}^{*m})/K}{\mathfrak{c}\boldsymbol{\alpha}}\right), k)
\end{aligned} \tag{5.1.19}$$

Let $W$ be a set of integral ideals prime to $\mathfrak{f}p$ such that $\{\left(\frac{K(\mathfrak{f}\mathfrak{p}\mathfrak{p}^{*m})/K}{\mathfrak{c}}\right); \mathfrak{c} \in W\}$ describes precisely $Gal(K(\mathfrak{f}\mathfrak{p}\mathfrak{p}^{*m})/K)$. Combining 5.1.13, 5.1.18 and 5.1.19, we get

$$\begin{aligned}
A_m &\equiv -12c\alpha^{-j}\left(\frac{\sqrt{d_K}}{2\pi}\right)^j \Omega^{-j-k}\pi^{*mj}\sum_{\mathfrak{c}\in W}\{\mathbf{N}\boldsymbol{\alpha} L(\bar{\psi}^{k+j}, \left(\frac{K(\mathfrak{f}\mathfrak{p}\mathfrak{p}^{*m})/K}{\mathfrak{c}}\right), k) \\
&\quad - \psi^k\bar{\psi}^{-j}(\boldsymbol{\alpha}) L(\bar{\psi}^{k+j}, \left(\frac{K(\mathfrak{f}\mathfrak{p}\mathfrak{p}^{*m})/K}{\mathfrak{c}\boldsymbol{\alpha}}\right), k)\} \mod \mathfrak{p}^{m-m_0-j-k} \\
&\equiv -12c\alpha^{-j}\left(\frac{\sqrt{d_K}}{2\pi}\right)^j \Omega^{-j-k}\pi^{*mj}(\mathbf{N}\boldsymbol{\alpha} - \psi^k\bar{\psi}^{-j}(\boldsymbol{\alpha}))L_{\mathfrak{f}\mathfrak{p}\mathfrak{p}^*}(\bar{\psi}^{k+j}, k) \mod \mathfrak{p}^{m-m_0-j-k}
\end{aligned} \tag{5.1.20}$$

Here $L_{\mathfrak{f}\mathfrak{p}\mathfrak{p}^*}(\bar{\psi}^{k+j}, s)$ is incomplete L-series. We see the only dependence of right hand side of 5.1.20 on $m$ is in the term $\pi^{*m}$. But since $\pi^* = \psi(\mathfrak{p}^*)$ is a $p$-adic unit, there exists $l$ such that $\pi^{*l} \equiv 1 \mod \mathfrak{p}$. Let

$$A = -12\alpha^{-j}\left(\frac{\sqrt{d_K}}{2\pi}\right)^{j}\Omega^{-j-k}(\mathbf{N}\boldsymbol{\alpha} - \psi^k\bar{\psi}^{-j}(\boldsymbol{\alpha}))L_{\mathfrak{f}\mathfrak{p}\mathfrak{p}^*}(\bar{\psi}^{k+j}, k) \tag{5.1.21}$$

which is independent of $m$. We have

$$cA \equiv A_{ml} \mod \mathfrak{p}^m$$

for all $m$ sufficently big (so that $(l-1)m > m_0 + j + k$), since $\pi^{*ml} \equiv 1 \pmod{\mathfrak{p}^m}$ by choice of $l$. This proves (2). $\square$

### *5.1.2 p is a bad reduction prime*

In this case the set of places of $F_\infty$ lying above $\mathfrak{p}$ is in 1-1 correspondence with $C \times D$ where $C = Gal(F/K)/Gal(F_v/K_\mathfrak{p})$ and $D = Gal(F_{0,M}/F)/Gal(F_{0,M,v}/F_v)$. Denote in this section $\hat{T}_\mathfrak{p} = T_\mathfrak{p}(E')^{\otimes(k+j)}(-j)$ and $\hat{T}'_\mathfrak{p} = T_\mathfrak{p}(E')^{\otimes(k+j)}$. Continue to use $T_\mathfrak{p}$ and $T'_\mathfrak{p}$ as defined in last section. We need the following commutative diagrams analogous to the the two diagrams we used in section 5.1:

$$\begin{array}{ccc}
H^1(K_\mathfrak{p}, T_\mathfrak{p}^*(1)/\mathfrak{p}^m) & \xrightarrow{\rho_1} & \sum_{\eta\in C} H^1(F_{m(\mathfrak{p}),m,v^\eta}, \hat{T}_\mathfrak{p}^{\,*}(1)/\mathfrak{p}^m) \\
\big\downarrow \theta_1 & & \big\downarrow \theta_2 \\
H^1(K_\mathfrak{p}, (T_\mathfrak{p}^* \otimes \frac{J_\infty}{J_\infty^{[k+1]}})/\mathfrak{p}^m) & \xrightarrow{\rho_3} & \sum_{\eta\in C} H^1(F_{m(\mathfrak{p}),m,v^\eta}, (\hat{T}_\mathfrak{p}^{\,*} \otimes \frac{J_\infty}{J_\infty^{[k+1]}})/\mathfrak{p}^m) \\
\cup\log\chi_G \big\uparrow & & \cup\log\chi_G = \cup\log\chi_{G'} \big\uparrow \\
H^0(K_\mathfrak{p}, (T_\mathfrak{p}^* \otimes \frac{J_\infty}{J_\infty^{[k+1]}})/\mathfrak{p}^m) & \xrightarrow{\rho_5} & \sum_{\eta\in C} H^0(F_{m(\mathfrak{p}),m,v^\eta}, (\hat{T}_\mathfrak{p}^{\,*} \otimes \frac{J_\infty}{J_\infty^{[k+1]}})/\mathfrak{p}^m) \\
\theta_4 \big\uparrow & & \theta_5 \big\uparrow \\
f^{-1}\{colie(G')^{\otimes(k+j)} \otimes colie(\mathbb{G}_m)^{\otimes(-j)}\}/\mathfrak{p}^m & \xrightarrow{\rho_7} & colie(G')^{\otimes(k+j)} \otimes colie(\mathbb{G}_m)^{\otimes(-j)} \otimes \mathcal{O}_m/\mathfrak{p}^m
\end{array}$$

$$\begin{array}{ccc}
\sum_{\eta\in C} H^1(F_{m,\upsilon^\eta}, \hat{T_{\mathfrak{p}}}^*(1)/\mathfrak{p}^m) & \xrightarrow{\rho_2} & \sum_{\eta\in C} H^1(F_{m,\upsilon^\eta}, \hat{T_{\mathfrak{p}}}^{'*}(1)/\mathfrak{p}^m) \\
\downarrow \theta_2 & & \downarrow \theta_3 \\
\sum_{\eta\in C} H^1(F_{m,\upsilon^\eta}, (\hat{T_{\mathfrak{p}}}^* \otimes \frac{J_\infty}{J_\infty^{[k+1]}})/\mathfrak{p}^m) & \xrightarrow{\rho_4} & \sum_{\eta\in C} H^1(F_{m,\upsilon^\eta}, (\hat{T_{\mathfrak{p}}}^{'*} \otimes \frac{J_\infty}{J_\infty^{[k+1]}})/\mathfrak{p}^m) \\
\cup\log\chi_G=\cup\log\chi_{G'} \uparrow & & \cup\log\chi_G \uparrow \\
\sum_{\eta\in C} H^0(F_{m,\upsilon^\eta}, (\hat{T_{\mathfrak{p}}}^* \otimes \frac{J_\infty}{J_\infty^{[k+1]}})/\mathfrak{p}^m) & \xrightarrow{\rho_6} & \sum_{\eta\in C} H^0(F_{m,\upsilon^\eta}, (\hat{T_{\mathfrak{p}}}^{'*} \otimes \frac{J_\infty}{J_\infty^{[k+1]}})/\mathfrak{p}^m) \\
\theta_5 \uparrow & & \theta_6 \uparrow \\
colie(G')^{\otimes(k+j)} \otimes colie(\mathbb{G}_m)^{\otimes(-j)} \otimes \mathcal{O}_m/\mathfrak{p}^m & \xrightarrow{\rho_8} & colie(G')^{\otimes k} \otimes \mathcal{O}_m/\mathfrak{p}^m
\end{array}$$

In the last line of the above two diagrams, we denote $\mathcal{O}_{m,\upsilon^\eta} = \mathcal{O}_{F_{m(\mathfrak{p}),m,\upsilon^\eta}}$, $\mathcal{O}_m = \sum_{\eta\in C} \mathcal{O}_{F_{m(\mathfrak{p}),m,\upsilon^\eta}}$ and denote $F_{m,\upsilon^\eta} = F_{m(\mathfrak{p}),m,\upsilon^\eta}$. The maps in the diagram are: $\rho_{2i-1}, i = 1, 2$ are $f \circ res$; $\rho_{2i-1}, i = 3, 4$ are $f \circ inclu$; $\rho_2$ is $\cup(\zeta_m^{'})^{\otimes(-j)}$ where $\zeta^{'}$ is the basis of $T_{\mathfrak{p}^*}(E^{'})$ we fixed in our study of coleman power series for elliptic units; $\rho_4 = \rho_6$ are composition of map $\cup {\zeta_m^{'}}^{\otimes(-j)}$ with the map induced by multiplication by $(\Omega_p^{'})^j(\alpha^{'})^{-j}$ and $\rho_8$ is tensoring with $\omega^{'\otimes(-j)} \otimes (ds/1+s)^{\otimes j}$; the maps $\theta_i$'s are the same as used before in section 5.1. Note $\rho_{2i}, i = 1, \ldots, 4$ are isomorphisms. By the same argument as in lemma 19 (2), we see that $ker(\rho_3)$ is killed by $p^e$ for some integer $e$ independent of $m$ for all $m$ big enough; and since $\hat{V_{\mathfrak{p}}} = V_{\mathfrak{p}}(\epsilon^{k+j})$, following the same proof as in lemma 19 (3), we can identify $p\times$ the source space of $\rho_7$ as the $\epsilon^{-k-j}$ part of $p\times$ the target space under the action of corresponding Galois group. $\forall\mathfrak{P} \mid \mathfrak{p}$ place of $K_\infty$, say $\mathfrak{P} = \upsilon^\tau, \tau \in D$,

$$\begin{aligned}
\rho_1(\overline{S_{0,m}}(u_{\mathfrak{P}})) &= \sum_{\sigma\in Gal(F_{m(\mathfrak{p}),m,\upsilon}/K_{\mathfrak{p}})\times C} f(\tau\sigma\overline{S_m}(u_\upsilon^{'})) \\
&= \sum_{\sigma\in Gal(F_{m(\mathfrak{p}),m,\upsilon}/K_{\mathfrak{p}})\times C} \epsilon^{k+j}(\tau\sigma)\tau\sigma\overline{S_m^{'}}(u_\upsilon^{'})
\end{aligned} \tag{5.1.22}$$

We also have $\rho_2(\sigma\overline{S_m^{'}}(u_\upsilon^{'})) = \bar{\psi}^{'j}(\sigma)\sigma\hat{S_m^{'}}(u_\upsilon^{'}), \forall\sigma$ where again we denote $\hat{S_m^{'}}(u_\upsilon^{'}) = \overline{S_m^{'}}(u_\upsilon^{'}) \cup (\zeta_m^{'})^{\otimes(-j)}$. Applying proposition 2 and 3.3.11 we get equality in $H^1(F_{m(\mathfrak{p}),m,\upsilon^\eta}, (\hat{T_{\mathfrak{p}}}^{'*}(1) \otimes J_\infty/J_\infty^{[k+1]})/\mathfrak{p}^m)$: $(\forall\eta \in C)$

$$c\hat{S_m^{'}}(\eta u_\upsilon^{'}) = -12\sum_{\gamma\in B} \rho\{\frac{c}{(k-1)!}\pi^{-m(\mathfrak{p})k}\pi^{*(-mk)}E_k^{\gamma\eta}(\frac{\Omega^{'}}{\pi^{m(\mathfrak{p})}\pi^{*m}}, L^{'}, \boldsymbol{\alpha})\omega^{'\otimes k}\} \tag{5.1.23}$$

where we denote $B$ as $Gal(K(\mathfrak{g}\mathfrak{p}^{*m})/K(E'_{\mathfrak{p}^{*m}}))$ which comes in with the coleman power series. Using the diagram we get equality in $\sum_{\eta\in C} H^1(F_{m(\mathfrak{p}),m,\upsilon^\eta}, (\hat{T}_{\mathfrak{p}}{}'^{*} \otimes J_\infty/J_\infty^{[k+1]})/\mathfrak{p}^m)$ (compare 5.1.9)

$$\begin{aligned}
&c\rho_4 \circ \rho_3 \circ \theta_1(\overline{S_{0,m}}(u_{\mathfrak{P}})) \\
&\quad = c\theta_3 \circ \rho_2 \circ \rho_1(\overline{S_{0,m}}(u_{\mathfrak{P}})) \\
&\quad = \theta_3(\sum_{\sigma\in Gal(F_{m(\mathfrak{p}),m,\upsilon}/K_{\mathfrak{p}})\times C} c\epsilon^{k+j}(\tau\sigma)\bar{\psi}'^{k+j}(\tau\sigma)\tau\sigma\hat{S'_m}(u'_\upsilon)) \\
&\quad = -12 \sum_{\substack{\sigma\in Gal(F_{m(\mathfrak{p}),m,\upsilon}/K_{\mathfrak{p}})\times C \\ \gamma\in B}} \rho\{(\Omega'_p)^j(\alpha')^{-j}\frac{c}{(k-1)!}(\pi')^{-m(\mathfrak{p})k}(\pi')^{*(-mk)} \\
&\quad \epsilon^{k+j}(\gamma\tau\sigma)\bar{\psi}'^{j}(\gamma\tau\sigma)E_k^{\gamma\tau\sigma}(\frac{\Omega'}{\pi'^{m(\mathfrak{p})}(\pi')^{*m}}, L', \boldsymbol{\alpha})\omega'^{\otimes k}\}
\end{aligned} \tag{5.1.24}$$

Hence we get equality in $\sum_{\eta\in C} H^1(F_{m(\mathfrak{p}),m,\upsilon^\eta}, (\hat{T}_{\mathfrak{p}}{}^{*} \otimes J_\infty/J_\infty^{[k+1]})/\mathfrak{p}^m)$ (compare 5.1.10):

$$\begin{aligned}
c\rho_3 \circ \theta_1(\overline{S_{0,m}}(u_{\mathfrak{P}})) &= \rho\{-12c\Omega_p'^{j}\alpha'^{-j}\frac{1}{(k-1)!}(\pi')^{-m(\mathfrak{p})k}(\pi')^{*(-mk)} \\
&\sum_{\sigma\in X}\epsilon^{k+j}(\tau\sigma)\bar{\psi}'^{j}(\tau\sigma)E_k^{\tau\sigma}(\frac{\Omega'}{\pi'^{m(\mathfrak{p})}(\pi')^{*m}}, L', \boldsymbol{\alpha})\varpi'\} \\
&\stackrel{\text{def}}{=} \rho(A_{m,\tau}\varpi')
\end{aligned} \tag{5.1.25}$$

where $X = Gal(F_{m(\mathfrak{p}),m,\upsilon}/K_{\mathfrak{p}}) \times B \times C$. It is easily seen that

$$\sigma A_{m,\tau} = \epsilon^{-k-j}(\sigma)A_{m,\tau},\ \forall\sigma \in Gal(F_{m(\mathfrak{p}),m,\upsilon}/K_{\mathfrak{p}}) \times C$$

Hence there is some element , namely:

$$f^{-1}(pA_{m,\tau}\varpi') \in \frac{f^{-1}\{colie(G')^{\otimes(k+j)} \otimes colie(\mathbb{G}_m)^{\otimes(-j)}\}}{\mathfrak{p}^m}$$

that maps to $pA_{m,\tau}$ under $\rho_7$. Then by left half of the commutative diagram and the fact that $ker(\rho_3)$ is killed by $p^e$, we get equality in $H^1(K_{\mathfrak{p}}, (T_{\mathfrak{p}}^* \otimes J_\infty/J_\infty^{[k+1]})/\mathfrak{p}^m)$

$$\sum_{\mathfrak{P}|\mathfrak{p}} p^{e+1}c\theta_1(\overline{S_{0,m}}(u_{\mathfrak{P}})) = \sum_{\tau\in D} f^{-1}(p^{e+1}A_{m,\tau}\varpi') \cup \log\chi_G \tag{5.1.26}$$

Now we rewrite $A_{m,\tau}$, applying the same method from last section to $E'$. Let $W$ be a set of integral ideals of $K$ prime to $\mathfrak{f}p$ such that $\{\left(\frac{K(\mathfrak{g}\mathfrak{p}^{m(\mathfrak{p})}\mathfrak{p}^{*m})/K}{\mathfrak{c}}\right); \mathfrak{c} \in W\}$ describes precisely $Gal(K(\mathfrak{g}\mathfrak{p}^{m(\mathfrak{p})}\mathfrak{p}^{*m})/K)$. The result is: (compare 5.1.20)

$$\begin{aligned}\sum_{\tau\in D} A_{m,\tau} &\equiv -12c\alpha'^{-j}\left(\frac{\sqrt{d_K}}{2\pi}\right)^{j}\Omega'^{(-j-k)}\pi'^{*mj}\\ &\sum_{\mathfrak{c}\in W}\epsilon^{k+j}(\mathfrak{c})\{\mathbf{N}\boldsymbol{\alpha}L(\bar{\psi}'^{k+j},\left(\frac{K(\mathfrak{g}\mathfrak{p}^{m(\mathfrak{p})}\mathfrak{p}^{*m})/K}{\mathfrak{c}}\right),k)\\ &-\psi'^{k}\bar{\psi}'^{-j}(\boldsymbol{\alpha})L(\bar{\psi}'^{k+j},\left(\frac{K(\mathfrak{g}\mathfrak{p}^{m(\mathfrak{p})}\mathfrak{p}^{*m})/K}{\mathfrak{c}\boldsymbol{\alpha}}\right),k)\}\mod \mathfrak{p}^{m-m_0-m(\mathfrak{p})(j+k)}\\ &\equiv -12c\alpha'^{-j}\left(\frac{\sqrt{d_K}}{2\pi}\right)^{j}\Omega'^{(-j-k)}\pi'^{*mj}(\mathbf{N}\boldsymbol{\alpha}-\psi^{k}\bar{\psi}^{-j}(\boldsymbol{\alpha}))L_{\mathfrak{g}\mathfrak{p}\mathfrak{p}^*}(\bar{\psi}^{k+j},k)\\ &\mod \mathfrak{p}^{m-m_0-m(\mathfrak{p})j-m(\mathfrak{p})k}\end{aligned} \tag{5.1.27}$$

here in the last equality we used the fact that $\psi' = \psi\epsilon$ and $\epsilon$ is a finite character, hence $\bar{\psi}'\epsilon = \bar{\psi}$ and $\epsilon^{-k-j}\psi'^{k}\bar{\psi}'^{-j} = \psi^{k}\bar{\psi}^{-j}$. Also when we apply the analogue of 5.1.17, we have to replace $m-j-k$ by $m-m(\mathfrak{p})j-m(\mathfrak{p})k$ since we are considering $E_{-j,k}(\Omega'/\pi'^{m(\mathfrak{p})}\pi'^{*m}, L')$ now. Recall $f^{*}\omega' = r\omega$, $f^{-1}(\Omega'_{\infty}) = r\Omega_{\infty}$, we see

$$\Omega'^{(-k-j)}f^{-1}(\varpi') = \Omega_{\infty}^{-j-k}g^{j+k}\varpi$$

The only dependence of the right hand side of 5.1.27 on $m$ is $\pi'^{*mj}$ term, and again we can apply the trick before to $\pi'$ and conclude

$$\sum_{\mathfrak{P}|\mathfrak{p}} exp^{*}(S_0(u_{\mathfrak{P}}))/\varpi = -12\alpha'^{-j}\left(\frac{\sqrt{d_K}}{2\pi}\right)^{j}\Omega_{\infty}^{(-j-k)}g^{k+j}(\mathbf{N}\boldsymbol{\alpha}-\psi^{k}\bar{\psi}^{-j}(\boldsymbol{\alpha}))L_{\mathfrak{g}\mathfrak{p}\mathfrak{p}^*}(\bar{\psi}^{k+j},k) \tag{5.1.28}$$

## 5.2 $p$ is nonsplit

**Theorem 5.** $\exists\alpha \in \mathcal{O}_K^{*}, d \in \mathcal{O}_{K_\mathfrak{p}}^{*}$

$$\sum_{\mathfrak{P}|\mathfrak{p}} exp^{*}(S_0(u_{\mathfrak{P}})) = d^{j}f^{k}\alpha^{-j}\left(\frac{\sqrt{d_K}}{2\pi}\right)^{-j}\Omega_{\infty}^{-k-j}(\mathbf{N}\boldsymbol{\alpha}-\psi^{k}\bar{\psi}^{-j}(\boldsymbol{\alpha}))L_{\mathfrak{f}\mathfrak{p}}(\bar{\psi}^{k+j},k) \tag{5.2.1}$$

*Proof.* Here we will just argue for the case $p$ is a good reduction prime for $E$. Notice in this case there is only one place of $K_\infty$ lying above $\mathfrak{p}$. We deal with the bad reduction case in exactly the same way as in the case $p$ splits and is a bad reduction prime, namely we apply what we do in section 5.1 for $F, E^{'}$, etc. We need the following two commutative diagrams:

$$
\begin{array}{ccc}
\tilde{H}^1(K_{\mathfrak{p}}, T_{\mathfrak{p}}^*(1)/\mathfrak{p}^n) & \xrightarrow{\rho_1} & \tilde{H}^1(K_{n,v}, T_{\mathfrak{p}}^*(1)/\mathfrak{p}^n) \\
 & & \downarrow{\scriptstyle \theta_2} \\
\downarrow{\scriptstyle \theta_1} & & H^1(K_{n,v}, (T_{\mathfrak{p}}^* \otimes \frac{J_\infty^{[k]}}{J_\infty^{[k+1]}})/\mathfrak{p}^n) \\
 & & \downarrow{\scriptstyle \theta_8} \\
H^1(K_{\mathfrak{p}}, (T_{\mathfrak{p}}^* \otimes \frac{J_\infty}{J_\infty^{[k+1]}})/\mathfrak{p}^n) & \xrightarrow{\rho_3} & H^1(K_{n,v}, (T_{\mathfrak{p}}^* \otimes \frac{J_\infty}{J_\infty^{[k+1]}})/\mathfrak{p}^n) \\
{\scriptstyle \cup \log \chi_G}\uparrow & & {\scriptstyle \cup \log \chi_G}\uparrow \\
H^0(K_{\mathfrak{p}}, (T_{\mathfrak{p}}^* \otimes \frac{J_\infty}{J_\infty^{[k+1]}})/\mathfrak{p}^n) & \xrightarrow{\rho_5} & H^0(K_{n,v}, (T_{\mathfrak{p}}^* \otimes \frac{J_\infty}{J_\infty^{[k+1]}})/\mathfrak{p}^n) \\
{\scriptstyle \theta_4}\uparrow & & {\scriptstyle \theta_7}\uparrow \\
colie(G)^{\otimes(k+j)} \otimes colie(\mathbb{G}_m)^{\otimes(-j)}/\mathfrak{p}^n & \xrightarrow{\rho_7} & colie(G)^{\otimes(k+j)} \otimes colie(\mathbb{G}_m)^{\otimes(-j)} \otimes \mathcal{O}_{K_{n,v}}/\mathfrak{p}^n
\end{array}
$$

$$
\begin{array}{ccc}
\tilde{H}^1(K_{n,v}, T_{\mathfrak{p}}^*(1)/\mathfrak{p}^n) & \xrightarrow{\rho_2} & \tilde{H}^1(K_{n,v}, T_{\mathfrak{p}}^{'*}(1)/\mathfrak{p}^n) \\
\downarrow{\scriptstyle \theta_2} & & \downarrow{\scriptstyle \theta_3} \\
H^1(K_{n,v}, (T_{\mathfrak{p}}^* \otimes \frac{J_\infty^{[k]}}{J_\infty^{[k+1]}})/\mathfrak{p}^n) & \xrightarrow{\rho_4} & H^1(K_{n,v}, (T_{\mathfrak{p}}^{'*} \otimes \frac{J_\infty^{[k]}}{J_\infty^{[k+1]}})/\mathfrak{p}^n) \\
{\scriptstyle \cup \log \chi_G}\uparrow & & {\scriptstyle \cup \log \chi_G}\uparrow \\
H^0(K_{n,v}, (T_{\mathfrak{p}}^* \otimes \frac{J_\infty^{[k]}}{J_\infty^{[k+1]}})/\mathfrak{p}^n) & \xrightarrow{\rho_6} & H^0(K_{n,v}, (T_{\mathfrak{p}}^{'*} \otimes \frac{J_\infty^{[k]}}{J_\infty^{[k+1]}})/\mathfrak{p}^n) \\
{\scriptstyle \theta_5}\uparrow & & {\scriptstyle \theta_6}\uparrow \\
colie(G)^{\otimes(k+j)} \otimes colie(\mathbb{G}_m)^{\otimes(-j)} \otimes \mathcal{O}_{K_{n,v}}/\mathfrak{p}^n & \xrightarrow{\rho_8} & colie(G)^{\otimes k} \otimes \mathcal{O}_{K_{n,v}}/\mathfrak{p}^n
\end{array}
$$

*Remark 14.* Note the subtle difference between this diagram and the one used before for $p$ splits case. By our calculation in the case $j = 0$, e.g 4.2.4 and 4.2.2, we see $exp^*(S_n(u_v))/\varpi \in K_{n,v}$ has $\mathfrak{p}$-valuation bounded independent of $n$, where $u$ is a generator of module of elliptic units. It follows we have (compare lemma 5)

**Lemma 20.** *Let $\tilde{H}^1(K_{n,v}, T^{\otimes(-k)}(1)/\pi^m)$ be the submodule of $H^1(K_{n,v}, T^{\otimes(-k)}(1)/\pi^m)$ generated by (image of) Soule elements $\hat{S_{n,m}}(u) \in H^1(K_{n,v}, T^{\otimes(-k)}(1)/\pi^m)$, $\forall u \in \mathfrak{C}_\infty$. Then there $\exists c$ independent of $n, m$, and $\exists \theta$, which make the following diagram commutative, for all $m \geq n$:*

$$
\begin{array}{ccc}
H^1(K_{n,v}, T^{\otimes(-k)}(1)/\pi^m) & \xrightarrow{c\iota_1} & H^1(K_{n,v}, (T^{\otimes(-k)} \otimes J_\infty/J_\infty^{[k+1]})/\pi^m) \\
\uparrow & & \uparrow \theta_8 \\
\tilde{H}^1(K_{n,v}, T^{\otimes(-k)}(1)/\pi^m) & \xrightarrow{\theta} & H^1(K_{n,v}, (T^{\otimes(-k)} \otimes J_\infty^{[k]}/J_\infty^{[k+1]})/\pi^m) \\
\downarrow \lambda & \nearrow \rho & \uparrow \log \chi_G \\
colie(G)^{\otimes k} \otimes \mathcal{O}_{K_{n,v}}/\pi^m & \xrightarrow{\iota_2} & H^0(K_{n,v}, (T^{\otimes(-k)} \otimes J_\infty^{[k]}/J_\infty^{[k+1]})/\pi^m)
\end{array}
$$

*where $\lambda((\hat{S_{n,m}}(u)) = c/(k-1)!\,\pi^{-nk}\omega^{\otimes k} \otimes \{(\frac{d}{\omega})^k \log(\phi^{-n}(g_{u,\xi}))\}(\xi_n)$ and $\theta_8$ is the natural map induced by $J_\infty^{[k]}/J_\infty^{[k+1]} \to J_\infty/J_\infty^{[k+1]}$.*

Here the map $\lambda$ can be thought of as the finite version of dual exponential map. We will only use the case $n = m$ of the above lemma in the following. If we were able to prove analogue of proposition 3 (c) in the case $p$ is nonsplit (i.e $c(n)$ in lemma 5 is independent of $n$), we do not have to use $\tilde{H}^1(K_\mathfrak{p}, T_\mathfrak{p}^*(1)/\pi^m)$. On the other hand, if $c(n) \to \infty$ (e.g grow like $p^n$), then it implies that the index of submodule generated by Soule type elements comming from elliptic units inside $H^1(K_{n,v}, T_\mathfrak{p}^*(1))/H_f^1(K_{n,v}, T_\mathfrak{p}^*(1))$ also tends to infinity as $n \to \infty$ and then we might not be able to define $\lambda$ on the whole $H^1(K_{n,v}, T_\mathfrak{p}^*(1)/\pi^m)$. Anyway we conclude following identity (compare with 3.3.11):

$$
\begin{aligned}
&\rho(\frac{c}{(k-1)!}\pi^{-nk}\omega^{\otimes k} \otimes \{(\frac{d}{\omega})^k \log(\phi^{-n}(g_{u,\xi}))\}(\xi_n)) \\
&\quad = \theta((\hat{S_{n,n}}(u)), \text{ in } H^1(K_{n,v}, (T^{\otimes(-k)} \otimes J_\infty^{[k]}/J_\infty^{[k+1]})/\pi^n)
\end{aligned} \tag{5.2.2}
$$

In the diagrams at the beginning of this section, $\theta_3$ is the composition of $\theta$ with multiplication by $\Omega_p^j\alpha^{-j}$ on $J_\infty^{[k]}/J_\infty^{[k+1]}/\mathfrak{p}^n$. $\theta_2$ is induced from $\theta$ to make the following diagram commutative:

$$\begin{array}{ccc} \tilde{H}^1(K_{n,v},T^{\otimes(-k-j)}(1+j)/\pi^n) & \xrightarrow{\theta_2} & H^1(K_{n,v},(T^{\otimes(-k-j)}(j)\otimes J_\infty^{[k]}/J_\infty^{[k+1]})/\pi^n) \\ \downarrow{\scriptstyle \cup\xi_n^{\otimes j}\epsilon_n^{\otimes(-j)}} & & \downarrow{\scriptstyle \cup\xi_n^{\otimes j}\epsilon_n^{\otimes(-j)}} \\ \tilde{H}^1(K_{n,v},T^{\otimes(-k)}(1)/\pi^n) & \xrightarrow{\theta} & H^1(K_{n,v},(T^{\otimes(-k)}\otimes J_\infty^{[k]}/J_\infty^{[k+1]})/\pi^n) \end{array}$$

$\rho_2=\cup\xi_n^{\otimes j}\epsilon_n^{\otimes(-j)}$; $\rho_4=\rho_6$ are the composition of $\cup\xi_n^{\otimes j}\epsilon_n^{\otimes(-j)}$ with the map induced by multiplication by $\Omega_p^j\alpha^{-j}$; $\rho_8$ is tensoring with $\omega^{\otimes(-j)}\otimes(ds/1+s)^{\otimes j}$. It is clear that $\rho_2$ and $\rho_8$ are isomorphisms while the kernel and cokernel of $\rho_4$ and $\rho_6$ are killed by some fixed power $p^d$ of $p$ (can take $d$ to be the $\mathfrak{p}$ valuation of $\Omega_p^j$). $\rho_1,\rho_3,\rho_5$ are restriction maps and $\rho_7$ is inclusion. $\theta_1$ is induced by $\mathbf{Z}_p(1)\hookrightarrow J_\infty$ (the same map we used when $p$ splits); $\theta_i, i=4,5,6,7$ are maps in 3.3.5 and 5.0.4.

Notice now we can also make sense of the map:

$$\text{mult by } \Omega_p^j \mod \mathfrak{p}^m : J_\infty^{[k]}/J_\infty^{[k+1]}/\mathfrak{p}^m \to J_\infty^{[k]}/J_\infty^{[k+1]}/\mathfrak{p}^m$$

as it is well known that $\forall k\in\mathbf{Z}, J_\infty^{[k]}/J_\infty^{[k+1]}$ has an $\mathcal{O}_{\mathbb{C}_p}$ module structure. Define $b_k$ so that the image of following map is $p^{b_k}\mathcal{O}_{\mathbb{C}_p}$.

$$J_\infty^{[k]}/J_\infty^{[k+1]} \xrightarrow{\lambda} \mathcal{O}_{\mathbb{C}_p}, \quad (\text{image of } t^k)\mapsto 1$$

**Lemma 21.** *(1): $ker(\rho_1)$ and $ker(\rho_3)$ are killed by $p^e$ where $e$ is an integer independent of $n$;*
*(2): after multiplying by $p^b$ where $b$ is independent of $n$, we can identify the source space of $\rho_7$ with the $Gal(K_{n,v}/K_\mathfrak{p})$ invariant part of the target space of $\rho_7$.*

*Proof.* Recall $T_\mathfrak{p}^*=\mathcal{O}_\mathfrak{p}(\psi_p^{-k-j}\chi_{\text{cyclo}}^j)$. Since $Gal(\bar{\mathbf{Q}}_p/K_{n,v})$ acts trivially on $T_\mathfrak{p}^*(1)/\mathfrak{p}^n$, we get

$$ker(\rho_1)=H^1(Gal(K_{n,v}/K_\mathfrak{p}),\mathcal{O}_\mathfrak{p}(\eta)/\mathfrak{p}^n)$$

where $\eta=\psi_p^{-k-j}\chi_{\text{cyclo}}^{1+j} \mod \mathfrak{p}^n$ which is a nontrivial character of $Gal(K_{n,v}/K_\mathfrak{p})$. Clearly this is killed by a power of $p$ independent of $n$ (only depends on $k,j,\psi$).

Hence (1) for $ker(\rho_1)$ follows.

We know

$$ker(\rho_3) = H^1(Gal(K_{n,v}/K_{\mathfrak{p}}), H^0(K_{n,v}, T^*_{\mathfrak{p}} \otimes J_\infty/J_\infty^{[k+1]}/\mathfrak{p}^n)) \tag{5.2.3}$$

Hence by obvious exact sequence, to prove the claim for $ker(\rho_3)$, we just need to prove that some power of $p$ independent of $n$ kills

$$H^1(Gal(K_{n,v}/K_{\mathfrak{p}}), H^0(K_{n,v}, T^*_{\mathfrak{p}} \otimes J_\infty^{[i]}/J_\infty^{[i+1]}/\mathfrak{p}^n))$$

for each $1 \leq i \leq k$. Recall $J_\infty^{[i]}/J_\infty^{[i+1]} \cong p^{b_i}\mathcal{O}_{\mathbb{C}_p}$ under the map $\lambda$ above. By definition, the action of $Gal(\bar{\mathbf{Q}}_p/K_{n,v})$ on $T^*_{\mathfrak{p}} \otimes J_\infty^{[i]}/J_\infty^{[i+1]}/\mathfrak{p}^n$ is trivial. Hence as $Gal(\bar{\mathbf{Q}}_p/K_{n,v})$ modules

$$T^*_{\mathfrak{p}} \otimes J_\infty^{[i]}/J_\infty^{[i+1]}/\mathfrak{p}^n \cong p^{b_i}\mathcal{O}_{\mathbb{C}_p}/\mathfrak{p}^n$$

By the following obvious exact sequence (note $\mathcal{O}_{\mathbb{C}_p}/\mathfrak{p}^n \cong \mathcal{O}_{\bar{\mathbf{Q}}_p}/\mathfrak{p}^n$):

$$H^0(K_{n,v}, p^{b_i}\mathcal{O}_{\mathbb{C}_p})/\mathfrak{p}^n \hookrightarrow H^0(K_{n,v}, T^*_{\mathfrak{p}} \otimes J_\infty^{[i]}/J_\infty^{[i+1]}/\mathfrak{p}^n) \to H^1(K_{n,v}, p^{b_i}\mathcal{O}_{\bar{\mathbf{Q}}_p})$$

together with $pH^1(K_{n,v}, \mathcal{O}_{\bar{\mathbf{Q}}_p}) = 0$ ([11] theorem 3), we conclude that

$$pH^0(K_{n,v}, T^*_{\mathfrak{p}} \otimes J_\infty^{[i]}/J_\infty^{[i+1]}/\mathfrak{p}^n) = p^{b_i+1}\mathcal{O}_{K_{n,v}}/\mathfrak{p}^n \tag{5.2.4}$$

as $H^0(L, \mathcal{O}_{\mathbb{C}_p}) = \mathcal{O}_L$ for any local field $L$, by [7] 2.2.4 (3).

Hence we just need to prove

$$p^d H^1(Gal(K_{n,v}/K_{\mathfrak{p}}), \mathcal{O}_{K_{n,v}}/\mathfrak{p}^n(\phi_i)) = 0, \forall 1 \leq i \leq k \tag{5.2.5}$$

for some $d$ independent of $n$, where $\phi_i = \psi^{-k-j}\chi_{\text{cyclo}}^{i+j} \mod \mathfrak{p}^n$ which is a nontrivial character of $Gal(K_{n,v}/K_{\mathfrak{p}})$. This follows immediately by applying lemma 11 to $K_{n,v}/K_{\mathfrak{p}}$. Hence (1) for $ker(\rho_3)$.

Finally, since we have the following obvious exact sequence:

$$0 \to \mathcal{O}_{\mathfrak{p}}/\mathfrak{p}^n \to (\mathcal{O}_{K_{n,v}}/\mathfrak{p}^n)^G \to H^1(G, \mathcal{O}_{K_{n,v}})$$

where $G = Gal(K_{n,v}/K_{\mathfrak{p}})$, we see to prove (2) in the lemma 21, we just need to prove that $p^2$ kills $H^1(Gal(K_{n,v}/K_{1,v}), \mathcal{O}_{K_{n,v}})$. This follows from lemma 10 we have proved before. □

*Remark 15.* Sen's result (i.e lemma 9) exactly deal with wildly ramified field extension case (as we are in now). The unramified or tamely ramified case is trivial.

Now we are ready to prove theorem 5, using the 2 commutative diagrams at the beginning of this section. By definition, we see

$$\begin{aligned}\rho_2\circ\rho_1(\overline{S_{0,n}}(u_v)) &= \rho_2(\sum_{\sigma\in Gal(K_{n,v}/K_{\mathfrak{p}})}\sigma\overline{S_{n,n}}(u_v))\\ &= \sum_{\sigma\in Gal(K_n/K)}\bar{\psi}^j(\sigma)\sigma\hat{S_{n,n}}(u_v)\end{aligned}\tag{5.2.6}$$

where we use $\hat{S_{n,n}}(u_v)$ to denote $\overline{S_{n,n}}(u_v)\cup\xi^{\otimes j}\epsilon^{-j}$ which is just the image of norm compatible units $\{u_n\}$ in $H^1(K_{n,v},T_{\mathfrak{p}}(E)^{\otimes(-k)}(1)/\mathfrak{p}^n)$ defined just as we did in definition 4 case $j=0$ and $\mathfrak{p}$ is a good reduction prime. Again using proposition 2 and 5.2.2, we have equality in $H^1(K_{n,v},(T_{\mathfrak{p}}'^{*}\otimes J_\infty^{[k]}/J_\infty^{[k+1]})/\mathfrak{p}^n)$:

$$\begin{aligned}\theta(\sigma\hat{S_{n,n}}(u_v)) &= \sum_{\tau\in C}\rho(\frac{c}{(k-1)!}\pi^{-nk}\omega^{\otimes k}\otimes\{(\frac{d}{\omega})^k\log(\phi^{-n}P(\lambda(t))^{\tau\sigma})\}(\xi_n))\\ &= \sum_{\tau\in C}\rho((\frac{c}{(k-1)!}\pi^{-nk}\omega^{\otimes k}\otimes\{(\frac{d}{\omega})^k\log(\phi^{-n}\Theta^{\tau\sigma}(\Omega-z,L,\boldsymbol{\alpha}))\}\mid_{z=v_n/\pi^n})\\ &= -12\sum_{\tau\in C}\rho((\frac{c}{(k-1)!}\pi^{-nk}E_k^{\tau\sigma}(\frac{\Omega}{\pi^n},L,\boldsymbol{\alpha})\omega^{\otimes k}),\ \forall\sigma\in Gal(K_n/K)\end{aligned}\tag{5.2.7}$$

Here we use $C$ to denote $Gal(K(E_{\mathfrak{f}\mathfrak{p}^n})/K(E_{\mathfrak{p}^n}))$. By 5.2.6 and 5.2.7 and using the above diagram, we get in $H^1(K_{n,v},(T_{\mathfrak{p}}'^{*}\otimes J_\infty^{[k]}/J_\infty^{[k+1]})/\mathfrak{p}^n)$

$$\begin{aligned}&\theta_3\circ\rho_2\circ\rho_1(\overline{S_{0,n}}(u_{\mathfrak{P}}))\\ &= \Omega_p^j\alpha^{-j}\sum_{\sigma\in Gal(K_n/K)}\bar{\psi}^j(\sigma)\theta_3(\sigma\hat{S_{n,n}}(u_v))\\ &= -12\sum_{\substack{\sigma\in Gal(K_n/K)\\ \tau\in C}}\rho\{c\Omega_p^j\alpha^{-j}\frac{1}{(k-1)!}\pi^{-nk}\bar{\psi}^j(\tau\sigma)E_k^{\tau\sigma}(\frac{\Omega}{\pi^n},L,\boldsymbol{\alpha})\omega^{\otimes k}\}\\ &= -12\sum_{\sigma\in Gal(K_{\mathfrak{f}\mathfrak{p}^n}/K)}\rho\{c\Omega_p^j\alpha^{-j}\frac{1}{(k-1)!}\pi^{-nk}\bar{\psi}^j(\sigma)E_k^{\sigma}(\frac{\Omega}{\pi^n},L,\boldsymbol{\alpha})\omega^{\otimes k}\}\end{aligned}\tag{5.2.8}$$

Since $ker(\rho_4)$ is killed by some $p^d$, with $d$ independent of $n$, it follows we have equality in $H^1(K_{n,v}, (T^*_{\mathfrak{p}} \otimes J_\infty/J_\infty^{[k+1]})/\mathfrak{p}^n)$,

$$\begin{aligned}
&\rho_3 \circ \theta_1(p^d\overline{S_{0,n}}(u_v)) \\
&\quad = \theta_8 \circ \theta_2 \circ \rho_1(p^d\overline{S_{0,n}}(u_v)) \\
&\quad = \rho\{-12cp^d\Omega_p^j\alpha^{-j}\frac{1}{(k-1)!}\pi^{-nk}\sum_{\sigma\in Gal(K_{\mathfrak{f}\mathfrak{p}^n}/K)}\bar{\psi}^j(\sigma)E_k^\sigma(\frac{\Omega}{\pi^n}, L, \boldsymbol{\alpha})\varpi\} \qquad (5.2.9)\\
&\quad \stackrel{\text{def}}{=} \rho(p^d A_n\varpi)
\end{aligned}$$

Since $\sigma\Omega_p = \bar{\psi}(\sigma)\Omega_p, \forall\sigma \in Gal(K_n/K)$, we see that $A_n \in \mathcal{O}_{K_{n,v}}/\mathfrak{p}^n$ is fixed by $Gal(K_n/K)$. Hence by lemma 21, (2), we see $p^bA_n\varpi$ is in the image of $\rho_7$ and by the commutative diagram, and also lemma 21 (1), we have equality in $H^1(K_{\mathfrak{p}}, (T^*_{\mathfrak{p}} \otimes J_\infty/J_\infty^{[k+1]})/\mathfrak{p}^n)$:

$$p^{e+d+b}\theta_1(c\overline{S_{0,n}}(u_v)) = p^{e+d+b}A_n\varpi \cup \log\chi_G \qquad (5.2.10)$$

**Lemma 22.** *Let $\mathfrak{p} \mid p$ and $p$ does not split. Suppose $(\mathfrak{f}, \mathfrak{p}) = 1$ and $v$ be a primitive $\mathfrak{f}\mathfrak{p}^n$ torsion point of $L$ which is the period lattice of $E$ with good reduction at $\mathfrak{p}$. Then $\pi^nE_1(v, L)$ and $\pi^jE_{-j,k}(v, L)$ are $\mathfrak{p}$ integral for $0 < j < k, k > 1$, where $\pi = \psi(\mathfrak{p})$ is a uniformizer of $K$ at $\mathfrak{p}$.*

*Remark 16.* Analogue of this in the case $p$ splits can be found in [3] 3.3 (i.v). We need this to make sense of congruence like 5.1.17.

*Proof.* First prove that $E_k(v, L)$ is $\mathfrak{p}$ integral, $\forall k \geq 3$. The point is one can easily verify from the definition that $E_k(v, L) = (-1)^k\mathcal{P}^{(k-2)}(v, L)$. But since $v$ is a primitive $\mathfrak{f}\mathfrak{p}^n$ torsion point of $L$, $\varepsilon(v, L)$ is not in the kernel of reduction modulo $\mathfrak{p}$ map, where

$$\varepsilon(-, L) : \mathbb{C}/L \to E, z \mapsto (x, y), x = \mathcal{P}(z), y = \mathcal{P}'(z)$$

Hence its $x$ and $y$ coordinates are integral at $\mathfrak{p}$ and by an easy induction we get $\mathcal{P}^{(k-2)}(v, L)$ is also integral at $\mathfrak{p}$, for all $k \geq 3$. Note this argument holds for $p$ splits as well as for $p$ which is nonsplit. Now when $p$ nonsplits, in [14], Yager proved by a careful study of congruence arising from formal group $\hat{E}_{\mathfrak{p}}$ that $p^nE_1(v, L)$ is $p$ integral (see cor 3.6 of [14]). To deduct this from his cor 3.6 just notice his $u_n = -2\mathcal{P}(\Omega_\infty/\pi^n)/\mathcal{P}'(\Omega_\infty/\pi^n)$ hence is $p$ integral.
For $k=2$ case, we can use a trick by de Shalit [3] P60, (10). Finally we use 4.1.9 to conclude $\pi^j E_{-j,k}(v,L)$ is $p$ integral, for $0 \le j < k, k > 1$. □

Yager has proved stronger results: Let $v_n = \Omega_\infty/f\pi^n = \Omega/\pi^n$, then

$$\bar{f}\pi^n E_1(v_n, L) \equiv \bar{f}\pi^{n+1} E_1(v_{n+1}, L) \mod \mathfrak{p}^{n-1}$$

and so $\gamma \overset{\text{def}}{=} \lim_{n\to\infty} \bar{f}\pi^n E_1(v_n, L)$ exists. In [14] prop 4.5 and th 5.5 Yager studied Galois action and $\mathfrak{p}$ valuation of $\gamma$. It turns out they coincide with those of $\Omega_p$. Thanks to his result, we know now $\exists d \in \mathcal{O}^*_{K_\mathfrak{p}}$ such that $\Omega_p = d\gamma$ and hence

$$\Omega_p \equiv d\bar{f}\pi^n E_1(\Omega/\pi^n, L) \mod \mathfrak{p}^{n-1} \tag{5.2.11}$$

By 4.1.9 and lemma 22, we get:

$$(\bar{f}\pi^n E_1(v_n, L))^j E_k(v_n, L) \equiv (\bar{f}\pi^n)^j E_{-j,k}(v_n, L) \mod \mathfrak{p}^n$$

$$\begin{aligned} &\frac{c}{(k-1)!}\pi^{-nk}\Omega_p^j \bar{\psi}^j(\sigma_\mathfrak{b}) E_k^{\sigma_\mathfrak{b}}(v_n, L) \\ &\equiv \frac{c}{(k-1)!}\pi^{-nk}(\Omega_p^j)^{\sigma_\mathfrak{b}} E_k^{\sigma_\mathfrak{b}}(v_n, L) \mod \mathfrak{p}^n \\ &\equiv \frac{c}{(k-1)!}\pi^{-nk} d^j (\bar{f}\pi^n)^j E_{-j,k}(\psi(\mathfrak{b})v_n, L) \mod \mathfrak{p}^n \\ &\equiv d^j f^{-j} c \left(\frac{\sqrt{d_K}}{2\pi}\right)^{-j} \Omega^{-k-j} L_{\mathfrak{f}\mathfrak{p}^n}(\bar{\psi}^{k+j}, \sigma_\mathfrak{b}, k) \mod \mathfrak{p}^n \end{aligned} \tag{5.2.12}$$

where in the last congruence we use 4.1.10, with $\Omega = v_n = \Omega_\infty/f\pi^n$, $\mathfrak{m} = \mathfrak{f}\mathfrak{p}^n$ and $\mathfrak{c} = \mathfrak{b}$. Now combining with 5.2.9, we get

$$\begin{aligned} A_n &\equiv c\alpha^{-j} d^j f^{-j} \sum_{\sigma \in Gal(K(\mathfrak{f}\mathfrak{p}^n)/K)} \left(\frac{\sqrt{d_K}}{2\pi}\right)^{-j} \Omega^{-k-j} \mathbf{N}(\boldsymbol{\alpha}) L_{\mathfrak{f}\mathfrak{p}^n}(\bar{\psi}^{k+j}, \sigma, k) \\ &\quad - (\psi^k \bar{\psi}^{-j})(\boldsymbol{\alpha}) L_{\mathfrak{f}\mathfrak{p}^n}(\bar{\psi}^{k+j}, \sigma\sigma_{\boldsymbol{\alpha}}, k) \mod \mathfrak{p}^n \\ &\equiv c\alpha^{-j} d^j f^k \left(\frac{\sqrt{d_K}}{2\pi}\right)^{-j} \Omega_\infty^{-k-j} (\mathbf{N}(\boldsymbol{\alpha}) - (\psi^k \bar{\psi}^{-j})(\boldsymbol{\alpha})) L_{\mathfrak{f}\mathfrak{p}}(\bar{\psi}^{k+j}, k) \mod \mathfrak{p}^n \end{aligned} \tag{5.2.13}$$

Call the right hand side to be $cA$ which is independent of $n$. By 5.2.10, we have:

$$exp^*(p^{e+d+b}cS_0(u_v))/\varpi \equiv p^{e+d+b}cA \mod \mathfrak{p}^n, \forall n > 0$$

where $b, c, d, e$ are constants independent of $n$. Hence we conclude:

$$exp^*(S_0(u_v))/\varpi = d^j\alpha^{-j}f^k\left(\frac{\sqrt{d_K}}{2\pi}\right)^{-j}\Omega_\infty^{-k-j}(\mathbf{N}(\boldsymbol{\alpha}) - (\psi^k\bar{\psi}^{-j})(\boldsymbol{\alpha}))L_{\mathfrak{f}\mathfrak{p}}(\bar{\psi}^{k+j}, k) \tag{5.2.14}$$

□

# CHAPTER 6
# SHAFAREVICH-TATE GROUP

Recall $p$-adic realization for $\mathcal{M}_{k,j}$ is

$$V_p = \begin{cases} \mathbf{Q}_p(\varphi_{\mathfrak{p}}) \oplus \mathbf{Q}_p(\varphi_{\mathfrak{p}^*}) & \text{if } p \text{ splits} \\ K_{\mathfrak{p}}(\varphi_{\mathfrak{p}}) & \text{if } p \text{ does not split} \end{cases} \tag{6.0.1}$$

We will consider $p$ splits case and $p$ nonsplit case separately. It does not matter whether $p$ is a good reduction prime for $E$ or not.

## 6.1 $p$ splits

Define

$$\begin{aligned} Sel(p) &= Ker\left(H^1(\mathbf{Q}, V_p/T_p) \to \prod_{l<\infty} \frac{H^1(\mathbf{Q}_l, T_p \otimes \mathbf{Q}_p/\mathbf{Z}_p)}{H^1_f(\mathbf{Q}_l, T_p) \otimes \mathbf{Q}_p/\mathbf{Z}_p}\right) \\ &= Ker\left(H^1(\mathbf{Q}, V_p/T_p) \to \prod_{l\neq p} H^1(\mathbf{Q}_l, V_p/T_p) \times \frac{H^1(\mathbf{Q}_p, V_p/T_p)}{H^1_f(\mathbf{Q}_p, T_p) \otimes \mathbf{Q}_p/\mathbf{Z}_p}\right) \\ &\cong Ker\left(H^1(K, V_p/T_p) \to \prod_{\nu\neq \mathfrak{p}} H^1(K_\nu, \mathbf{Q}_p/\mathbf{Z}_p(\varphi_{\mathfrak{p}})) \times \frac{H^1(K_{\mathfrak{p}}, \mathbf{Q}_p/\mathbf{Z}_p(\varphi_{\mathfrak{p}}))}{H^1_f(K_{\mathfrak{p}}, \mathbf{Z}_p(\varphi_{\mathfrak{p}})) \otimes \mathbf{Q}_p/\mathbf{Z}_p}\right) \end{aligned} \tag{6.1.1}$$

Here is the reason for the last equality:

- First we have restriction map $H^1(\mathbf{Q}, V_p/T_p) \cong H^1(K, V_p/T_p)^{Gal(K/\mathbf{Q})}$ and since

$$H^1(K, V_p/T_p) = H^1(K, \mathbf{Q}_p/\mathbf{Z}_p(\varphi_{\mathfrak{p}})) \oplus H^1(K, \mathbf{Q}_p/\mathbf{Z}_p(\varphi_{\mathfrak{p}^*}))$$

  where $\tau \in Gal(K/\mathbf{Q})$ acts by swapping the two factors, we have

$$H^1(K, \mathbf{Q}_p/\mathbf{Z}_p(\varphi_{\mathfrak{p}})) \cong H^1(K, V_p/T_p)^{Gal(K/\mathbf{Q})}$$

under which $\theta \mapsto \theta \oplus \tau\theta$. The same argument applies to local Galois cohomology groups at $\nu \mid l$ and $l$ is nonsplit:

$$H^1(K_\nu, V_p/T_p)^{Gal(K_\nu/\mathbf{Q}_l)} \cong H^1(K_\nu, \mathbf{Q}_p/\mathbf{Z}_p(\varphi_\mathfrak{p}))$$

Hence in such case we have the following commutative diagram:

$$\begin{array}{ccccc}
H^1(\mathbf{Q}, V_p/T_p) & \xrightarrow[\text{res}]{\cong} & H^1(K, V_p/T_p)^{Gal(K/\mathbf{Q})} & \xrightarrow{\cong} & H^1(K, \mathbf{Q}_p/\mathbf{Z}_p(\varphi_\mathfrak{p})) \\
\downarrow{\scriptstyle \text{res}} & & \downarrow{\scriptstyle \text{res}} & & \downarrow{\scriptstyle \text{res}} \\
H^1(\mathbf{Q}_l, V_p/T_p) & \xrightarrow[\text{res}]{\cong} & H^1(K_\nu, V_p/T_p)^{Gal(K_\nu/\mathbf{Q}_l)} & \xrightarrow{\cong} & H^1(K_\nu, \mathbf{Q}_p/\mathbf{Z}_p(\varphi_\mathfrak{p}))
\end{array}$$

$$\begin{aligned}
&Ker(H^1(\mathbf{Q}, V_p/T_p) \to H^1(\mathbf{Q}_l, V_p/T_p)) \\
&\quad \cong Ker(H^1(K, \mathbf{Q}_p/\mathbf{Z}_p(\varphi_\mathfrak{p})) \to H^1(K_\nu, \mathbf{Q}_p/\mathbf{Z}_p(\varphi_\mathfrak{p})))
\end{aligned} \tag{6.1.2}$$

- In the case $l \neq p$ and $l$ splits in $K/\mathbf{Q}$, say $l = \nu\nu^*$, we have the following commutative diagram:

$$\begin{array}{ccc}
H^1(\mathbf{Q}, V_p/T_p) & \xrightarrow{\cong} & H^1(K, \mathbf{Q}_p/\mathbf{Z}_p(\varphi_\mathfrak{p})) \\
\downarrow{\scriptstyle \text{res}} & & \downarrow{\scriptstyle (\text{res},\text{res})} \\
H^1(\mathbf{Q}_l, V_p/T_p) & \xrightarrow{\cong} & H^1(K_\nu, \mathbf{Q}_p/\mathbf{Z}_p(\varphi_\mathfrak{p})) \oplus H^1(K_{\nu^*}, \mathbf{Q}_p/\mathbf{Z}_p(\varphi_\mathfrak{p}))
\end{array}$$

where in the bottom row we use identification:

$$H^1(K_\nu, \mathbf{Q}_p/\mathbf{Z}_p(\varphi_{\mathfrak{p}^*})) \cong H^1(K_{\nu^*}, \mathbf{Q}_p/\mathbf{Z}_p(\varphi_\mathfrak{p}))$$

From this diagram it follows that

$$\begin{aligned}
&Ker(H^1(\mathbf{Q}, V_p/T_p) \to H^1(\mathbf{Q}_l, V_p/T_p)) \\
&\quad \cong \cap_{\omega=\nu,\nu^*} Ker(H^1(K, \mathbf{Q}_p/\mathbf{Z}_p(\varphi_\mathfrak{p})) \to H^1(K_\omega, \mathbf{Q}_p/\mathbf{Z}_p(\varphi_\mathfrak{p})))
\end{aligned} \tag{6.1.3}$$

- similarly we have commutative diagram:

$$\begin{array}{ccc}
H^1(\mathbf{Q}, V_p/T_p) & \xrightarrow{\cong} & H^1(K, \mathbf{Q}_p/\mathbf{Z}_p(\varphi_\mathfrak{p})) \\
\downarrow{\scriptstyle \text{res}} & & \downarrow{\scriptstyle (\text{res},\text{res})} \\
\frac{H^1(\mathbf{Q}_p, V_p/T_p)}{H^1_f(\mathbf{Q}_p, T_p) \otimes \mathbf{Q}_p/\mathbf{Z}_p} & \xrightarrow{\cong} & \frac{H^1(K_\mathfrak{p}, \mathbf{Q}_p/\mathbf{Z}_p(\varphi_\mathfrak{p}))}{H^1_f(K_\mathfrak{p}, \mathbf{Z}_p(\varphi_\mathfrak{p})) \otimes \mathbf{Q}_p/\mathbf{Z}_p} \oplus H^1(K^*_\mathfrak{p}, \mathbf{Q}_p/\mathbf{Z}_p(\varphi_\mathfrak{p}))
\end{array}$$

where the bottom isomorphism is because $H^1_f(\mathbf{Q}_p, \mathbf{Z}_p(\varphi_{\mathfrak{p}^*}))$ is finite torsion group.

We see from these that we also have the following:

$$S' \stackrel{\text{def}}{=} Ker\left(H^1(\mathbf{Q}, V_p/T_p) \to \prod_{l \neq p} \frac{H^1(\mathbf{Q}_l, T_p \otimes \mathbf{Q}_p/\mathbf{Z}_p)}{H^1_f(\mathbf{Q}_l, T_p) \otimes \mathbf{Q}_p/\mathbf{Z}_p} \times \frac{H^1(\mathbf{Q}_p, V_p/T_p)}{H^1(\mathbf{Q}_p, \mathbf{Q}_p/\mathbf{Z}_p(\varphi_{\mathfrak{p}}))}\right)$$
$$\cong Ker\left(H^1(K, \mathbf{Q}_p/\mathbf{Z}_p(\varphi_{\mathfrak{p}})) \to \prod_{\nu \neq \mathfrak{p}} H^1(K_\nu, \mathbf{Q}_p/\mathbf{Z}_p(\varphi_{\mathfrak{p}}))\right) \tag{6.1.4}$$

We have by definition :

$$0 \longrightarrow Sel(p) \longrightarrow S'(p) \longrightarrow \frac{A'_p}{A_p} \tag{6.1.5}$$

where $A_p = H^1_f(\mathbf{Q}_p, \mathbf{Z}_p(\varphi_{\mathfrak{p}})) \otimes \mathbf{Q}_p/\mathbf{Z}_p \subseteq A'_p = H^1(\mathbf{Q}_p, \mathbf{Q}_p/\mathbf{Z}_p(\varphi_{\mathfrak{p}}))$. We are going to relate $S'(p)$ to familiar objects in Iwasawa theory, and find the exact difference between $Sel(p)$ and $S'(p)$.

**Lemma 23.** *$S'(p) \cong Hom(\mathfrak{X}_\infty, \mathbf{Q}_p/\mathbf{Z}_p(\varphi_{\mathfrak{p}}))^{G_\infty}$ where $\mathfrak{X}_\infty$ is the Galois group over $K_\infty$ of the maximal $p$-abelian unramified outside places above $\mathfrak{p}$ extension of $K_\infty$.*

*Proof.* First note $H^1(G_\infty, \mathbf{Q}_p/\mathbf{Z}_p(\varphi_{\mathfrak{p}})) = H^2(G_\infty, \mathbf{Q}_p/\mathbf{Z}_p(\varphi_{\mathfrak{p}})) = 0$ as $p$-cohomological dimension of $G_\infty$ is 1, and that $1-\sigma$ acts invertibly on $\mathbf{Q}_p(\varphi_{\mathfrak{p}})$ where $\sigma$ is a generator of $G_\infty$. So restriction induces an isomorphism:

$$H^1(K, \mathbf{Q}_p/\mathbf{Z}_p(\varphi_{\mathfrak{p}})) \cong H^1(K_\infty, \mathbf{Q}_p/\mathbf{Z}_p(\varphi_{\mathfrak{p}}))^{G_\infty}$$

It is well known from criterion of Ogg-Neron-Shafarevich that $E$ has good reduction at all places of $K(E_{\mathfrak{p}})$ not lying above $\mathfrak{p}$. For place $\lambda \neq \mathfrak{p}$ of $K$, let $I_\lambda$ be the inertia subgroup of $Gal(\overline{K_\lambda}/K_\lambda)$ and let $I'_\lambda = I_\lambda \cap Gal(\bar{\mathbf{Q}}/K_\infty)$.Then $I'_\lambda$ is of finite index in $I_\lambda$ and $[I_\lambda : I_\lambda']$ is prime to $p$ since $\#((Gal(K(E_\mathfrak{p})/K))$ is prime to $p$. We have the following commutative diagram:

$$\begin{array}{ccccc}
H^1(K, \mathbf{Q}_p/\mathbf{Z}_p(\varphi_\mathfrak{p})) & \longrightarrow & H^1(K_\lambda, \mathbf{Q}_p/\mathbf{Z}_p(\varphi_\mathfrak{p})) & \xrightarrow{\rho_1} & H^1(I_\lambda, \mathbf{Q}_p/\mathbf{Z}_p(\varphi_\mathfrak{p})) \\
\downarrow{\scriptstyle \mathrm{res}} & & & & \downarrow{\scriptstyle \rho_2} \\
H^1(K_\infty, \mathbf{Q}_p/\mathbf{Z}_p(\varphi_\mathfrak{p}))^{G_\infty} & \xrightarrow{\cong} & Hom(X_\infty, \mathbf{Q}_p/\mathbf{Z}_p(\varphi_\mathfrak{p}))^{G_\infty} & \longrightarrow & H^1(I_\lambda', \mathbf{Q}_p/\mathbf{Z}_p(\varphi_\mathfrak{p})) \\
 & & & & \| \\
 & & & & Hom(I_\lambda', \mathbf{Q}_p/\mathbf{Z}_p(\varphi_\mathfrak{p}))
\end{array}$$

where $X_\infty$ is the Galois group over $K_\infty$ of the maximal $p$-abelian extension of $K_\infty$. The map $\rho_1$ is injective, since

$$Ker(\rho_1) \cong H^1(K_\lambda^{\mathrm{unr}}/K_\lambda, \mathbf{Q}_p/\mathbf{Z}_p(\varphi_\mathfrak{p})^{I_\lambda}) = 0$$

This is trival when $(\mathfrak{f}, \lambda) = 1$ as then $\mathbf{Q}_p/\mathbf{Z}_p(\varphi_\mathfrak{p})^{I_\lambda} = \mathbf{Q}_p/\mathbf{Z}_p(\varphi_\mathfrak{p})$ is divisible. When $(\mathfrak{f}, \lambda) \neq 1$, then since $\varphi_\mathfrak{p}$ factors through $K(\mathfrak{f}\mathfrak{p}^\infty)$ , so $\varphi_\mathfrak{p}(I_\lambda')$ is a nontrival finite subgroup of $\mathcal{O}_\mathfrak{p}^*$. Since $1 + \mathfrak{p}\mathcal{O}_\mathfrak{p}$ has no torsion, we see $\varphi_\mathfrak{p}(I_\lambda') \nsubseteq 1 + \mathfrak{p}\mathcal{O}_\mathfrak{p}$ which inplies again $Ker(\rho_1)$ is trivial. Note $\rho_2$ is also injective, as $\#(I_\lambda/I_\lambda')$ is prime to $p$. From this diagram, we see that

$$S'(p) \cong Hom(\mathfrak{X}_\infty, \mathbf{Q}_p/\mathbf{Z}_p(\varphi_\mathfrak{p}))^{G_\infty} \cong Hom(\Xi_\infty/\overline{\mathfrak{C}_\infty}, \mathbf{Q}_p/\mathbf{Z}_p(\varphi_\mathfrak{p}))^{G_\infty} \qquad (6.1.6)$$

where the second isomorphism comes from the class field exact sequence:

$$0 \longrightarrow \overline{\mathfrak{A}_\infty}/\overline{\mathfrak{C}_\infty} \longrightarrow \Xi_\infty/\overline{\mathfrak{C}_\infty} \longrightarrow \mathfrak{X}_\infty \longrightarrow A_\infty \longrightarrow 0 \qquad (6.1.7)$$

where $A_\infty = \lim_{F \subset K_\infty} Cl(F)$. Using Main conjecture in Iwasawa theory for $K$ as proved in [10] 4.4, applying the 1 variable version when $j = 0$ and 2 variable version when $j > 0$, we get the second isomorphism and we know from this $S'(p)$ is a finite torsion group as $\Xi_\infty/\overline{\mathfrak{C}_\infty}$ is a torsion $\Lambda$ module in this case ($p$ splits). Hence by 6.1.5, $Sel(p)$ is also finite and so $Sel(p) = \bigsqcup(p)$. By the way this also gives us $\mathbf{A}(\mathbb{Q})$ is finite torsion in our case, since $\mathbf{A}(\mathbb{Q}) \otimes \mathbf{Z}_p = H^1_{f,spec(\mathbf{Z})}(\mathbf{Q}, T_p)$ which is finite as

$$H^1_{f,spec(\mathbf{Z})}(\mathbf{Q}, T_p) \otimes \mathbf{Q}_p/\mathbf{Z}_p \hookrightarrow Sel(p) \qquad (6.1.8)$$

□

*Remark 17.* • Lemma 23 and its proof, except the second isomorphism in 6.1.6 holds word for word in the case $p$ is nonsplit.

- One can prove in the most general situation that $\bigsqcup(p)$ is finite. The hard part is to prove it is 0 for almost all $p$ and to link the order to special values of $L$ functions. Lemma 23 enables us to do so for $S^{'}(p)$ by the well known link between $\Xi_\infty/\overline{\mathfrak{C}_\infty}$ and $p$ adic $L$ functions of $K$. This had been done in both [5] and [6]. In our approach we do not do this, rather we note $\#(S^{'}(p))$ will cancel another term in Tamagawa number formula 2.2.6. Now we just need to pin down exactly the difference between $Sel(p)$ and $S^{'}(p)$.

For this purpose we need a more refined version of 6.1.5 as follows. First a general remark: Let $V_p$ be a $p$-adic de Rham representation of $Gal(\bar{\mathbf{Q}}/\mathbf{Q})$. Then

$$Sel(p) \cong Ker\left(H^1(G_{S_p}, V_p/T_p) \to \prod_{l\in S_p} \frac{H^1(\mathbf{Q}_l, V_p/T_p)}{H^1_f(\mathbf{Q}_l, T_p)\otimes \mathbf{Q}_p/\mathbf{Z}_p}\right)$$

where $S_p = \{l : V_p \text{ is ramified at } l\ \}\cup\{p\}$ and $G_{S_p}$ is the absolute Galois group of the maximal $p$ extension of $\mathbf{Q}$ unramified outside places in $S_p$. This follows from the fact that

$$\frac{H^1(\mathbf{Q}_l, V_p/T_p)}{H^1_f(\mathbf{Q}_l, T_p)\otimes \mathbf{Q}_p/\mathbf{Z}_p} \subseteq H^1(\mathbf{Q}_l^{\text{unr}}, V_p/T_p) = Hom(I_l, V_p/T_p),\ \forall l \notin S_p$$

by definition. Now back to our case $V_p$ as given in the beginning of this section and $p$ splits. Consider the commutative diagram:

$$\begin{array}{ccccccccc}
0 \to 0 \to & H^1(G_{S_p}, V_p/T_p) & \xrightarrow{id} & H^1(G_{S_p}, V_p/T_p) & \to 0 \\
& \downarrow \lambda & & \downarrow \text{res} & \\
0 \to \frac{A_p^{'}}{A_p} \to & \bigoplus_{l\in S_p} \frac{H^1(\mathbf{Q}_l,V_p/T_p)}{H^1_f(\mathbf{Q}_l,T_p)\otimes\mathbf{Q}_p/\mathbf{Z}_p} & \to & \bigoplus_{l\in S_p, l\neq p} \frac{H^1(\mathbf{Q}_l,V_p/T_p)}{H^1_f(\mathbf{Q}_l,T_p)\otimes\mathbf{Q}_p/\mathbf{Z}_p} \oplus \frac{H^1(\mathbf{Q}_p,V_p/T_p)}{A_p^{'}} & \to 0
\end{array}$$

It follows from this by snake lemma that

$$0 \longrightarrow Sel(p) \longrightarrow S^{'}(p) \longrightarrow \frac{A_p^{'}}{A_p} \longrightarrow coker(\lambda) \longrightarrow coker(res) \longrightarrow 0 \qquad (6.1.9)$$

We also have the following variant of the Poitou-Tate exact sequence:

$$\begin{array}{ccccccc} 0 & \longrightarrow & Sel(p) & \longrightarrow & H^1(G_{S_p}, V_p/T_p) & \xrightarrow{\lambda} & \oplus_{l\in S_p} \frac{H^1(\mathbf{Q}_l, V_p/T_p)}{H^1_f(\mathbf{Q}_l, T_p)\otimes \mathbf{Q}_p/\mathbf{Z}_p} \\ & & & & & & \downarrow \cong \\ 0 & = & H^2(G_{S_p}, V_p/T_p) & \longleftarrow & H^1_{f,spec(\mathbf{Z})}(\mathbf{Q}, T_p^*(1))^\vee & \xleftarrow[\rho^*]{} & \oplus_{l\in S_p} H^1(\mathbf{Q}_l, T_p^*(1))^\vee \end{array}$$

where $(-)^\vee$ means the Pontrayagin dual. If we apply our argument of finitness of $Sel(p)$ to similarly defined Selmer group for $T_p^*(1)$ everything still works, in particular we will have

$$S^{'}(T_p^*(1)) \cong Hom(\mathfrak{X}_\infty, \mathbf{Q}_p/\mathbf{Z}_p(\varphi_{\mathfrak{p}}^{-1}\chi_{\text{cyclo}}))^{G_\infty}$$

except in this case $\varphi_{\mathfrak{p}}^{-1}\chi_{\text{cyclo}} = \psi_{\mathfrak{p}}^{-k-j}\chi_{\text{cyclo}}^{1+j}$ so we have to use 2 variable version of Main conjecture even if $j = 0$. We conclude $Sel(T_p^*(1))$ is still finite and hence $H^1_{f,spec(\mathbf{Z})}(\mathbf{Q}, T_p^*(1))$ is also finite torsion whose order is $\#((V_p^*(1)/T_p^*(1))^{G_\mathbf{Q}})$. Since $\rho^*$ is onto, it induces an isomorphism of $coker(\lambda)$ onto $H^1_{f,spec(\mathbf{Z})}(\mathbf{Q}, T_p^*(1))^\vee$. In particular, $\#(coker(\lambda)) = \#((V_p^*(1)/T_p^*(1))^{G_\mathbf{Q}})$. Now $(A_p^{'}/A_p)^\vee \cong H^1_f(\mathbf{Q}_p, \mathbf{Z}_p(\varphi_{\mathfrak{p}}^{-1}\chi))$ is a finite torsion group which has order

$$\delta_p \overset{\text{def}}{=} \#(H^0(\mathbf{Q}_p, \mathbf{Q}_p/\mathbf{Z}_p(\psi_{\mathfrak{p}}^{-k-j}\chi^{1+j}))) \tag{6.1.10}$$

In order to prove $coker(res) = 0$, by 6.1.9, we just need to prove that any representative of $coker(\lambda)$ can be changed by an element of $A_p^{'}/A_p$ so as to lie in $Im(\lambda)$ which by the Poitou-Tate exact sequence, it is enough to show that $\rho^*$ maps $A_p^{'}/A_p$ onto $H^1_{f,spec(\mathbf{Z})}(\mathbf{Q}, T_p^*(1))^\vee$. But we see that

$$\rho : H^1_{f,spec(\mathbf{Z})}(\mathbf{Q}, T_p^*(1)) \hookrightarrow (A_p^{'}/A_p)^\vee$$

It follows that $\rho^*$ maps $A_p^{'}/A_p \subset \oplus_{l\in S_p} H^1(\mathbf{Q}_l, V_p/T_p)/(H^1_f(\mathbf{Q}_l, T_p)\otimes \mathbf{Q}_p/\mathbf{Z}_p)$ already onto $H^1_{f,spec(\mathbf{Z})}(\mathbf{Q}, T_p^*(1))^\vee$, hence the map "res" is actually onto. Hence by 6.1.9 we have

$$\#(Hom(\Xi_\infty/\overline{\mathfrak{C}_\infty}, \mathbf{Q}_p/\mathbf{Z}_p(\varphi_{\mathfrak{p}}))^{G_\infty})^{-1}\delta_p = \#((V_p^*(1)/T_p^*(1))^{G_\mathbf{Q}})\#(\bigsqcup(p))^{-1} \tag{6.1.11}$$

## 6.2 $p$ is nonsplit

Here the basic idea is still try to apply Iwasawa Main conjecture for $p$ non-split case ([10] 4.4 (ii)). The proof will be somewhat different from $p$ splits case and it is essentially due to Rubin who dealt with the case $k = 1, j = 0$. In this section, $V_p = K_{\mathfrak{p}}(\varphi_{\mathfrak{p}})$. First consider

$$S(p) \overset{\text{def}}{=} Ker\left(H^1(K, V_p/T_p) \to \prod_{\nu \subset \mathcal{O}_K} \frac{H^1(K_\nu, V_P/T_p)}{H^1_f(K_\nu, T_p) \otimes K_{\mathfrak{p}}/\mathcal{O}_{\mathfrak{p}})}\right)$$

**Lemma 24.** $Sel(p) = S(p)^{Gal(K/\mathbf{Q})}$

*Proof.* This follows from:

- if $l$ splits and $l \neq p$, say $l = \nu\nu^*$, the following diagram commutes:

$$\begin{array}{ccc} H^1(\mathbf{Q}, V_p/T_p) & \xrightarrow[\text{res}]{\cong} & H^1(K, V_p/T_p)^{Gal(K/\mathbf{Q})} \\ \Big\downarrow \text{res} & & \Big\downarrow \text{(res,res)} \\ H^1(\mathbf{Q}_l, V_p/T_p) & \xrightarrow{\rho} & H^1(K_\nu, K_{\mathfrak{p}}/\mathcal{O}_{\mathfrak{p}}(\varphi_{\mathfrak{p}})) \oplus H^1(K_{\nu^*}, K_{\mathfrak{p}}/\mathcal{O}_{\mathfrak{p}}(\varphi_{\mathfrak{p}})) \end{array}$$

  where $\rho(f) = (f, \tau(f))$ is injective.

- by definition of $H^1_f$, the restriction map gives isomorphism:

$$H^1(\mathbf{Q}_p, V_p/T_p) \cong H^1(K_{\mathfrak{p}}, V_p/T_p)^{Gal(K_{\mathfrak{p}}/\mathbf{Q}_p)}$$

  which induces:

$$H^1_f(\mathbf{Q}_p, T_p) \otimes \mathbf{Q}_p/\mathbf{Z}_p \cong H^1_f(K_{\mathfrak{p}}, T_p) \otimes K_{\mathfrak{p}}/\mathcal{O}_{\mathfrak{p}})^{Gal(K_{\mathfrak{p}}/\mathcal{O}_{\mathfrak{p}})}$$

  and the following diagram commutes:

$$\begin{array}{ccc} H^1(\mathbf{Q}, V_p/T_p) & \xrightarrow[\text{res}]{\cong} & H^1(K, V_p/T_p)^{Gal(K/\mathbf{Q})} \\ \Big\downarrow \text{res} & & \Big\downarrow \text{res} \\ H^1(\mathbf{Q}_p, V_p/T_p)/(H^1_f(\mathbf{Q}_p, T_p) \otimes \mathbf{Q}_p/\mathbf{Z}_p) & \xrightarrow[\text{res}]{\cong} & (\frac{H^1(K_{\mathfrak{p}}, V_p/T_p)}{H^1_f(K_{\mathfrak{p}}, T_p) \otimes K_{\mathfrak{p}}/\mathcal{O}_{\mathfrak{p}}})^{Gal(K_{\mathfrak{p}}/\mathcal{O}_{\mathfrak{p}})} \end{array}$$

- if $l$ is nonsplit, $\nu \mid l$, then

$$\begin{array}{ccc} H^1(\mathbf{Q}, V_p/T_p) & \xrightarrow[\text{res}]{\cong} & H^1(K, V_p/T_p)^{Gal(K/\mathbf{Q})} \\ \downarrow{\scriptstyle \text{res}} & & \downarrow \\ H^1(\mathbf{Q}_l, V_p/T_p) & \xrightarrow[\text{res}]{\cong} & H^1(K_\nu, K_{\mathfrak{p}}(\varphi_{\mathfrak{p}}))^{Gal(K_\nu/\mathbf{Q}_l)} \end{array}$$

□

Define

$$S^{'}(p) = Ker\left(H^1(K, V_p/T_p) \to \prod_{\nu \neq \mathfrak{p}} \frac{H^1(K_\nu, V_P/T_p)}{H^1_f(K_\nu, T_p) \otimes K_{\mathfrak{p}}/\mathcal{O}_{\mathfrak{p}})}\right)$$

As remarked in last section, restrction gives:

$$H^1(K, V_p/T_p) \cong H^1(K_\infty, V_p/T_p)^{G_\infty}$$

which induces $S^{'}(p) \cong Hom(\mathfrak{X}_\infty, K_{\mathfrak{p}}/\mathcal{O}_{\mathfrak{p}}(\varphi_{\mathfrak{p}}))^{G_\infty}$. It is clear that

$$S(p) = Ker\left(S^{'}(p) \to \frac{H^1(K_{\mathfrak{p}}, V_p/T_p)}{H^1_f(K_{\mathfrak{p}}, T_p) \otimes K_{\mathfrak{p}}/\mathcal{O}_{\mathfrak{p}}}\right) \qquad (6.2.1)$$

Now unlike in the case $p$ splits, $\mathfrak{X}_\infty$ is no longer a torsion $\Lambda$ module. Yet we can still apply Main conjecture as follows.

Write $G_\infty = \Gamma \times \Delta$, where $\Delta = Gal(K(E_{\mathfrak{p}})/K)$ is a finite group of order prime to $p$. Let $\chi = \varphi \mid_\Delta$. We have seen:

$$\begin{aligned} H^1(K_{\mathfrak{p}}, V_p/T_p) &\cong Hom(\Xi_\infty, K_{\mathfrak{p}}/\mathcal{O}_{\mathfrak{p}}(\varphi_{\mathfrak{p}}))^{G_\infty} \\ &\cong Hom(\Xi^\chi_\infty, K_{\mathfrak{p}}/\mathcal{O}_{\mathfrak{p}}(\varphi_{\mathfrak{p}}))^{\Gamma} \end{aligned}$$

Denote $\eta : H^1_f(K_{\mathfrak{p}}, T_p) \otimes K_{\mathfrak{p}}/\mathcal{O}_{\mathfrak{p}} \to Hom(\Xi^\chi_\infty, \mathcal{O}_{\mathfrak{p}}(\varphi_{\mathfrak{p}}))^{\Gamma}$. Recall $\Xi^\chi_\infty = \Xi^1_\infty \oplus \Xi^2_\infty$ where $\Xi^1_\infty$ and $\Xi^2_\infty$ are free rank 1 $\Lambda$ module with $\delta(\Xi^2_\infty) = 0$ for a generator $\delta$ of the image of $\eta$. By 6.2.1 we see that

$$S(p) = \{f \in Hom(\Xi^\chi_\infty, K_{\mathfrak{p}}/\mathcal{O}_{\mathfrak{p}}(\varphi_{\mathfrak{p}}))^{\Gamma} : f|_{\Xi_\infty} \in Image(\eta)\} \qquad (6.2.2)$$

Now since we have calculated in previous sections that $\delta(u) \neq 0$ for $u \in \overline{\mathfrak{C}_\infty}$ our chosen generator. Hence $\hat{\Xi} \stackrel{\text{def}}{=} \Xi_\infty^\chi / (\Xi^2, \overline{\mathfrak{C}_\infty})$ is a torsion $\Lambda$ module. Let $\hat{\mathfrak{X}}$ be the quotient of $\mathfrak{X}_\infty^\chi$ by the image of $\Xi_\infty^2$ under the class field theory map:

$$\Xi_\infty^2 \subseteq \Xi_\infty^\chi \to \Xi_\infty^\chi / \overline{\mathfrak{A}_\infty}^\chi \hookrightarrow \mathfrak{X}_\infty^\chi$$

then by 6.2.2, we have

$$S(p) = Hom(\hat{\mathfrak{X}}, K_{\mathfrak{p}}/\mathcal{O}_{\mathfrak{p}}(\varphi_{\mathfrak{p}}))^\Gamma$$

Now for out fixed $k, j$, assume nonsplit prime $p$ satisfies that $\chi = \psi_p^{k+j} \chi_{\text{cyclo}}^{-j}|_\Delta$ is nontrival, (note this is the only place in the paper where we use this assumption on $p$. One should be able to remove this condition in the nonsplit case of Iwasawa Main conjecture as proved by Rubin [10]) we can then apply [10] theorem 10.6 which gives via 6.1.7 (see [10] 11.13) $Char(\hat{\mathfrak{X}}) = Char(\hat{\Xi})$. By [10] theorem 5.3 (v) and theorem 11.12(ii), $\hat{\mathfrak{X}}$ and $\hat{\Xi}$ has no finite torsion submodule, hence

$$S(p) = Hom(\Xi_\infty^\chi / (\Xi^2, \overline{\mathfrak{C}_\infty}), K_{\mathfrak{p}}/\mathcal{O}_{\mathfrak{p}}(\varphi_{\mathfrak{p}}))^\Gamma \tag{6.2.3}$$

## 6.3 Putting everything together

Before we finally deduce our main theorem from all these computations of $exp^*$, let us remark on the constant $\mathbf{N}\boldsymbol{\alpha} - (\psi^k \bar{\psi}^{-j})(\boldsymbol{\alpha})$ that always comes up in our formulas: $\forall \alpha, (\alpha, \mathfrak{f}\mathfrak{p}) = 1$,

$$\mathbf{N}\boldsymbol{\alpha} - (\psi^k \bar{\psi}^{-j})(\boldsymbol{\alpha}) = -\psi^k \bar{\psi}^{-j}(\boldsymbol{\alpha})(1 - \psi^{-k} \bar{\psi}^j \chi_{\text{cyclo}}(\boldsymbol{\alpha}))$$

Since $\psi(\boldsymbol{\alpha})$ is a generator of $\boldsymbol{\alpha}$, hence $\psi^k \bar{\psi}^{-j}(\boldsymbol{\alpha})$ is a $\mathfrak{p}$ unit. It is clear we can choose an integral ideal $\boldsymbol{\alpha}$ such that

$$\left|\mathbf{N}\boldsymbol{\alpha} - (\psi^k \bar{\psi}^{-j})(\boldsymbol{\alpha})\right|_{\mathfrak{p}}^{-1} = \#(H^0(K_{\mathfrak{p}}, K_{\mathfrak{p}}/\mathcal{O}_{\mathfrak{p}}(\psi^{-k} \bar{\psi}^j \chi_{\text{cyclo}}))) \tag{6.3.1}$$

Now to finish the proof of the main theorem, or to verify both sides of 2.2.6 have the same $p$ valuation, for $p \neq 2, 3$ and $p \nmid d_K$, we just need to observe:

1. when $p$ splits, this follows from 3.3.2 4.2.4, 4.2.7, 5.1.3, 5.1.21, 5.1.28 and 6.1.11. Also notice the constants $f, g$ or $\alpha, \alpha^{'}$ in these equations all have $\mathfrak{p}$ valuation 0. Note

$$\frac{\mu_p(H^1_f(\mathbf{Q}_p, T_p))}{|P_p(V,1)|_p^{-1}} = |a|_p^{-1} \#(H^0(\mathbf{Q}_p, \mathbf{Q}_p/\mathbf{Z}_p(\varphi_{\mathfrak{p}^*})))$$

and $\frac{L_{\mathfrak{f}\mathfrak{p}}(\bar{\psi}^{k+j},k)}{P_p(V,1)} \times \#(H^0(\mathbf{Q}_p, \mathbf{Q}_p/\mathbf{Z}_p(\varphi_{\mathfrak{p}^*})))$ exactly takes care of $L(\bar{\psi}^{k+j}, k)$ and by above remark, $\left|\mathbf{N}\boldsymbol{\alpha} - (\psi^k\bar{\psi}^{-j})(\boldsymbol{\alpha})\right|_{\mathfrak{p}}^{-1}$ cancels the term $\delta_p$ that comes in $\#(V_p^*(1)/T_p^*(1))^{G_{\mathbf{Q}}})\#(\bigsqcup(p))^{-1})$.

2. when $p$ is inert, $\mathfrak{p} \mid p$, we have proved so far: under $exp_{K_{\mathfrak{p}}} : DR(K_{\mathfrak{p}}, K_{\mathfrak{p}}(\varphi_{\mathfrak{p}})) \to H^1_f(K_{\mathfrak{p}}, K_{\mathfrak{p}}(\varphi_{\mathfrak{p}}))$ we have

$$exp(L/L^0 \otimes \mathcal{O}_{\mathfrak{p}}) = aH^1_f(K_{\mathfrak{p}}, \mathcal{O}_{\mathfrak{p}}(\varphi_{\mathfrak{p}}))$$

and by 4.2.4 4.2.7 , 5.2.1 and 3.1.7, we see

$$\begin{aligned} |a|_{\mathfrak{p}}^{-1} &= |b|_{\mathfrak{p}}^{-1} \\ &= \left| \left(\frac{\sqrt{d_K}}{2\pi}\right)^{-j} \Omega_\infty^{-k-j} (\mathbf{N}\boldsymbol{\alpha} - \psi^k\bar{\psi}^{-j}(\boldsymbol{\alpha})) L_{\mathfrak{f}\mathfrak{p}}(\bar{\psi}^{k+j}, k) \right|_{\mathfrak{p}}^{-1} \qquad (6.3.2) \\ &\times \#(Hom_{G_\infty}(\Xi_\infty/(\Xi_\infty^2, \overline{\mathfrak{C}_\infty}), K_{\mathfrak{p}}/\mathcal{O}_{\mathfrak{p}}(\varphi_{\mathfrak{p}}))) \end{aligned}$$

Here again we do not care $\mathfrak{p}$ units $\alpha, d, f$ (or $g$).We have commutative diagram:

$$\begin{array}{ccc} DR(K_{\mathfrak{p}}, K_{\mathfrak{p}}(\varphi_{\mathfrak{p}})) & \xrightarrow{exp_{K_{\mathfrak{p}}}} & H^1_f(K_{\mathfrak{p}}, K_{\mathfrak{p}}(\varphi_{\mathfrak{p}})) \\ \text{inclu}\uparrow & & \uparrow\text{res} \\ DR(\mathbf{Q}_p, K_{\mathfrak{p}}(\varphi_{\mathfrak{p}})) & \xrightarrow{exp} & H^1_f(\mathbf{Q}_p, K_{\mathfrak{p}}(\varphi_{\mathfrak{p}})) \end{array}$$

By this diagram we view the second row as the $Gal(K_{\mathfrak{p}}/\mathbf{Q}_p) \cong Gal(K/\mathbf{Q})$ fixed part of the first row. Note

$$\left(\frac{\sqrt{d_K}}{2\pi}\right)^{-j} \Omega_\infty^{-k-j} L_{\mathfrak{f}\mathfrak{p}}(\bar{\psi}^{k+j}, k) \in \mathbf{Q}$$

hence we get

$$\mu_{\varpi,p}(H^1_f(K_{\mathfrak{p}},\mathcal{O}_{\mathfrak{p}}(\varphi_{\mathfrak{p}}))) = \left| \left(\frac{\sqrt{d_K}}{2\pi}\right)^{-j} \Omega_\infty^{-k-j}(\mathbf{N}\boldsymbol{\alpha} - \psi^k\bar{\psi}^{-j}(\boldsymbol{\alpha}))L_{\mathfrak{f}\mathfrak{p}}(\bar{\psi}^{k+j},k)\right|_p^{-1} \\ \times \#(Hom_{G_\infty}(\Xi_\infty/(\Xi^2_\infty,\overline{\mathfrak{C}_\infty}),K_{\mathfrak{p}}/\mathcal{O}_{\mathfrak{p}}(\varphi_{\mathfrak{p}}))^{Gal(K/\mathbf{Q})}) \tag{6.3.3}$$

By 6.2.3 and lemma 24 , we get

$$\#(\bigsqcup(p)) = \#(Sel(p)) = \#(Hom_{G_\infty}(\Xi_\infty/(\Xi^2_\infty,\overline{\mathfrak{C}_\infty}),K_{\mathfrak{p}}/\mathcal{O}_{\mathfrak{p}}(\varphi_{\mathfrak{p}}))^{Gal(K/\mathbf{Q})})$$

Hence we get both sides of 2.2.6 have the same $p$ valuation when $p$ is inert.

3. when $p$ ramifies, all our computations on $exp^*(u)$ and $\#(\bigsqcup(p))$ still holds. the only discrepency is the term $\sqrt{d_K}$ which seems to have different power on both sides of 2.2.6.